\documentclass[11pt,a4paper]{article}
\usepackage[T1]{fontenc}
\usepackage{lmodern}
\usepackage[margin=25mm]{geometry}
\usepackage{amsmath,amssymb,amsthm}
\usepackage{enumitem}
\usepackage{tikz}
\usepackage[colorlinks=true,linkcolor=blue,citecolor=red,urlcolor=blue]{hyperref}

\newcommand{\Z}{\mathbb{Z}}
\newcommand{\F}{\mathbb{F}}

\theoremstyle{plain}
\newtheorem{theorem}{Theorem}[section]
\newtheorem{lemma}[theorem]{Lemma}
\newtheorem{proposition}[theorem]{Proposition}
\newtheorem{corollary}[theorem]{Corollary}
\theoremstyle{definition}

\newtheorem{example}[theorem]{Example}
\newtheorem{question}[theorem]{Question}

\setlist[enumerate,1]{label=\textup{(\roman*)},leftmargin=*}
\title{\bf Good slopes, congruence preservation,\\
and gaps over residue rings\thanks{This work has been granted by the Juan de la Cierva 2024 grant (JDC2024-053082-I), funded by Ministerio de Ciencia, Innovación y Universidades (MICIU/AEI/10.13039/501100011033) and co-funded by the European Social Fund Plus (FSE+).}}
\author{
Dipak K. Bhunia\\
\small Departament de Matem\`atiques,\\ \small Universitat Polit\`ecnica de Catalunya,\\
\small Barcelona, Spain
}
\date{}

\begin{document}
\maketitle

\begin{abstract}
For $n\ge2$, a good slope of a permutation $f$ of $\mathbb{Z}_n=\mathbb{Z}/n\mathbb{Z}$
is an element $m\in\mathbb{Z}_n$ for which $x\mapsto f(x)+mx$ is a
permutation. We call $f$ congruence-preserving if it preserves
congruences modulo every divisor of $n$. For odd prime powers $n$, we
classify all permutations with more than $n/2$ good slopes. Among
permutations that are not congruence-preserving, we determine the
maximum number of good slopes for $n=3^e$ with $e\ge2$, and for
$n=25$, $125$ and $625$. Every permutation of $\mathbb{Z}_n$ has at most
$\varphi(n)$ good slopes, where $\varphi$ is Euler's totient function.
For several families of moduli with at least two distinct prime
divisors, we prove a sharp upper bound for the number of good slopes
of the permutations with fewer than $\varphi(n)$ good slopes. These
families include $9q$ and $3q^e$ for every prime $q\ge5$ and every
$e\ge1$. Together with a companion paper, this settles the bound for
every odd composite modulus below $105$.


\end{abstract}

\section{Introduction}

Let $p$ be an odd prime and $e\ge2$. Write $\Z_{p^e}=\Z/p^e\Z$,
let $\F_p$ be the field with $p$ elements, and let $\varphi$ be
Euler's totient function. For a permutation $f$ of
$\Z_{p^e}$, let $D_{p^e}(f)$ be the set of all $m\in\Z_{p^e}$ such that
$x\mapsto f(x)+mx$ is a permutation. We call its elements
\emph{good slopes}. The map $f$ is \emph{congruence-preserving} if
$x\equiv y\pmod{p^i}$ implies $f(x)\equiv f(y)\pmod{p^i}$ for
every $0\le i\le e$. For such a map, the $j$th output digit in
base $p$ depends only on the input digits in positions $0,\dots,j$.
Once the input digits in positions $0,\dots,j-1$ are fixed, the $j$th
output digit is a function of the $j$th input digit. These maps of
$\F_p$ are called \emph{digit maps}. A digit map is affine if it has the form
$t\mapsto\alpha t+\beta$ with $\alpha,\beta\in\F_p$. We call
$\alpha$ its multiplier.

In the companion paper \cite{Extremal}, we classify the permutations
of $\Z_{p^e}$ with the largest possible number
$\varphi(p^e)=(p-1)p^{e-1}$ of good slopes. They are
congruence-preserving, and their digit maps are affine with a common
multiplier. That paper also applies its bounds to a construction of
Ball and Simoens \cite{BS}, which produces mutually orthogonal quantum
Latin squares from a non-affine permutation and a set of its good
slopes with invertible differences. The main question of this paper is
whether many good slopes, but fewer than $\varphi(p^e)$, already force
congruence preservation. The results of \cite{Extremal} that we need
are stated explicitly in Section~\ref{sec:prelim}.

Over a finite field, direction theorems address the analogous
question. Good slopes are defined in the same way for a map $h$ from a
finite field to itself. A slope $a$ is \emph{determined} by the graph
of $h$ if $h(x)-h(y)=a(x-y)$ for some $x\ne y$. Thus, $a$ is not
determined precisely when $-a$ is a good slope of $h$. The direction
theorems of R\'edei, and later of Blokhuis, Ball, Brouwer, Storme and
Sz\H{o}nyi, show that sufficiently many good slopes impose strong
structural restrictions \cite{Ball,BBBSS,BBS,Redei}. De Beule surveys
direction problems in affine spaces \cite{DeBeule}. Over a residue
ring, the successive congruences provide additional structure that is
absent over a field. Problems of R\'edei type have also been studied
over finite chain rings, a class of rings that contains $\Z_{p^e}$.
For instance, Landjev and Boev \cite{LB} study blocking sets of
R\'edei type in projective Hjelmslev planes over finite chain rings.

For $n\ge2$, the complete mappings and the orthomorphisms of
$\Z_n=\Z/n\Z$ are the permutations $f$ of $\Z_n$ for which $f(x)+x$
and $f(x)-x$, respectively, are permutations. We define good slopes
and the set $D_n(f)$ for permutations $f$ of $\Z_n$ in the same way as
above. Thus, complete mappings and orthomorphisms are the permutations
for which $1$, respectively $-1$, is a good slope. The book of Evans
treats them for finite groups \cite{Evans}. Stones and Wanless
\cite{SW} study orthomorphisms of $\Z_n$. In particular, they consider
orthomorphisms that preserve congruences modulo a single
divisor or modulo every divisor of $n$. C\'egielski, Grigorieff and
Guessarian characterise the congruence-preserving maps between residue
rings by means of polynomials with rational coefficients \cite{CGG}.
Over the ring of $p$-adic integers, the maps that preserve congruences
modulo every power of $p$ are called compatible. Anashin studies these
maps, in particular those that are bijective modulo every power of $p$
\cite{Anashin}.

In this paper, we prove that a permutation with more than
$\varphi(p^e)/2$ good slopes preserves every congruence.
Its digit maps are affine, and at most $(p-1)/2$ distinct
multipliers occur. Conversely, every congruence-preserving permutation
with affine digit maps and at most $(p-1)/2$ distinct multipliers has
more than $p^e/2$ good slopes. Thus, the thresholds $\varphi(p^e)/2$
and $p^e/2$ give the same class of permutations. The digit
description also gives an explicit counting formula.

The proof bounds the number of good slopes of permutations that are not
congruence-preserving. For $p\ge5$, we obtain the upper bound
$(p-1)^2p^{e-2}/2$. For $p=3$, the exact maximum is $3^{e-1}$.
The main tool is a family of linear forms. A linear form $aX+bY$ is
\emph{good} on the graph of $f$ if $x\mapsto ax+bf(x)$ is a bijection
of $\Z_{p^e}$. For instance, $m$ is a good slope of $f$ exactly when
$mX+Y$ is good. We group good forms by their reductions modulo $p$.
The good slopes in a fixed residue class modulo $p$ are controlled by
certain maps of $\Z_{p^{e-1}}$, and relations between these maps
control the induction on $e$. Moment identities, obtained by summing
the values of a good form over the points at which it has a fixed
residue modulo $p$, sharpen these bounds for $\Z_{p^2}$ when $p\ge5$.
For $n=25$, $125$ and $625$, they determine the largest number of good
slopes of a permutation of $\Z_n$ that is not congruence-preserving,
namely $7$, $35$ and $175$.

Our constructions give permutations with $(2p-3)p^{e-2}$ good slopes
that are not congruence-preserving, for every odd prime $p$ and every
$e\ge2$. The exact maximum remains undetermined in general. For
$p=5$, the first exponent not settled by these results is $5$, where
our lower and upper bounds are $875$ and $880$.

We also study the gap below the maximum. Every permutation of $\Z_n$
has at most $\varphi(n)$ good slopes \cite{Extremal}. A permutation is
\emph{extremal} if it attains this bound. For odd composite $n$, let
$q$ be its largest prime divisor and put $G(n)=\varphi(n)(q-2)/(q-1)$.
We prove that every non-extremal permutation of $\Z_n$ has at most
$G(n)$ good slopes when $n$ is an odd prime power, when $n=9q$ with
$q\ge5$ prime, and when $n=pq^e$ with $e\ge1$ and odd primes $p<q$
such that $q\ge2p-1$. These bounds are sharp. Together with the results
recalled from \cite{Extremal}, they settle every odd composite modulus
below $105$. The general question remains open.

The paper is organised as follows. Section~\ref{sec:prelim} collects
the notation and preliminary results. Section~\ref{sec:main} states
the main bounds and classifies permutations with more than $p^e/2$
good slopes. Section~\ref{sec:nonCP} develops residue reductions and
constructions over prime powers. Section~\ref{sec:cube} studies good
linear forms and proves the reduction lemmas for the induction.
Section~\ref{sec:sharpsquares} obtains sharper bounds modulo $p^2$
from moment identities. Section~\ref{sec:higher} proves the bounds
for all odd prime powers. Section~\ref{sec:five} gives sharper
results for powers of $5$, including the exact maxima for $25$,
$125$ and $625$. Section~\ref{sec:gap} proves gap results for odd
composite moduli. Finally,
Section~\ref{sec:questions} presents the remaining questions.

\section{Preliminaries}\label{sec:prelim}

In this section, we fix our notation and recall the results of
\cite{Extremal} that we use. We also recall a form of R\'edei's
theorem for sets of points and the digit description of
congruence-preserving permutations.

For $n\ge2$, we write $\Z_n=\Z/n\Z$ and let $\Z_n^\times$ denote
its group of units, that is, its invertible elements. We let
$\varphi$ be Euler's totient function and $p_0(n)$ the least prime
divisor of $n$. We write $p^e\,\|\,n$ when $p$ is prime, $e\ge1$,
$p^e\mid n$ and $p^{e+1}\nmid n$. Unless otherwise specified,
$p$ is an odd prime. For a prime power $q$, we write $\F_q$ for the
field with $q$ elements. We identify $\Z_p$ with $\F_p$. For $e\ge2$,
the rings $\Z_{p^e}$ and $\F_{p^e}$ are not isomorphic. We put $[P]=1$
if a statement $P$ holds and $[P]=0$ otherwise. For a subset $S$ of a
ring and elements $a$ and $u$ of the ring, we write
$S+a=\{s+a:s\in S\}$, $S-a=S+(-a)$, $uS=\{us:s\in S\}$ and
$-S=\{-s:s\in S\}$.

For a permutation $f$ of $\Z_n$, we put
$$
D_n(f)=\{m\in\Z_n:x\mapsto f(x)+mx\text{ is a permutation of }\Z_n\}.
$$
For a map $h:\F_q\to\F_q$ that need not be a permutation, we write
$D_{\F_q}(h)=\{m\in\F_q:x\mapsto h(x)+mx\text{ is a permutation}\}$.
A subset of $D_n(f)$ is \emph{admissible} if the difference of any
two distinct elements is a unit. We write $\mu_n(f)$ for the largest
size of such a subset. A permutation is \emph{extremal} if
$|D_n(f)|=\varphi(n)$. A map $g$ of $\Z_n$ is congruence-preserving if
$x\equiv y\pmod d$ implies $g(x)\equiv g(y)\pmod d$ for every divisor
$d$ of $n$. A map $x\mapsto ax+b$ of $\Z_n$, with $a,b\in\Z_n$, is
\emph{affine}. It is a permutation precisely when $a$ is a unit.
For an affine permutation $f(x)=ax+b$, we have
$D_n(f)=\Z_n^\times-a=\{u-a:\ u\in\Z_n^\times\}$, because
$x\mapsto f(x)+mx$ is a permutation exactly when $a+m$ is a unit.

For $z\ne0$ in $\Z_{p^e}$, let $\nu(z)$ be the largest $v<e$
for which $z\in p^v\Z_{p^e}$. A bar denotes reduction modulo $p$.
We identify the additive groups $\Z_{p^{e-i}}$ and
$p^i\Z_{p^e}$ by $t\mapsto p^it$. Write
$x=x_0+px_1+\cdots+p^{e-1}x_{e-1}$ with
$x_j\in\{0,\dots,p-1\}$, and identify these digits with elements
of $\F_p$.

The following normalisation, counting and decomposition results
are proved in \cite{Extremal}.

\begin{lemma}[{\cite[Lemma~2.2]{Extremal}}]\label{lem:normal}
Let $f$ be a permutation of $\Z_n$, let $a\in D_n(f)$ and $b,c\in\Z_n$, and
let $u,v\in\Z_n^\times$. Then, $x\mapsto f(x)+ax+b$ and
$x\mapsto uf(vx+c)$ are permutations of $\Z_n$, and
$$
D_n\bigl(f(x)+ax+b\bigr)=D_n(f)-a,\quad
D_n\bigl(uf(vx+c)\bigr)=uv\,D_n(f).
$$
Both operations preserve the number of good slopes, the value of $\mu_n$,
and the property of being affine.
\end{lemma}

\begin{lemma}[{\cite[Lemma~2.3]{Extremal}}]\label{lem:admissible}
\begin{enumerate}
\item Let $m,m'\in\Z_n$. Then, $m-m'\in\Z_n^\times$ if and only if
$m\not\equiv m'\pmod p$ for every prime divisor $p$ of $n$.
\item For every prime divisor $p$ of $n$, let $T_p$ be a nonempty
subset of the set of residues modulo $p$, and let $D$ be the set of $m\in\Z_n$ with
$m\bmod p\in T_p$ for every prime $p\mid n$. Then,
$|D|=n\prod_{p\mid n}(|T_p|/p)$. Moreover, the largest subsets of $D$ whose
pairwise differences are units have exactly $\min_{p\mid n}|T_p|$
elements.
\item Let $p$ be prime, let $e\ge1$, and let $f$ be a permutation
of $\Z_{p^e}$. Then,
$\mu_{p^e}(f)=\bigl|\{m\bmod p:m\in D_{p^e}(f)\}\bigr|.$
\end{enumerate}
\end{lemma}

\begin{lemma}[{\cite[Lemma~2.4]{Extremal}}]\label{lem:split-general}
Let $n=dk$, where $d,k\ge2$, and let $f$ be a permutation of $\Z_n$
that preserves congruences modulo $d$. Define $g:\Z_d\to\Z_d$ by
$g(x\bmod d)=f(x)\bmod d$. For $0\le r<d$, define
$f_r:\Z_k\to\Z_k$ by
$f(r+dt)-f(r)=d f_r(t),$
where multiplication by $d$ identifies $\Z_k$ with $d\Z_n$.
Then, $g$ and all $f_r$ are permutations. Moreover,
$$
m\in D_n(f)\quad\Longleftrightarrow\quad
m\bmod d\in D_d(g)\ \text{and}\
m\bmod k\in\bigcap_{r=0}^{d-1}D_k(f_r).
$$
\end{lemma}

We call $g$ the \emph{quotient permutation} of $f$, and we call the
maps $f_r$ the \emph{residue-class permutations} of $f$.

\begin{lemma}[{\cite[Lemma~2.6]{Extremal}}]\label{lem:powersum}
Let $p$ be a prime with $p^e\,\|\,n$, where $e\ge1$, and put
$s_i=\sum_{x=0}^{n-1}x^i$ for $i\ge1$.
\begin{enumerate}
\item If $p$ is odd, then
$s_i\equiv\frac np\sum_{x=0}^{p-1}x^i\pmod{p^e}$. Consequently,
$s_i\equiv0\pmod{p^e}$ for $1\le i\le p-2$, and
$s_{p-1}\equiv-n/p\pmod{p^e}$.
\item If $p=2$, then $s_1\equiv-n/2\pmod{2^e}$.
\end{enumerate}
In particular, $s_{p-1}\not\equiv0\pmod{p^e}$ in both cases.
\end{lemma}

\begin{theorem}[{\cite[Theorem~3.1]{Extremal}}]\label{thm:universal}
Let $f$ be a permutation of $\Z_n$, and let $p$ be a prime divisor of $n$.
Then, $D_n(f)$ does not meet every residue class modulo $p$. Consequently,
the following hold.
\begin{enumerate}
\item There is $a\in\Z_n^\times$ such that $D_n(f)\subseteq\Z_n^\times-a$.
In particular, $|D_n(f)|\le\varphi(n)$.
\item $\mu_n(f)\le p_0(n)-1$.
\item If $|D_n(f)|=\varphi(n)$, then $D_n(f)=\Z_n^\times-a$ for the unit
$a$ in (i), and $\mu_n(f)=p_0(n)-1$.
\end{enumerate}
\end{theorem}

The next finite-field bound follows from the direction theorem of
Ball \cite{Ball}, as shown in \cite[Theorem~4.6]{Extremal}. A map is additive if it
sends a sum to the sum of its values. Over a prime field, an additive
map is linear.

\begin{theorem}[{\cite{Ball}}]\label{thm:fields}
Let $q$ be a prime power and let $h:\F_q\to\F_q$ be a map.
\begin{enumerate}
\item If $q$ is odd and $h$ is not the sum of an additive map and a
constant, then $|D_{\F_q}(h)|\le(q-3)/2$.
\item If $|D_{\F_q}(h)|>q/2$, then $h$ is the sum of an additive map and a
constant.
\end{enumerate}
\end{theorem}

We also need a form of R\'edei's theorem for sets of points that are not
given as graphs. It is proved in \cite[Lemma~4.7]{Extremal}, and we include the
short proof. Let $p$ be a prime. A \emph{direction} of $\F_p^2$ is a
one-dimensional subspace $\langle v\rangle$, and the lines with direction
$\langle v\rangle$ are its cosets $P+\langle v\rangle$. The directions
$\langle(0,1)\rangle$ and $\langle(1,0)\rangle$ are called
\emph{vertical} and \emph{horizontal}, respectively, and so are the lines
with these directions. Every non-vertical direction is
$\langle(1,-m)\rangle$ for a unique $m\in\F_p$, and two points $(r,u)$ and
$(r',u')$ lie on a common line with this direction if and only if
$u+mr=u'+mr'$. A set $S\subseteq\F_p^2$ of $p$ points is \emph{injective
in the direction} $\langle v\rangle$ if every line with direction
$\langle v\rangle$ contains exactly one point of $S$. A set of $p$ points
that is contained in a line is equal to this line. A set of points
\emph{determines} a direction if some line with this direction
contains at least two of its points. Thus, a set of $p$ points is
injective in a direction exactly when it does not determine this
direction. Part (ii) of the next lemma is therefore the theorem of
R\'edei and Megyesi, which states that a set of $p$ points of
$\F_p^2$ that is not a line determines at least $(p+3)/2$ directions
\cite{Redei}. Somlai has recently given a new proof of this
theorem \cite{Somlai}, and Lov\'asz and Schrijver described the sets
that determine exactly $(p+3)/2$ directions \cite{LS}.

\begin{lemma}\label{lem:pointsets}
Let $p$ be an odd prime, and let $S\subseteq\F_p^2$ be a set of $p$ points
that is not a line.
\begin{enumerate}
\item If $S$ meets every vertical line in exactly one point, then $S$ is
injective in at most $(p-3)/2$ non-vertical directions.
\item The set $S$ is injective in at most $(p-1)/2$ directions.
\end{enumerate}
\end{lemma}

\begin{proof}
In (i), $S$ is the graph $\{(r,h(r)):r\in\F_p\}$ of a map
$h:\F_p\to\F_p$. By the description of the lines above, $S$ is injective
in the direction $\langle(1,-m)\rangle$ if and only if
$r\mapsto h(r)+mr$ is a permutation, that is, if and only if
$m\in D_{\F_p}(h)$. The graph of an affine map is a line. Hence, $h$ is
not affine. Theorem~\ref{thm:fields}(i) gives (i).

For (ii), we may assume that $S$ is injective in some direction
$\langle v\rangle$. Choose $A\in\mathrm{GL}_2(\F_p)$ with
$Av\in\langle(0,1)\rangle$. The map $A$ permutes the lines and the
directions of $\F_p^2$, and $AS$ is not a line. Moreover, $AS$ is
injective in the vertical direction, which means that it meets every
vertical line in exactly one point. By (i), $AS$ is injective in at most
$(p-3)/2$ non-vertical directions. Together with the vertical direction,
this gives at most $(p-1)/2$ directions for $AS$, and hence for $S$.
\end{proof}

Let $g$ be a congruence-preserving map of $\Z_{p^e}$, and write
$g(x)=y_0+py_1+\cdots+p^{e-1}y_{e-1}$. For each $j$ and each
tuple $u=(x_0,\dots,x_{j-1})$, there is a map
$\sigma_{j,u}:\F_p\to\F_p$ with $y_j=\sigma_{j,u}(x_j)$.
These are the digit maps of $g$. For $j=0$, the tuple is empty.
There are $J=1+p+\cdots+p^{e-1}$ digit maps.
We recall the corresponding slope criterion and count from
\cite{Extremal}.

\begin{lemma}[{\cite[Lemma~7.2]{Extremal}}]\label{lem:digits}
Let $p$ be a prime, let $e\ge1$, and let $g$ be a congruence-preserving
map of $\Z_{p^e}$ with digit maps $\sigma_{j,u}$.
\begin{enumerate}
\item The map $g$ is a permutation if and only if every $\sigma_{j,u}$ is a
permutation of $\F_p$.
\item Let $g$ be a permutation, and let $T(g)$ be the intersection of the
sets $D_p(\sigma_{j,u})$ over all digit maps. Then,
$D_{p^e}(g)=\{m\in\Z_{p^e}:m\bmod p\in T(g)\}$. Moreover,
$|D_{p^e}(g)|=p^{e-1}|T(g)|$ and $\mu_{p^e}(g)=|T(g)|$.
\end{enumerate}
\end{lemma}

\begin{proposition}[{\cite[Proposition~7.3]{Extremal}}]\label{prop:fullclass}
Let $p$ be a prime, let $e\ge1$, and let $f$ be a permutation of
$\Z_{p^e}$. The following are equivalent.
\begin{enumerate}
\item $f$ is congruence-preserving.
\item $p\Z_{p^e}\subseteq D_{p^e}(f)$.
\item $D_{p^e}(f)$ contains a full residue class modulo $p$.
\end{enumerate}
Consequently, if $f$ is not congruence-preserving, then
$|D_{p^e}(f)|\le(p-1)(p^{e-1}-1)$.
\end{proposition}

\begin{corollary}[{\cite[Corollary~7.6]{Extremal}}]\label{cor:counthalf}
Let $p$ be an odd prime, let $e\ge1$, put $J=1+p+\dots+p^{e-1}$, and let
$s$ be an integer with $1\le s\le(p+1)/2$. The number of
congruence-preserving permutations of $\Z_{p^e}$ with $p^{e-1}(p-s)$ good
slopes is
$$
\binom{p-1}{s}\,p^J\sum_{k=0}^{s}(-1)^k\binom sk(s-k)^J.
$$
Each such permutation $f$ satisfies $\mu_{p^e}(f)=p-s$.
\end{corollary}

\begin{proposition}[{\cite[Theorem~5.2 and Corollary~5.3]{Extremal}}]\label{thm:classification}\label{cor:count}
A permutation of $\Z_{p^e}$ is extremal if and only if it is
congruence-preserving and all its digit maps are affine with a
common multiplier $\alpha\in\F_p^\times$. There are
$(p-1)p^{1+p+\cdots+p^{e-1}}$ such permutations, of which
$(p-1)p^{2e-1}$ are affine.
\end{proposition}

Let $p$ be an odd prime, let $e\ge1$, and let $f$ be a permutation
of $\Z_{p^e}$. For $1\le i<e$, we say that $f$ \emph{loses congruences
at level $i$} if it does not preserve congruences modulo $p^i$.
Thus, $f$ is congruence-preserving precisely when it loses
congruences at no level. For $c\in\F_p$, put
$$
d_c(f)=|\{m\in D_{p^e}(f):\bar m=c\}|,
\quad N(f)=|D_{p^e}(f)|+d_0(f^{-1}).
$$
Section~\ref{sec:cube} interprets $N(f)$ as the number of classes
of primitive linear forms that are bijective on the graph of $f$.

\section{Rigidity from many good slopes}\label{sec:main}

In this section, we state the main bounds for permutations of
$\Z_{p^e}$ that are not congruence-preserving. We then derive the
classification of the permutations with more than $p^e/2$ good slopes.
The proof of Theorem~\ref{thm:allpowers} is completed in
Section~\ref{sec:higher}.

\begin{theorem}\label{thm:allpowers}
Let $p$ be an odd prime, let $e\ge2$, and let $f$ be a permutation
of $\Z_{p^e}$ that is not congruence-preserving. Then,
$|D_{p^e}(f)|\le\varphi(p^e)/2$. More precisely, the following hold.
\begin{enumerate}
\item If $p=3$, then $|D_{3^e}(f)|\le3^{e-1}$ and
$N(f)\le4\cdot3^{e-2}$. Both bounds are attained for every $e\ge2$.
\item If $p\ge5$, then
$|D_{p^e}(f)|\le\frac{(p-1)^2}2\,p^{e-2}$ and  $N(f)\le\frac{p(p-1)}2\,p^{e-2}=\frac{\varphi(p^e)}2.$
\end{enumerate}
\end{theorem}

The set of good slopes of a congruence-preserving permutation is a
union of residue classes modulo $p$, by Lemma~\ref{lem:digits}.
Combining this fact with Theorem~\ref{thm:allpowers} gives the
following classification.

\begin{corollary}\label{cor:allhalf}
Let $p$ be an odd prime, let $e\ge2$, and let $f$ be a permutation
of $\Z_{p^e}$. The following are equivalent.
\begin{enumerate}
\item $|D_{p^e}(f)|>p^e/2$.
\item $|D_{p^e}(f)|>\varphi(p^e)/2$.
\item The permutation $f$ is congruence-preserving, every digit map
is affine, and the number $s$ of distinct multipliers satisfies
$s\le(p-1)/2$.
\end{enumerate}
In this case,
$D_{p^e}(f)=\{m\in\Z_{p^e}:m\bmod p\notin-\Lambda\}$, $|D_{p^e}(f)|=p^{e-1}(p-s)$, and $\mu_{p^e}(f)=p-s$,
where $\Lambda\subseteq\F_p^\times$ is the set of multipliers of the
digit maps.
Corollary~\ref{cor:counthalf} counts these permutations.
\end{corollary}

\begin{proof}
Clearly, (i) implies (ii). Under (ii), Theorem~\ref{thm:allpowers}
forces $f$ to preserve congruences. By Lemma~\ref{lem:digits},
$|D_{p^e}(f)|=p^{e-1}|T(f)|$, and hence $|T(f)|\ge(p+1)/2$.
Every digit map has at least this many good slopes.
By Theorem~\ref{thm:fields}(i), every digit map is affine.
Their common good slopes are exactly $\F_p\setminus(-\Lambda)$.
Thus, $p-s=|T(f)|\ge(p+1)/2$, proving (iii).
Conversely, (iii) gives
$|D_{p^e}(f)|=p^{e-1}(p-s)\ge p^{e-1}(p+1)/2>p^e/2$.
The formulas for $D_{p^e}(f)$ and $\mu_{p^e}(f)$ follow from
Lemma~\ref{lem:digits}.
\end{proof}

For $p\ge5$, the sharper bound in Theorem~\ref{thm:allpowers}(ii) also
covers the boundary case $|D_{p^e}(f)|=\varphi(p^e)/2$.

\begin{corollary}\label{cor:fullthreshold}
Let $p\ge5$ be prime, let $e\ge2$, and let $f$ be a permutation
of $\Z_{p^e}$. Put $U_e=(p-1)^2p^{e-2}/2$. Then,
$|D_{p^e}(f)|>U_e$ if and only if $f$ is congruence-preserving and all
its digit maps are affine with at most $(p+1)/2$ distinct multipliers.
For $s$ such multipliers, $|D_{p^e}(f)|=p^{e-1}(p-s)$ and
$\mu_{p^e}(f)=p-s$. The number of these permutations for each $s$
is given by Corollary~\ref{cor:counthalf}.
\end{corollary}

\begin{proof}
We have $U_e/p^{e-1}=(p-1)^2/(2p)$, and this number lies strictly
between $(p-3)/2$ and $(p-1)/2$. If $|D_{p^e}(f)|>U_e$, then
Theorem~\ref{thm:allpowers} shows that $f$ is congruence-preserving. By
Lemma~\ref{lem:digits}, we have $p^{e-1}|T(f)|>U_e$, and hence
$|T(f)|\ge(p-1)/2>(p-3)/2$.
By Theorem~\ref{thm:fields}(i), every digit map is affine.
Thus, $|T(f)|=p-s$ and $s\le(p+1)/2$.
Conversely, these conditions give
$|D_{p^e}(f)|\ge\varphi(p^e)/2>U_e$.
The remaining assertions follow from Lemma~\ref{lem:digits} and
Corollary~\ref{cor:counthalf}.
\end{proof}

\begin{corollary}\label{cor:threehalf}
Let $e\ge2$. A permutation $f$ of $\Z_{3^e}$ has more than $3^{e-1}$
good slopes if and only if it is extremal. There are
$2\cdot3^{(3^e-1)/2}$ such permutations, of which $2\cdot3^{2e-1}$ are
affine. Every such permutation $f$ satisfies $|D_{3^e}(f)|=2\cdot3^{e-1}$
and $\mu_{3^e}(f)=2$.
\end{corollary}

\begin{proof}
For $p=3$, we have $3^{e-1}=\varphi(3^e)/2$. Hence,
Corollary~\ref{cor:allhalf} applies, and it allows only $s=1$.
By Proposition~\ref{thm:classification}, these are precisely the
extremal permutations, and the same proposition gives the counts.
Finally, Theorem~\ref{thm:universal}(iii) gives $\mu_{3^e}(f)=2$.
\end{proof}

For powers of $5$, sharper bounds determine the exact maxima
$7$, $35$ and $175$ over $\Z_{25}$, $\Z_{125}$ and $\Z_{625}$
among permutations that are not congruence-preserving.
These values are proved in Theorems~\ref{thm:exact25},
\ref{thm:exact125} and~\ref{thm:exact625}.
Corollary~\ref{cor:fivebound} treats larger exponents. For every odd
prime $p$ and every $e\ge2$, Corollary~\ref{cor:twolinesLift} gives a
permutation of $\Z_{p^e}$ with $(2p-3)p^{e-2}$ good slopes that is not
congruence-preserving.

\section{Residue reductions and permutations modulo \texorpdfstring{$p^2$}{p squared}}\label{sec:nonCP}

Throughout this section, $p$ is an odd prime and $e\ge2$.
We first reduce a permutation to one that loses congruences at every
level of a smaller ring. We then study permutations of $\Z_{p^2}$
and construct examples over higher powers of $p$.

\begin{lemma}\label{lem:split}
Let $e\ge2$, and let $f$ be a permutation of $\Z_{p^e}$ that preserves
congruences modulo $p^i$, where $1\le i\le e-1$. Let $g$ be its
quotient permutation of $\Z_{p^i}$, and let $f_r$, where $0\le r<p^i$,
be its residue-class permutations of $\Z_{p^{e-i}}$, as in
Lemma~\ref{lem:split-general} with $d=p^i$ and $k=p^{e-i}$.
\begin{enumerate}
\item An element $m\in\Z_{p^e}$ lies in $D_{p^e}(f)$ if and only if
$m\bmod p^i\in D_{p^i}(g)$ and $m\bmod p^{e-i}\in D_{p^{e-i}}(f_r)$ for
every $r$.
\item Let $1\le j<i$. Then, $f$ loses congruences at level $j$ if and only
if $g$ does. Let $1\le j<e-i$. Then, $f$ loses congruences at level $i+j$
if and only if some $f_r$ loses congruences at level $j$.
\item We have $|D_{p^e}(f)|\le p^{e-i}|D_{p^i}(g)|$, and
$|D_{p^e}(f)|\le p^i|D_{p^{e-i}}(f_r)|$ for every $r$. If $f$ is not
congruence-preserving, then $g$ or some $f_r$ is not
congruence-preserving.
\end{enumerate}
\end{lemma}

\begin{proof}
Part (i) is Lemma~\ref{lem:split-general} with $d=p^i$ and
$k=p^{e-i}$. Recall that $f(r+p^it)-f(r)=p^if_r(t)$.

For $1\le j<i$, the residue of $f(x)$ modulo $p^j$ is the residue of
$g(x\bmod p^i)$ modulo $p^j$. This gives the first assertion of (ii). Now
let $1\le j<e-i$. Two elements that are congruent modulo $p^{i+j}$ have
the form $x=r+p^it$ and $y=r+p^it'$ with $t\equiv t'\pmod{p^j}$, and
$f(x)-f(y)=p^i\bigl(f_r(t)-f_r(t')\bigr)$. Hence,
$f(x)\equiv f(y)\pmod{p^{i+j}}$ if and only if
$f_r(t)\equiv f_r(t')\pmod{p^j}$. As $x$ and $y$ vary, the triple
$(r,t,t')$ runs over all triples with $t\equiv t'\pmod{p^j}$. This gives
the second assertion of (ii).

By (i), reduction modulo $p^i$ maps $D_{p^e}(f)$ into $D_{p^i}(g)$, and
every fibre of this reduction has $p^{e-i}$ elements. This gives the first
inequality in (iii), and the second one follows in the same way. Finally,
suppose that $f$ loses congruences at some level $j$. Then, $j\ne i$, and
(ii) shows that $g$ or some $f_r$ is not congruence-preserving.
\end{proof}

For a permutation $f$ of $\Z_{p^e}$ that is not congruence-preserving, let
$L(f)\ge1$ be the largest number of consecutive integers in
$\{1,\dots,e-1\}$ at which $f$ loses congruences.

\begin{proposition}\label{prop:block}
Let $e\ge2$, and let $f$ be a permutation of $\Z_{p^e}$ that is not
congruence-preserving.
Then, there are an integer $k$ with $2\le k\le L(f)+1$ and a permutation
$h$ of $\Z_{p^k}$ that loses congruences at every level
$1\le i\le k-1$, such that
$|D_{p^e}(f)|\le p^{e-k}|D_{p^k}(h)|.$
\end{proposition}

\begin{proof}
We use induction on $e$. If $f$ loses congruences at every level
$1\le i\le e-1$, we take $h=f$ and $k=e=L(f)+1$. Otherwise, $f$ preserves
congruences modulo $p^i$ for some $1\le i\le e-1$, and we use the
notation of Lemma~\ref{lem:split}. Suppose first that $g$ is not
congruence-preserving. Then, $i\ge2$, because every permutation of $\Z_p$
is congruence-preserving. By Lemma~\ref{lem:split}(ii), the levels at which
$g$ loses congruences are the levels below $i$ at which $f$ does. Hence,
$L(g)\le L(f)$. The induction hypothesis gives an integer $k$ with
$2\le k\le L(g)+1$ and a permutation $h$ of $\Z_{p^k}$ with the required
property and $|D_{p^i}(g)|\le p^{i-k}|D_{p^k}(h)|$. By
Lemma~\ref{lem:split}(iii),
$|D_{p^e}(f)|\le p^{e-i}|D_{p^i}(g)|\le p^{e-k}|D_{p^k}(h)|$.

Otherwise, Lemma~\ref{lem:split}(iii) gives some $r$ for which $f_r$ is
not congruence-preserving, and then $e-i\ge2$. If $f_r$ loses congruences
at consecutive levels $j,j+1,\dots,j+s$, then $f$ loses congruences at the
levels $i+j,\dots,i+j+s$, by Lemma~\ref{lem:split}(ii). Hence,
$L(f_r)\le L(f)$, and we conclude as before, using the inequality
$|D_{p^e}(f)|\le p^i|D_{p^{e-i}}(f_r)|$.
\end{proof}

Now let $e=2$, and let $f$ be a permutation of $\Z_{p^2}$.
We retain the notation $d_c(f)$ for the number of good slopes
with residue $c$ modulo $p$. We identify the additive groups
$\F_p$ and $p\Z_{p^2}$ through $u\mapsto pu$.
The \emph{residue pattern} of $f$ is the function $w_f:\F_p^2\to\Z$ given
by
$$
w_f(r,a)=\bigl|\{x\in\Z_{p^2}:\bar x=r,\ \overline{f(x)}=a\}\bigr|.
$$
Every residue class modulo $p$ has $p$ elements, and $f$ is a permutation.
Hence, $\sum_{a}w_f(r,a)=\sum_{r}w_f(r,a)=p$ for all $r$ and $a$. The
permutation $f$ is congruence-preserving if and only if it preserves
congruences modulo $p$, that is, if and only if $w_f$ takes only the values
$0$ and $p$. We call $f$ \emph{spread} if $w_f(r,a)=1$ for all
$(r,a)\in\F_p^2$. Equivalently, $x\mapsto(\bar x,\overline{f(x)})$ is a
bijection from $\Z_{p^2}$ onto $\F_p^2$.

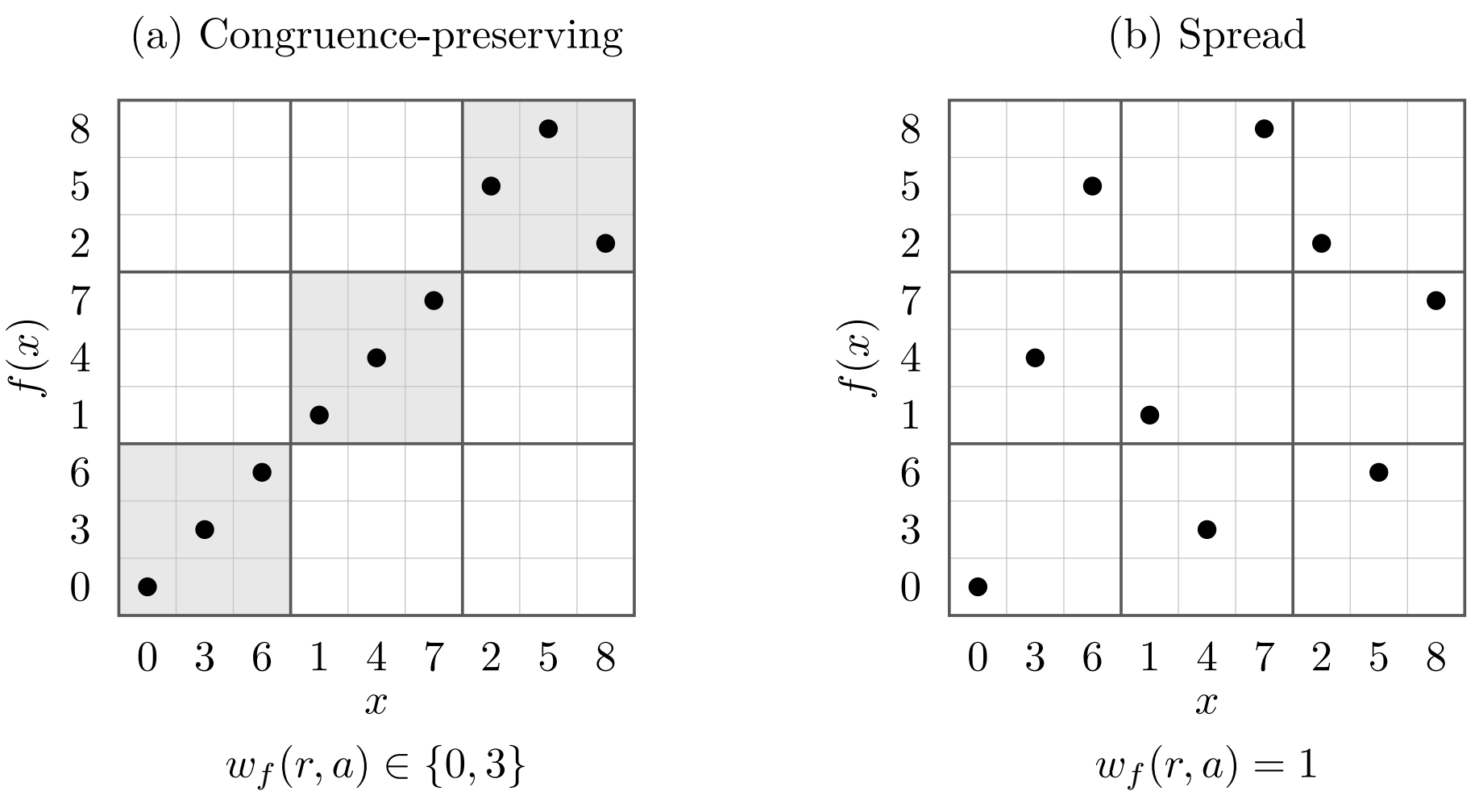
\begin{figure}[htbp]
\centering
\begin{tikzpicture}[x=.47cm,y=.47cm,every node/.style={font=\small}]
\begin{scope}[shift={(0,0)}]
\fill[black!9] (0,0) rectangle (3,3);
\fill[black!9] (3,3) rectangle (6,6);
\fill[black!9] (6,6) rectangle (9,9);
\draw[step=1,black!22,line width=.25pt] (0,0) grid (9,9);
\draw[black!65,line width=.7pt] (0,0) rectangle (9,9);
\foreach \k in {3,6}{
  \draw[black!65,line width=.7pt] (\k,0)--(\k,9);
  \draw[black!65,line width=.7pt] (0,\k)--(9,\k);
}
\foreach \v [count=\i from 0] in {0,3,6,1,4,7,2,5,8}{
  \node[below] at (\i+.5,-.2) {$\v$};
  \node[left] at (-.2,\i+.5) {$\v$};
}
\node at (4.5,-1.55) {$x$};
\node[rotate=90] at (-1.6,4.5) {$f(x)$};
\node at (4.5,10.1) {(a) Congruence-preserving};
\foreach \c/\r in {0/0,3/3,6/7,1/1,4/4,7/8,2/2,5/5,8/6}{
  \fill (\c+.5,\r+.5) circle[radius=2.2pt];
}
\node at (4.5,-2.65) {$w_f(r,a)\in\{0,3\}$};
\end{scope}
\begin{scope}[shift={(14.5,0)}]
\draw[step=1,black!22,line width=.25pt] (0,0) grid (9,9);
\draw[black!65,line width=.7pt] (0,0) rectangle (9,9);
\foreach \k in {3,6}{
  \draw[black!65,line width=.7pt] (\k,0)--(\k,9);
  \draw[black!65,line width=.7pt] (0,\k)--(9,\k);
}
\foreach \v [count=\i from 0] in {0,3,6,1,4,7,2,5,8}{
  \node[below] at (\i+.5,-.2) {$\v$};
  \node[left] at (-.2,\i+.5) {$\v$};
}
\node at (4.5,-1.55) {$x$};
\node[rotate=90] at (-1.6,4.5) {$f(x)$};
\node at (4.5,10.1) {(b) Spread};
\foreach \c/\r in {0/0,3/3,6/6,1/4,4/1,7/2,2/7,5/8,8/5}{
  \fill (\c+.5,\r+.5) circle[radius=2.2pt];
}
\node at (4.5,-2.65) {$w_f(r,a)=1$};
\end{scope}
\end{tikzpicture}
\caption{Residue patterns over $\Z_9$. Both axes follow the order
$0,3,6,1,4,7,2,5,8$. The thick lines separate residue classes modulo
$3$. In panel (a), each of the three shaded blocks
contains three graph points. In panel (b), every block contains
exactly one graph point.}
\label{fig:pattern}
\end{figure}

Figure~\ref{fig:pattern} illustrates these notions for $p=3$. The rows
and the columns are ordered by residues modulo $3$, and the thick lines
separate the residue classes. Hence, the nine blocks correspond to the
pairs $(r,a)\in\F_3^2$, and $w_f(r,a)$ is the number of points in the
block of $(r,a)$. Panel (a) shows the congruence-preserving permutation
$f(r+3t)=r+3(t+[r=2])$ of $\Z_9$, where $0\le r,t<3$.
Each of three blocks contains three points. Panel (b) shows
the spread permutation of Example~\ref{ex:Z9spread}, which has one point
in each block.

The next lemma describes the good slopes of $f$ in a fixed residue class
modulo $p$.

\begin{lemma}\label{lem:fibres}
Let $f$ be a permutation of $\Z_{p^2}$, let $c\in\F_p$, and fix
$\hat c\in\Z_{p^2}$ with residue $c$. For $z\in\F_p$, put
$B_z=\{x\in\Z_{p^2}:\overline{f(x)}+c\bar x=z\}.$
Fix $\hat z\in\Z_{p^2}$ with residue $z$, and define $u_z(x)\in\F_p$ for
$x\in B_z$ by $f(x)+\hat cx=\hat z+pu_z(x)$. Let $h\in\F_p$. Then,
$\hat c+ph\in D_{p^2}(f)$ if and only if, for every $z\in\F_p$, the set
$B_z$ has $p$ elements and the map $x\mapsto u_z(x)+h\bar x$ is injective on
$B_z$. Equivalently, for every $z$, the $p$ points $(\bar x,u_z(x))$ with
$x\in B_z$ are distinct, and the set $S_z$ formed by them is injective in
the direction $\langle(1,-h)\rangle$.
\end{lemma}

\begin{proof}
Put $m=\hat c+ph$. For $x\in B_z$, we have $phx=p(h\bar x)$, and therefore
$f(x)+mx=\hat z+p\bigl(u_z(x)+h\bar x\bigr)$. Since
$\overline{f(x)+mx}=\overline{f(x)}+c\bar x$, the sets $B_z$ are the
preimages of the residue classes modulo $p$ under $x\mapsto f(x)+mx$. This
map is a permutation if and only if it maps every $B_z$ bijectively onto
the residue class of $\hat z$, which has $p$ elements. By the formula
above, this happens if and only if $|B_z|=p$ and $x\mapsto u_z(x)+h\bar x$
is injective on $B_z$, for every $z$. The last assertion follows from the
description of the lines with direction $\langle(1,-h)\rangle$ before
Lemma~\ref{lem:pointsets}.
\end{proof}

As $h$ runs over $\F_p$, the element $\hat c+ph$ runs over the residue
class of $c$. Hence, $d_c(f)$ is the number of $h\in\F_p$ satisfying the
conditions of Lemma~\ref{lem:fibres}. The next lemma shows that a large
value of $d_c(f)$ forces $f$ to be spread.

\begin{lemma}\label{lem:dichotomy}
Let $f$ be a permutation of $\Z_{p^2}$ that is not congruence-preserving.
If $d_c(f)>(p-1)/2$ for some $c\in\F_p$, then $f$ is spread.
\end{lemma}

\begin{proof}
We use the notation of Lemma~\ref{lem:fibres}, and let $z\in\F_p$. By
Lemma~\ref{lem:fibres}, $S_z$ is a set of $p$ points, and it is injective
in the directions $\langle(1,-h)\rangle$ for $d_c(f)>(p-1)/2$ values of
$h$. By Lemma~\ref{lem:pointsets}(ii), $S_z$ is a line. The elements of
$B_z$ satisfy $\overline{f(x)}=z-c\bar x$. Hence, the number of $x\in B_z$ with
$\bar x=r$ is $w_f(r,z-cr)$. If $S_z$ is a vertical line, then all elements
of $B_z$ lie in the residue class of one element $r_z\in\F_p$. As
$|B_z|=p$, the set $B_z$ is this whole class, and $w_f(r_z,z-cr_z)=p$. If
$S_z$ is not vertical, then
$B_z$ contains exactly one element of every residue class, and
$w_f(r,z-cr)=1$ for all $r\in\F_p$.

Suppose that $S_z$ is vertical and $S_{z'}$ is not, for some $z\neq z'$.
Then, the residue class of $r_z$ is contained in $B_z$. Hence, it contains
no element of $B_{z'}$. This contradicts $w_f(r_z,z'-cr_z)=1$. Hence, either
all sets $S_z$ are vertical, or none of them is. In the first case, the
sets $B_z$ are the residue classes modulo $p$, and $f$ maps each of them
into a residue class. Then, $f$ preserves congruences modulo $p$, which is
excluded. In the second case, $w_f(r,a)=1$ for all $(r,a)$, because every
point $(r,a)\in\F_p^2$ satisfies $a+cr=z$ for some $z$. Thus, $f$ is
spread.
\end{proof}

\begin{corollary}\label{cor:notspread}
Let $f$ be a permutation of $\Z_{p^2}$ that is neither
congruence-preserving nor spread. Then, $d_c(f)\le(p-1)/2$ for every
$c\in\F_p$, and $|D_{p^2}(f)|\le(p-1)^2/2$.
\end{corollary}

\begin{proof}
The first assertion is Lemma~\ref{lem:dichotomy}. By
Theorem~\ref{thm:universal}, $d_c(f)=0$ for some $c\in\F_p$. Each of the
other $p-1$ residues contributes at most $(p-1)/2$.
\end{proof}

Now let $f$ be spread, let $c\in\F_p$, and keep the notation of
Lemma~\ref{lem:fibres}. Let $z,r\in\F_p$. The elements $x\in B_z$ with
$\bar x=r$ satisfy $\overline{f(x)}=z-cr$, and there is exactly one of them
because $w_f(r,z-cr)=1$. We define the \emph{local map}
$\psi_{c,z}:\F_p\to\F_p$ by $\psi_{c,z}(\bar x)=u_z(x)$ for $x\in B_z$.
Thus, $S_z$ is the graph of $\psi_{c,z}$, and $S_z$ is injective in the
direction $\langle(1,-h)\rangle$ if and only if $h\in D_{\F_p}(\psi_{c,z})$.
By Lemma~\ref{lem:fibres},
$$
d_c(f)=\Bigl|\bigcap_{z\in\F_p}D_{\F_p}(\psi_{c,z})\Bigr|.
$$
Replacing $\hat c$ by $\hat c+pk$ adds the map $r\mapsto kr$ to
$\psi_{c,z}$, and replacing $\hat z$ by $\hat z+pk$ subtracts the
constant $k$. We call $c$ \emph{rigid} for $f$ if all local maps
$\psi_{c,z}$ with $z\in\F_p$ are affine. By the preceding remark, this
does not depend on the choice of $\hat c$ and $\hat z$. The
\emph{multiplier} of an affine map $r\mapsto\lambda r+\beta$ of $\F_p$ is
$\lambda$.

\begin{lemma}\label{lem:rigidcount}
Let $f$ be a spread permutation of $\Z_{p^2}$, and let $c\in\F_p$. If $c$
is not rigid for $f$, then $d_c(f)\le(p-3)/2$. If $c$ is rigid for $f$ and
$K_c$ is the set of multipliers of the local maps $\psi_{c,z}$ with
$z\in\F_p$, then $d_c(f)=p-|K_c|\le p-1$.
\end{lemma}

\begin{proof}
If some $\psi_{c,z}$ is not affine, Theorem~\ref{thm:fields}(i) gives
$|D_{\F_p}(\psi_{c,z})|\le(p-3)/2$, and the formula for $d_c(f)$ gives
$d_c(f)\le(p-3)/2$. If $\psi_{c,z}(r)=\lambda r+\beta$, then
$\psi_{c,z}(r)+hr$ is a permutation of $\F_p$ exactly when
$h\neq-\lambda$. Hence, $D_{\F_p}(\psi_{c,z})=\F_p\setminus\{-\lambda\}$.
If $c$ is rigid, the formula for $d_c(f)$ therefore gives
$d_c(f)=p-|K_c|$.
\end{proof}

To compare different residues $c$, we use the geometry of
$\Z_{p^2}^2$. A \emph{line} of $\Z_{p^2}^2$ is a set
$\{P+tv:t\in\Z_{p^2}\}$, where $P,v\in\Z_{p^2}^2$ and $v$ has a unit
coordinate. We call $v$ a \emph{direction vector} of this line. The
reduction modulo $p$ of this line is the line of $\F_p^2$ through $\bar P$ with
direction $\langle\bar v\rangle$, where bars denote coordinatewise
reduction. A \emph{lifted plane} is a set $\Gamma\subseteq\Z_{p^2}^2$ of
$p^2$ points whose reductions modulo $p$ are the $p^2$ points of
$\F_p^2$. A direction $\delta$ of $\F_p^2$ is \emph{rigid} for $\Gamma$
if, for every line $\ell$ of $\F_p^2$ with direction $\delta$, the $p$
points of $\Gamma$ that reduce into $\ell$ lie on a line of
$\Z_{p^2}^2$. We write the points of $\Z_{p^2}^2$ as $(X,Y)$.

\begin{lemma}\label{lem:liftedplane}
Let $f$ be a spread permutation of $\Z_{p^2}$, and let $c\in\F_p$. Then,
the graph $\Gamma_f=\{(x,f(x)):x\in\Z_{p^2}\}$ is a lifted plane.
Moreover, $c$ is rigid for $f$ if and only if the direction
$\langle(1,-c)\rangle$ is rigid for $\Gamma_f$. In this case, for every
$z\in\F_p$, the points $(x,f(x))$ with $x\in B_z$ lie on a line with
direction vector $(1,-\hat c+p\lambda_z)$, where $\lambda_z$ is the
multiplier of $\psi_{c,z}$.
\end{lemma}

\begin{proof}
The first assertion is the definition of a spread permutation. The lines
of $\F_p^2$ with direction $\langle(1,-c)\rangle$ are the sets
$\ell_z=\{(r,a):a+cr=z\}$ with $z\in\F_p$. The points of $\Gamma_f$ that
reduce into $\ell_z$ are the points $(x,f(x))$ with $x\in B_z$. Suppose
that $\psi_{c,z}(r)=\lambda r+\beta$. For $x\in B_z$, we have
$f(x)+\hat cx=\hat z+p(\lambda\bar x+\beta)$ and $p\lambda\bar x=p\lambda
x$. Hence,
$f(x)+(\hat c-p\lambda)x=\hat z+p\beta,$
and these points lie on the line through $(0,\hat z+p\beta)$ with
direction vector $(1,-\hat c+p\lambda)$. Conversely, suppose that they lie
on a line of $\Z_{p^2}^2$. The reduction of this line contains the $p$
points of $\ell_z$. Hence, it equals $\ell_z$. After multiplying its direction
vector by a unit, we may assume that it is $(1,-\hat c+p\lambda)$ for some
$\lambda\in\F_p$. Then, $Y+(\hat c-p\lambda)X$ is constant on the line.
Thus, $f(x)+\hat cx-p\lambda\bar x$ is constant on $B_z$, which means that
$u_z(x)-\lambda\bar x$ is constant on $B_z$. Hence, $\psi_{c,z}$ is affine.
\end{proof}

The next proofs use multiplicative representatives of the residue
classes modulo $p$ in $\Z_{p^2}$. Their sums are given by an explicit
polynomial over $\F_p$.

\begin{lemma}\label{lem:carry}
Let $p$ be an odd prime. For $a\in\F_p$, let $\tau(a)\in\Z_{p^2}$ be
the residue of $\hat a^p$ modulo $p^2$, where $\hat a$ is an integer with
residue $a$ modulo $p$. Then, $\tau(a)$ does not depend on the choice of $\hat a$,
and its residue modulo $p$ is $a$. Moreover, $\tau(ab)=\tau(a)\tau(b)$,
$\tau(-a)=-\tau(a)$ and $\sum_{a\in\F_p}\tau(a)=0$. Finally, for all
$a,b\in\F_p$,
$$
\tau(a)+\tau(b)=\tau(a+b)+pE(a,b),\quad
E(a,b)=\sum_{k=1}^{p-1}\frac{(-1)^k}{k}\,a^kb^{p-k},
$$
where $E(a,b)$ is computed in $\F_p$.
\end{lemma}

\begin{proof}
By the binomial theorem, $(\hat a+pt)^p\equiv\hat a^p\pmod{p^2}$ for
every integer $t$. Thus, $\tau(a)$ is well defined. Its residue modulo
$p$ is $a$ by Fermat's theorem. Multiplicativity follows from
$(\hat a\hat b)^p=\hat a^p\hat b^p$, and $\tau(-a)=-\tau(a)$ holds
because $p$ is odd. Pairing $a$ with $-a$ gives
$\sum_a\tau(a)=0$, since $\tau(0)=0$. For $1\le k\le p-1$, we have
$$
\frac1p\binom pk=\frac{(p-1)(p-2)\cdots(p-k+1)}{k!}
\equiv\frac{(-1)^{k-1}(k-1)!}{k!}=\frac{(-1)^{k-1}}k\pmod p.
$$
Expanding $(\hat a+\hat b)^p$ with the binomial theorem therefore gives
the formula for $\tau(a)+\tau(b)$.
\end{proof}

For a lifted plane $\Gamma$ and $i,j\in\F_p$, we write $P_{ij}$ for the
point of $\Gamma$ that reduces to $(i,j)$. If the vertical and the
horizontal directions are rigid, the points $P_{ij}$ have the following
normal form.

\begin{lemma}\label{lem:normalform}
Let $\Gamma$ be a lifted plane for which the directions
$\langle(0,1)\rangle$ and $\langle(1,0)\rangle$ are rigid. Then, there
are $t_i,\nu_i,s_j,\mu_j\in\F_p$ such that
$P_{ij}=\bigl(\tau(i)+p(t_ij+\nu_i),\,\tau(j)+p(s_ji+\mu_j)\bigr)$
for all $i,j\in\F_p$.
\end{lemma}

\begin{proof}
Write $P_{ij}=(\tau(i)+px_{ij},\tau(j)+py_{ij})$, where
$x_{ij},y_{ij}\in\F_p$. Fix $i$. The points $P_{ij}$ with $j\in\F_p$
lie on a line of $\Z_{p^2}^2$ whose reduction is the vertical line
through $(i,0)$. Hence, this line has a direction vector $(pt_i,1)$ with
$t_i\in\F_p$, and $X-pt_iY$ is constant on it. At $P_{ij}$, we have
$pt_iY=pt_ij$. Hence, $x_{ij}-t_ij$ does not depend on $j$, that is,
$x_{ij}=t_ij+\nu_i$ for some $\nu_i\in\F_p$. In the same way, the
rigidity of $\langle(1,0)\rangle$ gives $y_{ij}=s_ji+\mu_j$ with
$s_j,\mu_j\in\F_p$.
\end{proof}

The next lemma is the main step for spread permutations. Together with
Lemma~\ref{lem:liftedplane}, it shows that at most two residues are rigid
for a spread permutation of $\Z_{p^2}$ when $p\ge5$.

\begin{lemma}\label{lem:tworigid}
Let $p\ge5$ be a prime. Then, every lifted plane has at most two rigid
directions.
\end{lemma}

\begin{proof}
We use the representatives $\tau$ and the polynomial $E$ of
Lemma~\ref{lem:carry}. Suppose that a lifted plane $\Gamma$ has three
rigid directions. The group
$\mathrm{GL}_2(\F_p)$ acts transitively on ordered pairs of distinct
directions. Choose $\bar A\in\mathrm{GL}_2(\F_p)$ that maps two of the
rigid directions to $\langle(0,1)\rangle$ and $\langle(1,0)\rangle$, and
let $A$ be a matrix over $\Z_{p^2}$ with reduction $\bar A$. Then, $A$ is
invertible, and it maps lines of $\Z_{p^2}^2$ to lines. Hence, $A\Gamma$ is
a lifted plane, and $\bar A$ maps the rigid directions of $\Gamma$ to those
of $A\Gamma$. Thus, we may assume that $\langle(0,1)\rangle$,
$\langle(1,0)\rangle$ and $\langle(1,-c)\rangle$ are rigid for $\Gamma$,
where $c\in\F_p^\times$.

We write the points of $\Gamma$ as in Lemma~\ref{lem:normalform}, and
we put $x_{ij}=t_ij+\nu_i$ and $y_{ij}=s_ji+\mu_j$.

Let $z\in\F_p$, and consider the points $P_{ij}$ with $j+ci=z$. Since
$\tau$ is multiplicative and $\tau(c)$ has residue $c$, at these points we
have
$$
Y+\tau(c)X=\tau(j)+\tau(ci)+p(y_{ij}+cx_{ij})=\tau(z)+p\Psi(i,j),
$$
where $\Psi(i,j)=E(j,ci)+y_{ij}+cx_{ij}$. By rigidity, these points lie on
a line of $\Z_{p^2}^2$ whose reduction has direction
$\langle(1,-c)\rangle$. As in the proof of Lemma~\ref{lem:liftedplane},
this line has a direction vector $(1,-\tau(c)-p\lambda)$ with
$\lambda\in\F_p$, and $Y+(\tau(c)+p\lambda)X$ is constant on it. Since
$p\lambda X=p\lambda i$ at $P_{ij}$, the map $i\mapsto\Psi(i,z-ci)+\lambda
i$ is constant. Thus, for every $z\in\F_p$, the map
$i\mapsto\Psi(i,z-ci)$ is affine.

For an integer $b\ge0$, put $\sigma_b=\sum_{a\in\F_p}a^b$, where $0^0=1$.
Summing over the cyclic group $\F_p^\times$ gives
$$
\sigma_b=\begin{cases}-1&\text{if $b\ge1$ and $p-1$ divides $b$},\\
0&\text{otherwise}.\end{cases}
$$
Let $k\in\{1,2\}$ and
$l=p-2-k$. Since $p\ge5$, we have $1\le l\le p-3$. Consider
$$
\Lambda_k=\sum_{i,j\in\F_p}\Psi(i,j)\,i^k(j+ci)^l.
$$
Substituting $j=z-ci$, we obtain
$\Lambda_k=\sum_zz^l\sum_i\Psi(i,z-ci)\,i^k$. If $\Psi(i,z-ci)=\alpha
i+\beta$, then the inner sum is $\alpha\sigma_{k+1}+\beta\sigma_k=0$,
because $1\le k<k+1\le p-2$. Hence, $\Lambda_k=0$.

On the other hand, we have
$\Psi(i,j)=E(j,ci)+s_ji+\mu_j+ct_ij+c\nu_i$ and
$(j+ci)^l=\sum_{r=0}^l\binom lrc^ri^rj^{l-r}$. We compute the
contributions of the five terms of $\Psi$ to $\Lambda_k$. The term $\mu_j$
contributes $\sum_j\mu_j\sum_r\binom lrc^rj^{l-r}\sigma_{k+r}$, which
vanishes because $1\le k+r\le p-2$. The term $c\nu_i$ contributes
$c\sum_i\nu_ii^k\sum_j(j+ci)^l=c\,\sigma_l\sum_i\nu_ii^k=0$. The term
$ct_ij$ contributes
$c\sum_it_ii^k\sum_r\binom lrc^ri^r\sigma_{l-r+1}=0$, because
$1\le l-r+1\le l+1\le p-2$. The term $s_ji$ contributes
$$
\sum_js_j\sum_{r=0}^l\binom lrc^rj^{l-r}\sigma_{k+1+r}=-c^lS,
\quad S=\sum_{j\in\F_p}s_j,
$$
because $1\le k+1+r\le p-1$, with equality only for $r=l$. Finally, the
term
$E(j,ci)=\sum_{a=1}^{p-1}\frac{(-1)^a}{a}c^{p-a}j^ai^{p-a}$ contributes
$$
\sum_{a=1}^{p-1}\frac{(-1)^a}{a}c^{p-a}\sum_{r=0}^l\binom lrc^r
\sigma_{a+l-r}\,\sigma_{p-a+k+r}.
$$
The two indices of $\sigma$ are positive, and their sum is
$p+k+l=2(p-1)$. Hence, the product of the two power sums is $1$ if both
indices equal $p-1$, and $0$ otherwise. Both indices equal $p-1$ exactly
when $r=a-k-1$, which requires $a\ge k+1$. In this case,
$c^{p-a+r}=c^{l+1}$. Writing $a=b+k+1$ with $0\le b\le l$, we see that the
contribution of $E(j,ci)$ is
$$
c^{l+1}\sum_{b=0}^l\frac{(-1)^{b+k+1}}{b+k+1}\binom lb.
$$
Over the rationals, this sum equals
$(-1)^{k+1}\int_0^1x^k(1-x)^l\,dx=(-1)^{k+1}k!\,l!/(k+l+1)!$. All
denominators are prime to $p$, because $b+k+1\le k+l+1=p-1$. Hence, the
identity also holds in $\F_p$. By Wilson's theorem, $(k+l+1)!=(p-1)!=-1$
in $\F_p$. Moreover, $(p-1)!=k!\prod_{i=1}^{p-1-k}(p-i)$ gives
$k!\,(p-1-k)!=(-1)^{k+1}$ in $\F_p$, and $(p-1-k)!=(p-1-k)\,l!=-(k+1)\,l!$
in $\F_p$. Therefore, $k!\,l!=(-1)^k/(k+1)$, and the contribution of
$E(j,ci)$ equals $c^{l+1}/(k+1)$.

Combining these computations, we obtain
$$
\Lambda_k=\frac{c^{p-1-k}}{k+1}-c^{p-2-k}S.
$$
For $k=1$ and $k=2$, the equations $\Lambda_k=0$ give
$c^{p-2}=2c^{p-3}S$ and $c^{p-3}=3c^{p-4}S$. Multiplying the second
equation by $2c$ and using the first one, we get $2c^{p-2}=3c^{p-2}$.
Hence, $c^{p-2}=0$, which contradicts $c\neq0$.
\end{proof}

We next rule out two rigid residues whose local maps each have a
single multiplier. The following elementary obstruction applies to
every odd prime.

\begin{lemma}\label{lem:grid}
Let $U,V\subseteq\Z_{p^2}$ be sets of representatives of the residue
classes modulo $p$, and let $\kappa\in\Z_{p^2}^\times$.
Then, $(u,v)\mapsto v-\kappa u$ is not injective on $U\times V$.
\end{lemma}

\begin{proof}
Suppose that it is injective. It is then a bijection onto $\Z_{p^2}$.
By Lemma~\ref{lem:powersum},
$$
\sum_{u\in U}\sum_{v\in V}(v-\kappa u)^{p-1}
\equiv -p\pmod{p^2}.
$$
For $1\le j\le p-2$, both $\sum_{u\in U}u^j$ and
$\sum_{v\in V}v^j$ are divisible by $p$. Thus, every mixed term in the
binomial expansion of the left-hand side vanishes modulo $p^2$.
The two remaining terms are
$$
p\sum_{v\in V}v^{p-1}
+p\kappa^{p-1}\sum_{u\in U}u^{p-1}.
$$
Each power sum is $-1$ modulo $p$, and
$\kappa^{p-1}\equiv1\pmod p$. Hence, their total is $-2p$ modulo $p^2$.
This contradicts the preceding congruence.
\end{proof}

Lemma~\ref{lem:grid} has the following consequence for lifted planes.

\begin{lemma}\label{lem:parallel}
Let $\Gamma\subseteq\Z_{p^2}^2$ be a lifted plane, and let
$v_1,v_2\in\Z_{p^2}^2$ be vectors whose reductions $\bar v_1,\bar v_2$
are linearly independent over $\F_p$. Suppose that, for $k=1,2$, the set
$\Gamma$ is contained in a union of $p$ lines with direction vector
$v_k$. Let $w\in\Z_{p^2}^2$ be a vector whose reduction lies in neither
$\langle\bar v_1\rangle$ nor $\langle\bar v_2\rangle$. Then, some line
with direction vector $w$ contains two points of $\Gamma$.
\end{lemma}

\begin{proof}
The determinant of the matrix with columns $v_1$ and $v_2$ has a
nonzero reduction modulo $p$, and therefore it is a unit. Hence, there is
$A\in\mathrm{GL}_2(\Z_{p^2})$ with $Av_1=(1,0)$ and $Av_2=(0,1)$. The
map $A$ sends every line with direction vector $v$ to a line with
direction vector $Av$. A line with direction vector $(1,0)$ has
constant second coordinate, and a line with direction vector $(0,1)$
has constant first coordinate. Thus, the points of $A\Gamma$ have at
most $p$ first coordinates and at most $p$ second coordinates. Since
there are $p^2$ points, $A\Gamma=U\times V$ for sets $U,V$ of size $p$.
The reduction of $A$ is invertible, and the points of $\Gamma$ reduce
to all points of $\F_p^2$. The same therefore holds for $A\Gamma$, and
$U$ and $V$ are sets of representatives of the residue classes modulo
$p$.

Write $Aw=(w_1,w_2)$. The reduction of $Aw$ lies in neither
$\langle(1,0)\rangle$ nor $\langle(0,1)\rangle$. Hence, $w_1$ and $w_2$
are units. Put $\kappa=w_2w_1^{-1}$. Two points $(u,v)$ and $(u',v')$ lie on
a common line with direction vector $Aw$ exactly when
$v-\kappa u=v'-\kappa u'$. By Lemma~\ref{lem:grid}, there are two
distinct points of $U\times V$ with this property. Their images under
$A^{-1}$ are two points of $\Gamma$ on a common line with direction
vector $w$.
\end{proof}

We now return to spread permutations.

\begin{lemma}\label{lem:singlemultiplier}
Let $f$ be a spread permutation of $\Z_{p^2}$. At most one residue
$c\in\F_p$ satisfies $d_c(f)=p-1$.
\end{lemma}

\begin{proof}
By Lemma~\ref{lem:rigidcount}, a residue $c$ with $d_c(f)=p-1$ is
rigid, and its local maps have a common multiplier $\lambda$. By
Lemma~\ref{lem:liftedplane}, the graph $\Gamma_f$ is then contained in
a union of $p$ lines with direction vector $(1,-\hat c+p\lambda)$, one
for each $z\in\F_p$. Suppose that two distinct residues $c_1$ and $c_2$
have this property. The vectors $(1,-c_1)$ and $(1,-c_2)$ of $\F_p^2$
are linearly independent, and $(0,1)$ lies in neither of the subspaces
that they span. By Lemma~\ref{lem:parallel}, some vertical line contains
two points of $\Gamma_f$. This is impossible, because $\Gamma_f$ is the
graph of a map.
\end{proof}

For an odd prime $p$, put $C_p=(p^2-2p+3)/2$. For instance, $C_3=3$
and $C_5=9$. Since $p(p-1)/2-C_p=(p-3)/2$, we have
$C_p\le p(p-1)/2$, with equality only for $p=3$.

\begin{theorem}\label{thm:p2}
Let $p$ be an odd prime, and let $f$ be a permutation of $\Z_{p^2}$
that is not congruence-preserving. Then,
$|D_{p^2}(f)|\le C_p\le\frac{p(p-1)}2.$
If $f$ is not spread, the stronger bound
$|D_{p^2}(f)|\le(p-1)^2/2$ holds.
\end{theorem}

\begin{proof}
The assertion for a permutation that is not spread is
Corollary~\ref{cor:notspread}. Suppose that $f$ is spread.
By Theorem~\ref{thm:universal}, $d_{c_0}(f)=0$ for some $c_0\in\F_p$.

Let $p\ge5$. By Lemmas~\ref{lem:liftedplane} and~\ref{lem:tworigid},
at most two residues are rigid. By Lemma~\ref{lem:rigidcount}, a rigid
residue contributes at most $p-1$ good slopes, and a non-rigid residue
contributes at most $(p-3)/2$. By Lemma~\ref{lem:singlemultiplier}, at
most one residue contributes exactly $p-1$. Among the $p-1$ residues
different from $c_0$, at most two are rigid, and they contribute at
most $p-1$ and $p-2$. Each of the remaining $p-3$ residues contributes
at most $(p-3)/2$. Since $p-2>(p-3)/2$, we obtain
$|D_{p^2}(f)|\le(p-1)+(p-2)+(p-3)\frac{p-3}{2}
=\frac{p^2-2p+3}{2}=C_p.$

For $p=3$, Lemma~\ref{lem:rigidcount} gives $d_c(f)=0$ for a
non-rigid residue, and $d_c(f)\le2$ for every residue. At most two
residues contribute, and Lemma~\ref{lem:singlemultiplier} rules out
two contributions equal to $2$. Hence, $|D_9(f)|\le3=C_3$.
\end{proof}

\begin{example}\label{ex:Z9spread}
Let $f$ be the involution $(3\ 4)(5\ 6)(7\ 8)$ of $\Z_9$, which fixes $0$,
$1$ and $2$. Its values are as follows.
$$
\begin{array}{|c|ccccccccc|}
\hline
x&0&1&2&3&4&5&6&7&8\\ \hline
f(x)&0&1&2&4&3&6&5&8&7\\ \hline
\end{array}
$$
The pairs $(\bar x,\overline{f(x)})$ are $(0,0)$, $(1,1)$, $(2,2)$,
$(0,1)$, $(1,0)$, $(2,0)$, $(0,2)$, $(1,2)$ and $(2,1)$. Hence, $f$ is
spread. It is not congruence-preserving, since $0\equiv3\pmod3$ while
$f(0)=0$ and $f(3)=4$. We take $\hat c=c$ and $\hat z=z$ in
$\{0,1,2\}$. Table~\ref{tab:local} lists the local maps $\psi_{c,z}$ as
triples $(\psi_{c,z}(0),\psi_{c,z}(1),\psi_{c,z}(2))$.

\begin{table}[ht]
\centering
\begin{tabular}{|c|c|c|c|}
\hline
 & $z=0$ & $z=1$ & $z=2$\\ \hline
$c=0$ & $(0,1,2)$ & $(1,0,2)$ & $(1,2,0)$\\
$c=1$ & $(0,2,2)$ & $(2,2,1)$ & $(0,0,0)$\\
$c=2$ & $(0,1,2)$ & $(0,1,2)$ & $(2,0,1)$\\ \hline
\end{tabular}
\caption{The local maps $\psi_{c,z}$ of Example~\ref{ex:Z9spread}.}
\label{tab:local}
\end{table}

The residue $1$ is not rigid, since $\psi_{1,0}$ is not affine. The
residue $0$ is rigid with multipliers $1$ and $2$, and the residue $2$ is
rigid with the single multiplier $1$. By Lemma~\ref{lem:rigidcount},
$d_0(f)=1$, $d_1(f)=0$ and $d_2(f)=2$. Indeed, $D_9(f)=\{0,2,5\}$. Thus,
Theorem~\ref{thm:p2} is sharp for $p=3$.
\end{example}

For $p\ge5$, the bound $C_p$ is not attained, as
Proposition~\ref{prop:spreadimproved} will show. The next construction
gives exactly $2p-3$ good slopes for every odd prime $p$. Its residue
$0$ is rigid with a single multiplier, and its residue $1$ is rigid with
two multipliers. For $p=3$, the number $2p-3$ equals $C_3$.

\begin{proposition}\label{prop:twolines}
Let $p$ be an odd prime. For $r,w\in\{0,1,\dots,p-1\}$, put
$\varepsilon(r,w)=[w<r]$, and let $s(r,w)$ and $t(r,w)$ be the elements
of $\{0,1,\dots,p-1\}$ with $s(r,w)\equiv w-r\pmod p$ and
$$
t(r,w)\equiv\begin{cases}
w-\varepsilon(r,w)&\text{if }w\le p-2,\\
r-1&\text{if }w=p-1
\end{cases}\pmod p.
$$
Then, there is a unique permutation $f$ of $\Z_{p^2}$ with
$f\bigl(r+p\,t(r,w)\bigr)=pr+s(r,w)$
for all $r,w\in\{0,1,\dots,p-1\}$. This permutation is spread, and
$$
D_{p^2}(f)=\{ph:h=0,\dots,p-2\}\cup
\{1+ph:h=0,\dots,p-3\}.
$$
In particular, $|D_{p^2}(f)|=2p-3$ and $\mu_{p^2}(f)=2$.
\end{proposition}

\begin{proof}
If $r=0$, then $t(0,w)=w$ for every $w$. Let $r\ge1$. As $w$ runs over
$\{r,\dots,p-2\}$, over $\{0,\dots,r-1\}$ and over $\{p-1\}$, the
number $t(r,w)$ runs over $\{r,\dots,p-2\}$, over
$\{p-1,0,\dots,r-2\}$ and over $\{r-1\}$, respectively. Hence,
$w\mapsto t(r,w)$ is a bijection of $\{0,1,\dots,p-1\}$ for every $r$,
and the elements $r+p\,t(r,w)$ are the $p^2$ elements of $\Z_{p^2}$.
The map $w\mapsto s(r,w)$ is also a bijection for every $r$. Hence, the
elements $pr+s(r,w)$ are also the $p^2$ elements of $\Z_{p^2}$, and $f$
is a well-defined permutation. For $x=r+p\,t(r,w)$, we have $\bar x=r$
and $\overline{f(x)}=w-r$. Thus, $x\mapsto(\bar x,\overline{f(x)})$ is a
bijection onto $\F_p^2$, and $f$ is spread.

Let $x=r+p\,t(r,w)$, and let $h$ be an integer. Since $phx=phr$ in
$\Z_{p^2}$, we have
$f(x)+phx=s(r,w)+p(1+h)r.$
If $h\not\equiv-1\pmod p$, this element determines $s(r,w)$ and $r$,
hence $w$, and therefore $x$. Thus, $ph$ is a good slope. If
$h\equiv-1\pmod p$, the element equals $s(r,w)$, which takes only $p$
values. Thus, $ph$ is not good.

Next, $s(r,w)+r=w+p\,\varepsilon(r,w)$ as integers. Hence,
$$
f(x)+(1+ph)x=w+p\bigl(\varepsilon(r,w)+t(r,w)+(1+h)r\bigr).
$$
If $w\le p-2$, then $\varepsilon(r,w)+t(r,w)\equiv w\pmod p$, and the
right-hand side equals $w+p\bigl(w+(1+h)r\bigr)$. If $w=p-1$, then
$\varepsilon(r,w)=0$ and $t(r,w)\equiv r-1\pmod p$. In this case, the
right-hand side equals $p-1+p\bigl((2+h)r-1\bigr)$. The element
$f(x)+(1+ph)x$ determines its residue $w$ modulo $p$. For fixed $w$, the
coefficient of $p$ is, modulo $p$, an affine function of $r$ with
multiplier $1+h$ if $w\le p-2$, and $2+h$ if $w=p-1$. Hence, $1+ph$ is a
good slope if and only if neither multiplier vanishes, that is, if and
only if $h\not\equiv-1,-2\pmod p$.

Finally, let $c\in\{2,\dots,p-1\}$. We show that no element of the
residue class of $c$ is a good slope. For $x=r+p\,t(r,w)$, we have
$f(x)+cx\equiv w+(c-1)r\pmod p$. For $r\in\{0,1,\dots,p-1\}$, let
$w_r\in\{0,1,\dots,p-1\}$ be the residue of $-1-(c-1)r$ modulo $p$.
Put $s_r=s(r,w_r)$, $t_r=t(r,w_r)$, $\varepsilon_r=\varepsilon(r,w_r)$
and $x_r=r+p\,t_r$. Then, the elements $x$ with
$f(x)+cx\equiv-1\pmod p$ are exactly $x_0,\dots,x_{p-1}$. As $r$ runs
over $\{0,1,\dots,p-1\}$, so do $w_r$ and $s_r$, because
$s_r\equiv-1-cr\pmod p$ and $c-1$ and $c$ are units modulo $p$. Summing
the integer identities $s_r+r=w_r+p\,\varepsilon_r$ therefore gives
$\sum_r\varepsilon_r=(p-1)/2$. Moreover, $w_r=p-1$ only for $r=0$, and
then $t_0=p-1=w_0-\varepsilon_0$. Hence, $t_r\equiv w_r-\varepsilon_r$
for every $r$, and
$$
\sum_rt_r\equiv\sum_rw_r-\sum_r\varepsilon_r\equiv-\frac{p-1}2
\equiv\frac12\pmod p.
$$

Now let $m=c+ph$ with $h\in\Z$. For every $r$, the integer
$s_r+cr+1$ is divisible by $p$. We write it as $p\kappa_r$. Then,
$f(x_r)+mx_r+1=p\kappa_r+p\,(r+ct_r+hr)$
in $\Z_{p^2}$. Since $s_r$ and $r$ both run over $\{0,1,\dots,p-1\}$,
we have
$$
\sum_r\kappa_r=\frac1p\sum_r(s_r+cr+1)=\frac{(1+c)(p-1)}2+1
\equiv\frac{1-c}2\pmod p.
$$
Consequently,
$\sum_r\bigl(\kappa_r+r+ct_r+hr\bigr)\equiv\frac{1-c}2+\frac c2
=\frac12\pmod p.$
If $m$ were a good slope, then $x\mapsto f(x)+mx$ would map the $p$
elements $x_r$ onto the $p$ elements of $\Z_{p^2}$ congruent to $-1$
modulo $p$. The residues of $\kappa_r+r+ct_r+hr$ modulo $p$ would then
run over $\F_p$, and their sum would be $0$. This contradiction proves
the formula for $D_{p^2}(f)$.

The good slopes lie in the residue classes of $0$ and $1$ modulo $p$.
Hence, an admissible set has at most two elements. The set $\{0,1\}$ is
admissible, and therefore $\mu_{p^2}(f)=2$.
\end{proof}

\begin{example}\label{ex:Z25seven}
For $p=5$, the permutation $f$ of Proposition~\ref{prop:twolines} is
given by the following table. The entry in row $t$ and column $r$ is
$f(r+5t)$.
$$
\begin{array}{|c|ccccc|}
\hline
 & r=0 & r=1 & r=2 & r=3 & r=4\\ \hline
t=0 & 0 & 8 & 14 & 18 & 22\\
t=1 & 1 & 5 & 12 & 19 & 23\\
t=2 & 2 & 6 & 10 & 16 & 24\\
t=3 & 3 & 7 & 11 & 15 & 20\\
t=4 & 4 & 9 & 13 & 17 & 21\\ \hline
\end{array}
$$
By Proposition~\ref{prop:twolines},
$$
D_{25}(f)=\{0,1,5,6,10,11,15\},\quad \mu_{25}(f)=2.
$$
The permutation does not preserve congruences modulo $5$, since
$f(0)=0$ and $f(5)=1$. Theorem~\ref{thm:exact25} will show that no
permutation of $\Z_{25}$ that is not congruence-preserving has more
than seven good slopes.
\end{example}

Figure~\ref{fig:slopes25} arranges these good slopes by their residues
modulo $5$. Neither of the two classes containing good slopes is full.
The circled pair $\{0,1\}$ is admissible and attains $\mu_{25}(f)=2$.

\begin{figure}[htbp]
\centering
\begin{tikzpicture}[x=.95cm,y=.72cm,every node/.style={font=\small}]
\fill[black!12] (-0.5,3.5) rectangle (0.5,4.5);
\fill[black!12] (0.5,3.5) rectangle (1.5,4.5);
\fill[black!12] (-0.5,2.5) rectangle (0.5,3.5);
\fill[black!12] (0.5,2.5) rectangle (1.5,3.5);
\fill[black!12] (-0.5,1.5) rectangle (0.5,2.5);
\fill[black!12] (0.5,1.5) rectangle (1.5,2.5);
\fill[black!12] (-0.5,0.5) rectangle (0.5,1.5);
\foreach \i in {0,...,5}{
  \draw[black!45,line width=.35pt] ({\i-.5},-.5)--({\i-.5},4.5);
  \draw[black!45,line width=.35pt] (-.5,{\i-.5})--(4.5,{\i-.5});
}
\draw[black!65,line width=.5pt] (-.5,-.5) rectangle (4.5,4.5);
\foreach \c in {0,...,4}{\node at (\c,5.1) {$\c$};}
\foreach \h in {0,...,4}{\node at (-1,{4-\h}) {$\h$};}
\node at (2,6) {$c=m\bmod5$};
\node[rotate=90] at (-1.75,2) {$h$};
\node at (0,4) {$0$};
\node at (1,4) {$1$};
\node at (2,4) {$2$};
\node at (3,4) {$3$};
\node at (4,4) {$4$};
\node at (0,3) {$5$};
\node at (1,3) {$6$};
\node at (2,3) {$7$};
\node at (3,3) {$8$};
\node at (4,3) {$9$};
\node at (0,2) {$10$};
\node at (1,2) {$11$};
\node at (2,2) {$12$};
\node at (3,2) {$13$};
\node at (4,2) {$14$};
\node at (0,1) {$15$};
\node at (1,1) {$16$};
\node at (2,1) {$17$};
\node at (3,1) {$18$};
\node at (4,1) {$19$};
\node at (0,0) {$20$};
\node at (1,0) {$21$};
\node at (2,0) {$22$};
\node at (3,0) {$23$};
\node at (4,0) {$24$};
\draw[black,line width=.8pt] (0,4) circle[radius=.23cm];
\draw[black,line width=.8pt] (1,4) circle[radius=.23cm];
\node at (2,-1.45) {$(d_0(f),\dots,d_4(f))=(4,3,0,0,0)$};
\end{tikzpicture}
\caption{The complete slope set for the permutation of $\Z_{25}$ in
Example~\ref{ex:Z25seven}. The cell in column $c$ and row $h$ contains
$m=c+5h$. Shading marks the seven good slopes. The circles mark
an admissible set of maximum size.}
\label{fig:slopes25}
\end{figure}

The construction extends to every larger exponent. More generally,
every permutation of $\Z_{p^k}$ lifts to $\Z_{p^e}$ with a proportional
number of good slopes. Recall from Theorem~\ref{thm:universal} that,
for every permutation $g$ of $\Z_{p^k}$, there is an integer $c$, prime
to $p$, with $D_{p^k}(g)\subseteq\Z_{p^k}^\times-c$.

\begin{proposition}\label{prop:slopeLift}
Let $2\le k<e$, and let $g$ be a permutation of $\Z_{p^k}$. Let $c$ be
an integer, prime to $p$, with $D_{p^k}(g)\subseteq\Z_{p^k}^\times-c$. For
$x\in\Z_{p^e}$, let $\hat x\in\{0,1,\dots,p^k-1\}$ be the residue of $x$
modulo $p^k$, and put
$f(x)=g(\hat x)+c(x-\hat x),$
where each value of $g$ is represented by its unique integer in
$\{0,1,\dots,p^k-1\}$. Then, $f$ is a
permutation of $\Z_{p^e}$, and
$$
D_{p^e}(f)=\{m\in\Z_{p^e}:m\bmod p^k\in D_{p^k}(g)\}.
$$
In particular, $|D_{p^e}(f)|=p^{e-k}|D_{p^k}(g)|$ and
$\mu_{p^e}(f)=\mu_{p^k}(g)$. If $g$ is not congruence-preserving, then
neither is $f$.
\end{proposition}

\begin{proof}
If $x\equiv y\pmod{p^k}$, then $\hat x=\hat y$ and $f(x)-f(y)=c(x-y)$. Since
$c$ is a unit, $f$ maps the residue class of $\hat x$ modulo $p^k$
bijectively onto the residue class of $g(\hat x)$. As $g$ is a permutation,
so is $f$. Moreover, $f$ preserves congruences modulo $p^k$. In the
notation of Lemma~\ref{lem:split} with $i=k$, its quotient permutation
is $g$, and every residue-class permutation is $t\mapsto ct$. By
Lemma~\ref{lem:split}(i), an element $m$ is a good slope of $f$ if and
only if $m\bmod p^k\in D_{p^k}(g)$ and $t\mapsto(c+m)t$ is a
permutation of $\Z_{p^{e-k}}$. The second condition means that $c+m$ is
a unit. It follows from the first condition, by the choice of $c$. This
proves the formula for $D_{p^e}(f)$ and its cardinality.

Two elements with a unit difference have distinct residues modulo $p$.
Hence, reduction modulo $p^k$ maps every admissible subset of
$D_{p^e}(f)$ injectively onto an admissible subset of $D_{p^k}(g)$.
Conversely, any lifts of the elements of an admissible subset of
$D_{p^k}(g)$ form an admissible subset of $D_{p^e}(f)$. This proves
$\mu_{p^e}(f)=\mu_{p^k}(g)$. Finally, suppose that $g$ loses
congruences at level $j<k$. Choose $x,y\in\{0,1,\dots,p^k-1\}$ with
$x\equiv y\pmod{p^j}$ and $g(x)\not\equiv g(y)\pmod{p^j}$. Then,
$f(x)=g(x)$ and $f(y)=g(y)$, and $f$ also loses congruences at level
$j$.
\end{proof}

\begin{corollary}\label{cor:twolinesLift}
Let $p$ be an odd prime, and let $e\ge2$. There is a permutation $f$ of
$\Z_{p^e}$ that is not congruence-preserving and satisfies
$$
|D_{p^e}(f)|=(2p-3)p^{e-2},\quad \mu_{p^e}(f)=2.
$$
\end{corollary}

\begin{proof}
For $e=2$, we take the permutation of Proposition~\ref{prop:twolines}.
Its good slopes are congruent to $0$ or $1$ modulo $p$, and hence they
lie in $\Z_{p^2}^\times-1$. For $e\ge3$, Proposition~\ref{prop:slopeLift}
with $k=2$ and $c=1$ gives the assertion.
\end{proof}

\begin{example}\label{ex:lift}
Let $f_2$ be the permutation of Example~\ref{ex:Z9spread}. Its good
slopes $0$, $2$ and $5$ lie in $\Z_9^\times-2$. Let $e\ge3$, and apply
Proposition~\ref{prop:slopeLift} with $g=f_2$, $k=2$ and $c=2$. For
$x\in\Z_{3^e}$, let $\hat x\in\{0,1,\dots,8\}$ be the residue of $x$ modulo
$9$. The resulting permutation of $\Z_{3^e}$ is
$f(x)=f_2(\hat x)+2(x-\hat x),$
and it satisfies
$$
D_{3^e}(f)=\{m\in\Z_{3^e}:m\bmod9\in\{0,2,5\}\},\quad
|D_{3^e}(f)|=3^{e-1}.
$$
If $x\equiv y\pmod9$, then $\hat x=\hat y$ and $f(x)-f(y)=2(x-y)$. Thus, $f$
preserves congruences modulo $3^i$ for $2\le i\le e-1$. It loses
congruences at level $1$, because $f(0)=0$ and $f(3)=4$. Together with
Example~\ref{ex:Z9spread}, this shows that the bound of
Theorem~\ref{thm:allpowers} is attained for $p=3$ and every $e\ge2$.
\end{example}

A different family can lose congruences at every level. We can still
determine its good slopes by an elementary digit criterion.
For this purpose, we extend the valuation $\nu$ of
Section~\ref{sec:prelim} to all of $\Z_{p^e}$ by setting $\nu(0)=e$.

\begin{lemma}\label{lem:digitcoeff}
Let $e\ge2$, and let $c_0,\dots,c_{e-1}\in\Z_{p^e}$. For $0\le x_i<p$,
consider the map
$(x_0,\dots,x_{e-1})\longmapsto\sum_{i=0}^{e-1}c_i x_i.$
This map is a bijection onto $\Z_{p^e}$ if and only if the valuations
$\nu(c_i)$ are $0,1,\dots,e-1$ in some order.
\end{lemma}

\begin{proof}
Suppose first that the valuations are distinct and none is $e$.
If two digit tuples had the same image, choose an index $i$ for which
$x_i\ne x'_i$ and $\nu(c_i)$ is least. Then,
$\nu(c_i(x_i-x'_i))=\nu(c_i)$, because $0<|x_i-x'_i|<p$.
All other nonzero summands have larger valuation. Hence, their sum
cannot cancel this term, and the map is injective. As there are $p^e$
digit tuples, it is a bijection.

Conversely, suppose that the map is bijective. For $1\le k\le e$,
let $\zeta$ be a primitive complex $p^k$th root of unity. Summing
$\zeta^{\sum_i c_i x_i}$ over the digit tuples gives
$$
0=\sum_{y\in\Z_{p^e}}\zeta^y
=\prod_{i=0}^{e-1}\left(\sum_{a=0}^{p-1}\zeta^{c_i a}\right).
$$
If $\zeta^{c_i}=1$, the corresponding factor is $p$.
Otherwise, the geometric series vanishes exactly when
$(\zeta^{c_i})^p=1$. Thus, a factor vanishes exactly when
$\nu(c_i)=k-1$. Each of $0,\dots,e-1$ must therefore occur among the
$e$ valuations, which proves necessity.
\end{proof}

\begin{proposition}\label{prop:digitreorder}
Let $e\ge2$, let $\sigma$ be a nonidentity permutation of
$\{0,\dots,e-1\}$, and
define a permutation of $\Z_{p^e}$ by
$f(x)=\sum_i p^{\sigma(i)}x_i$, where the $x_i$ are the base-$p$ digits
of $x$. Then, $f$ is not congruence-preserving. An element
$m\in\Z_{p^e}$ is a good slope of $f$ if and only if the valuations
$\nu(p^{\sigma(i)}+mp^i)$, where $0\le i<e$, are pairwise distinct and
smaller than $e$.
Moreover, $D_{p^e}(f)\subseteq p\Z_{p^e}$.
Consequently, $|D_{p^e}(f)|\le p^{e-1}$ and $\mu_{p^e}(f)=1$.
\end{proposition}

\begin{proof}
The formula permutes the digit positions. Hence, it defines a
permutation.
Let $j$ be the least index moved by $\sigma$, and put $k=\sigma^{-1}(j)$.
Then, $k>j$ and $\sigma(j)>j$. The inputs $0$ and $p^k$ agree modulo
$p^{j+1}$, whereas their images $0$ and $p^j$ do not.
Thus, $f$ is not congruence-preserving.

The criterion for $D_{p^e}(f)$ is Lemma~\ref{lem:digitcoeff} applied
to the coefficients $p^{\sigma(i)}+mp^i$. If $m$ is a unit, both the
coefficient at $i=j$ and the coefficient at $i=k$ have valuation $j$.
Hence, $m$ is not good. All good slopes are therefore divisible by $p$.
There are at most $p^{e-1}$ of them, and no two have a unit difference.
Since $0$ is good, $\mu_{p^e}(f)=1$.
\end{proof}

\begin{example}\label{ex:rotation}
Let $e\ge2$, and write $x=r+p^{e-1}t$, where
$0\le r<p^{e-1}$ and $0\le t<p$. The cyclic digit permutation
$f(x)=pr+t$ loses congruences at every level, because $f(0)=0$ and
$f(p^{e-1})=1$. By Proposition~\ref{prop:digitreorder},
a good slope must have the form $m=pu$. For such a slope,
$f(x)+mx=p(1+u)r+t\pmod{p^e}.$
The residue modulo $p$ determines $t$. The remaining map on $r$ is
bijective exactly when $1+u$ is a unit modulo $p$. Hence,
$$
D_{p^e}(f)=\{m\in\Z_{p^e}:p\mid m,\ m\not\equiv-p\pmod{p^2}\},
\quad |D_{p^e}(f)|=(p-1)p^{e-2}.
$$
For $p=3$ and $e=3$, this gives
$D_{27}(f)=\{0,3,9,12,18,21\}$. Sections~\ref{sec:cube}
to~\ref{sec:higher} develop bounds for arbitrary permutations that lose
congruences at every level.
\end{example}

\section{Good slopes in residue classes}\label{sec:cube}

In this section, we study the good slopes of a permutation of $\Z_{p^e}$
that lie in a fixed residue class modulo $p$. They are controlled by
certain maps of $\Z_{p^{e-1}}$, as Lemma~\ref{lem:cubefibres} shows. We also count
good slopes in a way that does not depend on the choice of coordinates,
by means of direction classes. The resulting reduction lemmas hold for
every exponent. Section~\ref{sec:sharpsquares} sharpens the bounds for $\Z_{p^2}$.
Section~\ref{sec:higher} then proves Theorem~\ref{thm:allpowers}
by induction on the exponent.

Throughout this section, $p$ is an odd prime and $e\ge1$.
For a permutation $f$ of $\Z_{p^e}$, the counts $d_c(f)$ satisfy
$\mu_{p^e}(f)=|\{c\in\F_p:d_c(f)>0\}|$ by
Lemma~\ref{lem:admissible}(iii). We continue to use a bar for
reduction modulo $p$. For $e\ge2$, we identify the additive groups
$\Z_{p^{e-1}}$ and $p\Z_{p^e}$ through $u\mapsto pu$.
We extend the residue pattern from Section~\ref{sec:nonCP} to
every exponent by putting
$$
w_f(r,a)=\bigl|\{x\in\Z_{p^e}:\bar x=r,\ \overline{f(x)}=a\}\bigr|,
$$
where $r,a\in\F_p$. Its row sums and column sums equal $p^{e-1}$. The
permutation $f$ preserves congruences modulo $p$ if and only if $w_f$
takes only the values $0$ and $p^{e-1}$. For $e\ge2$, the function $w_f$
is constant if and only if $w_f(r,a)=p^{e-2}$ for all $r$ and $a$.

We need an analogue of good slopes for maps that are not permutations.
For $k\ge1$ and a map $\gamma:\Z_{p^k}\to\Z_{p^k}$, let $H(\gamma)$ be
the set of
$h\in\Z_{p^k}$ for which $u\mapsto u+h\gamma(u)$ is a permutation of
$\Z_{p^k}$. It always contains $0$. Let $f$ be a permutation of
$\Z_{p^k}$, and let $m\in\Z_{p^k}$. The substitution $x=f^{-1}(u)$ shows
that $x\mapsto f(x)+mx$ is a permutation exactly when
$u\mapsto u+mf^{-1}(u)$ is one. Hence,
$D_{p^k}(f)=H(f^{-1}).$
We use the valuation $\nu$ of Section~\ref{sec:prelim} on
$\Z_{p^k}$, with $\nu(0)=k$.

\begin{lemma}\label{lem:Hmaps}
Let $k\ge1$, and let $\gamma:\Z_{p^k}\to\Z_{p^k}$ be a map.
\begin{enumerate}
\item We have $p\Z_{p^k}\subseteq H(\gamma)$ if and only if $\gamma$ is
congruence-preserving.
\item If $H(\gamma)$ contains a residue class modulo $p$, then $\gamma$ is
congruence-preserving.
\item Let $\gamma$ be congruence-preserving, with digit maps
$\sigma_{j,u}$ as in Section~\ref{sec:prelim}. Then, $h\in H(\gamma)$ if
and only if $t\mapsto t+\bar h\sigma_{j,u}(t)$ is a permutation of $\F_p$
for all $j$ and $u$. In particular, $H(\gamma)$ is a union of residue
classes modulo $p$. Moreover, $\bar h\in H(\bar\gamma)$ for every
$h\in H(\gamma)$, where $\bar\gamma:\F_p\to\F_p$ is the map induced by
$\gamma$ modulo $p$.
\end{enumerate}
\end{lemma}

\begin{proof}
Suppose that $\gamma$ is congruence-preserving, and let
$h\in p\Z_{p^k}$. Assume that $u+h\gamma(u)=v+h\gamma(v)$ for some
$u\neq v$, and put $w=\nu(u-v)$. Then,
$\gamma(u)\equiv\gamma(v)\pmod{p^w}$. Hence,
$u-v=h\bigl(\gamma(v)-\gamma(u)\bigr)$ is divisible by $p^{w+1}$, which
contradicts the definition of $w$. Thus, $h\in H(\gamma)$.

Conversely, suppose that $\gamma$ is not congruence-preserving. Then,
there are $j\ge1$ and $y,z\in\Z_{p^k}$ with $y\equiv z\pmod{p^j}$ and
$\gamma(y)\not\equiv\gamma(z)\pmod{p^j}$. In particular, $y\neq z$. Put
$v=\nu(\gamma(y)-\gamma(z))$ and $w=\nu(y-z)$. Then, $v<j\le w$. Write
$\gamma(y)-\gamma(z)=p^vc_0$ and $y-z=p^wd_0$, where $c_0$ is prime to
$p$. The element $h=-p^{w-v}d_0c_0^{-1}$ lies in $p\Z_{p^k}$, and it
satisfies $h\bigl(\gamma(y)-\gamma(z)\bigr)=-(y-z)$. Thus,
$y+h\gamma(y)=z+h\gamma(z)$, and $h\notin H(\gamma)$. This proves (i).

For (ii), let $a+p\Z_{p^k}\subseteq H(\gamma)$. If $a\in p\Z_{p^k}$, we
apply (i). Otherwise, $a$ is a unit, and the inverses of the elements of
$a+p\Z_{p^k}$ are the elements of $a^{-1}+p\Z_{p^k}$. For a unit $h$, we
have $u+h\gamma(u)=h\bigl(\gamma(u)+h^{-1}u\bigr)$. Hence,
$u\mapsto\gamma(u)+su$ is a permutation for every
$s\in a^{-1}+p\Z_{p^k}$. Put $\pi(u)=\gamma(u)+a^{-1}u$. Then, $\pi$ is
a permutation with $p\Z_{p^k}\subseteq D_{p^k}(\pi)$. By
Proposition~\ref{prop:fullclass}, $\pi$ is congruence-preserving, and so
is $\gamma(u)=\pi(u)-a^{-1}u$.

For (iii), let $h\in\Z_{p^k}$. The map $x\mapsto x+h\gamma(x)$ is
congruence-preserving. Fix $j$ and
$u=(x_0,\dots,x_{j-1})$, and put $U=\sum_{i<j}x_ip^i$.
By the definition of the digit maps, we have
$\gamma(x)\equiv G_u+p^j\sigma_{j,u}(x_j)\pmod{p^{j+1}}$, where $G_u$
depends only on $u$. Write $U+hG_u=A+p^jB$ with integers $0\le A<p^j$
and $B$. Then,
$$
x+h\gamma(x)\equiv A+p^j\bigl(x_j+\bar h\sigma_{j,u}(x_j)+B\bigr)
\pmod{p^{j+1}}.
$$
Thus, the digit map of $x+h\gamma(x)$ for $j$ and $u$ is
$t\mapsto t+\bar h\sigma_{j,u}(t)+B$. By Lemma~\ref{lem:digits}(i),
$x+h\gamma(x)$ is a permutation if and only if
$t\mapsto t+\bar h\sigma_{j,u}(t)$ is a permutation of $\F_p$ for all
$j$ and $u$. This condition depends only on $\bar h$. For $j=0$, the
digit map is $\bar\gamma$, and the condition for $j=0$ says that
$\bar h\in H(\bar\gamma)$.
\end{proof}

Let $k\ge1$, and let $\pi$ be a permutation of $\Z_{p^k}$.
The number $N(\pi)$ defined in Section~\ref{sec:prelim}
has the following interpretation and invariance under changes of
coordinates. A pair $(a,b)\in\Z_{p^k}^2$ is \emph{primitive} if at
least one of its coordinates is a unit. We identify two primitive pairs
when one is a unit multiple of the other, and we write $[a:b]$ for the
class of $(a,b)$. Such a class is \emph{good} for $\pi$ if
$x\mapsto ax+b\pi(x)$ is a permutation. Every class has a unique
representative of the form $(m,1)$ with $m\in\Z_{p^k}$, or $(1,t)$ with
$t\in p\Z_{p^k}$. A class $[m:1]$ is good exactly when
$m\in D_{p^k}(\pi)$. Substituting $y=\pi(x)$, we see that a class
$[1:t]$ is good exactly when $t\in D_{p^k}(\pi^{-1})$. Thus, $N(\pi)$ is
the number of good classes.

Reduction modulo $p$ maps the classes onto the $p+1$ points of the
projective line over $\F_p$. Each fibre has $p^{k-1}$ elements, and we
call these fibres \emph{direction classes}. For $c\in\F_p$, the
\emph{finite} direction class of $c$ consists of the classes $[m:1]$ with
$\bar m=c$. The remaining direction class consists of the classes
$[1:t]$ with $t\in p\Z_{p^k}$, and we call it \emph{vertical}. The
\emph{count} of a direction class is the number of good classes in it.
Thus, the finite direction class of $c$ has count $d_c(\pi)$, and the
vertical class has count $d_0(\pi^{-1})$. For $k\ge2$, a direction class
is \emph{regular} if it contains a whole fibre of the reduction modulo
$p^2$ consisting of good classes. Such a fibre has $p^{k-2}$ elements.

For the involution $f$ of Example~\ref{ex:Z9spread},
Figure~\ref{fig:projective9} lists all twelve classes of primitive
linear forms over $\Z_9$. The first three columns contribute the
three ordinary good slopes. The vertical column contributes one
further good class. Thus, $N(f)=3+1=4$.

\begin{figure}[htbp]
\centering
\begin{tikzpicture}[x=2.8cm,y=.85cm,every node/.style={font=\small}]
\fill[black!12] (-0.5,1.5) rectangle (0.5,2.5);
\fill[black!12] (1.5,1.5) rectangle (2.5,2.5);
\fill[black!12] (1.5,0.5) rectangle (2.5,1.5);
\fill[black!12] (2.5,1.5) rectangle (3.5,2.5);
\foreach \i in {0,...,4}{
  \draw[black!45,line width=.35pt] ({\i-.5},-.5)--({\i-.5},2.5);
}
\foreach \j in {0,...,3}{
  \draw[black!45,line width=.35pt] (-.5,{\j-.5})--(3.5,{\j-.5});
}
\draw[black!65,line width=.5pt] (-.5,-.5) rectangle (3.5,2.5);
\draw[black!65,line width=.8pt] (2.5,-.5)--(2.5,2.5);
\foreach \c in {0,1,2}{\node at (\c,3.05) {$c=\c$};}
\node at (3,3.05) {vertical};
\foreach \h in {0,1,2}{\node[left] at (-.57,{2-\h}) {$h=\h$};}
\node at (0,2) {$[0:1]$};
\node at (0,1) {$[3:1]$};
\node at (0,0) {$[6:1]$};
\node at (1,2) {$[1:1]$};
\node at (1,1) {$[4:1]$};
\node at (1,0) {$[7:1]$};
\node at (2,2) {$[2:1]$};
\node at (2,1) {$[5:1]$};
\node at (2,0) {$[8:1]$};
\node at (3,2) {$[1:0]$};
\node at (3,1) {$[1:3]$};
\node at (3,0) {$[1:6]$};
\node at (0,-1.15) {$d_0(f)=1$};
\node at (1,-1.15) {$d_1(f)=0$};
\node at (2,-1.15) {$d_2(f)=2$};
\node at (3,-1.15) {$d_0(f^{-1})=1$};
\end{tikzpicture}
\caption{Good linear forms for the permutation of
Example~\ref{ex:Z9spread}. A class $[a:b]$ represents the form
$ax+bf(x)$ up to multiplication by a unit. The finite columns
contain $[c+3h:1]$, and the vertical column contains $[1:3h]$.
Shading marks the four good classes.}
\label{fig:projective9}
\end{figure}

\begin{lemma}\label{lem:projective}
Let $k\ge2$, and let $\pi$ be a permutation of $\Z_{p^k}$.
\begin{enumerate}
\item At most $p$ direction classes have a nonzero count.
\item Let $A\in\mathrm{GL}_2(\Z_{p^k})$, and suppose that $A$ maps the
graph of $\pi$ onto the graph of a permutation $\rho$. Then,
$N(\rho)=N(\pi)$. Moreover, $A$ induces a bijection between the
direction classes of $\pi$ and of $\rho$ that preserves counts and
regularity. The residue pattern of $\rho$ is constant if and only if
that of $\pi$ is constant. Finally, for every $j$, the permutation
$\rho$ preserves congruences modulo $p^j$ if and only if $\pi$ does.
\item Suppose that $\pi$ preserves congruences modulo $p^i$, where
$1\le i<k$. Let $h$ be its quotient permutation or one of its
residue-class permutations from Lemma~\ref{lem:split}, and let $\Z_{p^j}$
be the ring on which $h$ acts. Then, $N(\pi)\le p^{k-j}N(h)$.
\end{enumerate}
\end{lemma}

\begin{proof}
Part (i) follows from Theorem~\ref{thm:universal}, since some finite
direction class contains no good slope.

For (ii), let $\Gamma$ be the graph of $\pi$. Then, $A\Gamma$ is the
graph of $\rho$. A primitive pair $(a,b)$ defines the linear form
$(X,Y)\mapsto aX+bY$. Its class is good for $\pi$ exactly when the form
is injective on $\Gamma$. This happens exactly when the form defined by
the row vector $(a,b)A^{-1}$ is injective on $A\Gamma$. As $A^{-1}$ is
invertible modulo $p$, this row vector is again primitive. We obtain a
bijection between the good classes for $\pi$ and for $\rho$. It is
compatible with reduction modulo $p$ and modulo $p^2$, which proves the
assertions about $N$, counts and regularity. The reduction of $A$
permutes the points of $\F_p^2$, and the points of $\Gamma$ above a
point $P$ are mapped onto the points of $A\Gamma$ above $\bar AP$. This
proves the assertion about residue patterns.

Write the coordinates of $A(x,\pi(x))$ as $X(x)=ax+b\pi(x)$ and
$Y(x)=cx+d\pi(x)$. Both maps are permutations, and $\rho=Y\circ X^{-1}$.
If $\pi$ preserves congruences modulo $p^j$, then so do $X$ and $Y$. A
permutation that preserves congruences modulo $p^j$ maps every residue
class modulo $p^j$ onto a residue class. Hence, its inverse preserves the
same congruences. Thus, $X^{-1}$ and $\rho$ preserve congruences modulo
$p^j$. Applying this argument to $A^{-1}$ proves the converse.

For (iii), let $[a:b]$ be good for $\pi$. The form $ax+b\pi(x)$ is a
permutation that preserves congruences modulo $p^i$. It induces the
permutation $ax+bg(x)$ of $\Z_{p^i}$, where $g$ is the quotient
permutation. On a residue class $x=r+p^it$, we have
$ax+b\pi(x)=ar+b\pi(r)+p^i\bigl(at+b\pi_r(t)\bigr)$, where
$\pi_r$ is the corresponding residue-class permutation. Hence, the
injectivity of the form on this class implies that
$t\mapsto at+b\pi_r(t)$ is injective. Thus, the
reduction of $[a:b]$ modulo $p^j$ is good for $h$, and it is again
primitive. The reduction of classes from modulus $p^k$ to modulus $p^j$
has fibres of size $p^{k-j}$, as the representatives above show. This
proves (iii).
\end{proof}

Let two different direction classes of $\pi$ have nonzero counts, and
choose a good form from each of them. Their coefficient matrix is
invertible, because the two pairs are not proportional modulo $p$.
Taking these two forms as new coordinates, we obtain a permutation
$\rho$ as in Lemma~\ref{lem:projective}(ii). In the new coordinates, the
first chosen class becomes the vertical class, and the second becomes
the finite class of $0$.

We now bound $N(\pi)$ when $\pi$ is a permutation of $\Z_{p^2}$ that is
not congruence-preserving.

\begin{proposition}\label{prop:directions}
Let $\pi$ be a permutation of $\Z_{p^2}$ that is not
congruence-preserving. Then, $N(\pi)\le p(p-1)/2$ if $p\ge5$, and
$N(\pi)\le4$ if $p=3$.
\end{proposition}

\begin{proof}
A permutation of $\Z_{p^2}$ preserves congruences modulo $p$ if and only
if it maps every residue class modulo $p$ onto a residue class. Hence,
$\pi^{-1}$ is not congruence-preserving. Moreover,
$w_{\pi^{-1}}(a,r)=w_\pi(r,a)$ for all $r,a\in\F_p$. Hence, $\pi^{-1}$ is
spread if and only if $\pi$ is spread.

Suppose first that $\pi$ is not spread. Corollary~\ref{cor:notspread},
applied to $\pi$ and to $\pi^{-1}$, gives
$|D_{p^2}(\pi)|\le(p-1)^2/2$ and $d_0(\pi^{-1})\le(p-1)/2$. Hence,
$N(\pi)\le p(p-1)/2$. For $p=3$, this bound equals $3$.

Now suppose that $\pi$ is spread, and let $\Gamma$ be its graph. By
Lemma~\ref{lem:liftedplane}, $\Gamma$ is a lifted plane. Exchanging the
coordinates maps lines to lines, and it maps $\Gamma$ onto the graph of
$\pi^{-1}$. Hence, Lemma~\ref{lem:liftedplane}, applied to $\pi$ and to
$\pi^{-1}$, shows the following. The direction $\langle(1,-c)\rangle$
is rigid for $\Gamma$ if and only if $c$ is rigid for $\pi$, and the
vertical direction is rigid for $\Gamma$ if and only if $0$ is rigid for
$\pi^{-1}$. As for direction classes, we call $d_c(\pi)$ the
\emph{count} of $\langle(1,-c)\rangle$, and we call $d_0(\pi^{-1})$ the
count of the vertical direction. Then, $N(\pi)$ is the sum of the $p+1$
counts. By
Lemma~\ref{lem:rigidcount}, every count is at most $p-1$, and the count
of a direction that is not rigid is at most $(p-3)/2$. By
Theorem~\ref{thm:universal}, the count of some direction
$\langle(1,-c_0)\rangle$ vanishes.

Suppose that two directions have count $p-1$. By
Lemma~\ref{lem:singlemultiplier}, one of them is vertical. We denote
the other one by $\langle(1,-c)\rangle$. As in the proof of
Lemma~\ref{lem:singlemultiplier}, $\Gamma$ is contained in a union of
$p$ lines with direction vector $(1,-\hat c+p\lambda)$ for some
$\lambda\in\F_p$. In the same way, the graph of $\pi^{-1}$ is contained in
a union of $p$ lines with direction vector $(1,p\lambda')$ for some
$\lambda'\in\F_p$. Exchanging the coordinates, we see that $\Gamma$ is
contained in a union of $p$ lines with direction vector $(p\lambda',1)$.
Now let $m\in\Z_{p^2}$ with $\bar m\neq c$. By
Lemma~\ref{lem:parallel}, applied with $w=(1,-m)$, some line with
direction vector $(1,-m)$ contains two points of $\Gamma$. Then, $m$ is
not a good slope of $\pi$. Hence, all other counts vanish, and
$N(\pi)=2(p-1)$. This is at most $p(p-1)/2$ for $p\ge5$, and it equals
$4$ for $p=3$.

Otherwise, at most one count equals $p-1$. Let $p\ge5$. By
Lemma~\ref{lem:tworigid}, at most two directions are rigid. Hence, at
most two of the $p$ counts different from the count of
$\langle(1,-c_0)\rangle$ exceed $(p-3)/2$, and these two counts are at
most $p-1$ and $p-2$. Therefore,
$$
N(\pi)\le(p-1)+(p-2)+(p-2)\frac{p-3}2=\frac{p(p-1)}2.
$$
For $p=3$, the counts of the directions that are not rigid vanish. At
most one count equals $2$, the other counts are at most $1$, and the
count of $\langle(1,-c_0)\rangle$ is $0$. Hence, $N(\pi)\le2+1+1=4$.
\end{proof}

Bounds on $N(\pi)$ also bound the sets $H(\gamma)$ of maps $\gamma$ that
are not congruence-preserving.

\begin{lemma}\label{lem:Hnormal}
Let $k\ge2$, and let $\gamma:\Z_{p^k}\to\Z_{p^k}$ be a map that is not
congruence-preserving. Then, $|H(\gamma)|\le p^{k-1}-1$, or there is a
permutation $\pi$ of $\Z_{p^k}$ that is not congruence-preserving and
satisfies $|H(\gamma)|\le N(\pi)$.
\end{lemma}

\begin{proof}
If $H(\gamma)$ lies in a single residue class modulo $p$, then it is not
the whole class by Lemma~\ref{lem:Hmaps}(ii), and
$|H(\gamma)|\le p^{k-1}-1$. Otherwise, choose $h_1,h_2\in H(\gamma)$ such
that $h_2-h_1$ is a unit. Put $\pi_i(u)=u+h_i\gamma(u)$ for $i=1,2$.
These maps are permutations, and so is $\pi=\pi_2\circ\pi_1^{-1}$. For
$h\in\Z_{p^k}$, put $\beta=(h-h_1)(h_2-h_1)^{-1}$ and $\alpha=1-\beta$.
Then,
$$
u+h\gamma(u)=\alpha\pi_1(u)+\beta\pi_2(u).
$$
Substituting $u=\pi_1^{-1}(y)$, we see that $h\in H(\gamma)$ if and only
if $y\mapsto\alpha y+\beta\pi(y)$ is a permutation. If $\beta$ is a unit,
then $\alpha y+\beta\pi(y)=\beta\bigl(\pi(y)+(\beta^{-1}-1)y\bigr)$, and
this happens exactly when $\beta^{-1}-1\in D_{p^k}(\pi)$. If $\beta$ is
not a unit, then $\alpha$ is a unit, and
$\alpha y+\beta\pi(y)=\alpha\bigl(y+\beta(1-\beta)^{-1}\pi(y)\bigr)$. In
this case, $h\in H(\gamma)$ exactly when
$\beta(1-\beta)^{-1}\in H(\pi)=D_{p^k}(\pi^{-1})$, and this element lies
in $p\Z_{p^k}$. The map $h\mapsto\beta$ is bijective. The map
$\beta\mapsto\beta^{-1}-1$ is injective on the units, and
$\beta\mapsto\beta(1-\beta)^{-1}$ is injective on $p\Z_{p^k}$.
Consequently, $|H(\gamma)|\le|D_{p^k}(\pi)|+d_0(\pi^{-1})=N(\pi)$.

Suppose that $\pi$ is congruence-preserving. Then,
$p\Z_{p^k}\subseteq D_{p^k}(\pi)$ by Proposition~\ref{prop:fullclass}.
For $h\in h_2+p\Z_{p^k}$, we have $\beta\in1+p\Z_{p^k}$. Then, $\beta$ is
a unit and $\beta^{-1}-1\in p\Z_{p^k}$. Hence,
$h_2+p\Z_{p^k}\subseteq H(\gamma)$, and Lemma~\ref{lem:Hmaps}(ii) shows
that $\gamma$ is congruence-preserving. This contradiction proves that
$\pi$ is not congruence-preserving.
\end{proof}

\begin{lemma}\label{lem:Hbound}
Let $\gamma:\Z_{p^2}\to\Z_{p^2}$ be a map that is not
congruence-preserving. Then, $|H(\gamma)|\le p(p-1)/2$ if $p\ge5$,
and $|H(\gamma)|\le4$ if $p=3$.
\end{lemma}

\begin{proof}
Apply Lemma~\ref{lem:Hnormal}. Its first alternative gives
$|H(\gamma)|\le p-1$. In its second alternative,
Proposition~\ref{prop:directions} gives the asserted bound.
\end{proof}

The next lemma extends Lemma~\ref{lem:fibres} to every exponent. It
describes the good slopes of $f$ in a fixed residue class modulo $p$
through the sets $H(\gamma)$.

\begin{lemma}\label{lem:cubefibres}
Let $e\ge2$, let $f$ be a permutation of $\Z_{p^e}$, let $m_0\in D_{p^e}(f)$, and put
$f_0(x)=f(x)+m_0x$. For $z\in\{0,1,\dots,p-1\}$, define
$\gamma_z:\Z_{p^{e-1}}\to\Z_{p^{e-1}}$ by
$$
\gamma_z(u)=f_0^{-1}(z+pu)\bmod p^{e-1}.
$$
Then, for every $h\in\Z_{p^{e-1}}$, the element $m_0+ph$ lies in
$D_{p^e}(f)$ if and only if $h\in\bigcap_zH(\gamma_z)$. In particular,
$d_c(f)=\bigl|\bigcap_zH(\gamma_z)\bigr|$ for $c=\bar m_0$.
\end{lemma}

\begin{proof}
By Lemma~\ref{lem:normal}, $m_0+ph\in D_{p^e}(f)$ if and only if
$ph\in D_{p^e}(f_0)$. Substituting $x=f_0^{-1}(y)$, we see that this
holds if and only if $y\mapsto y+phf_0^{-1}(y)$ is a permutation. This
map preserves every residue class modulo $p$. For $y=z+pu$ with
$u\in\Z_{p^{e-1}}$, we have
$y+phf_0^{-1}(y)=z+p\bigl(u+h\gamma_z(u)\bigr),$
because $phf_0^{-1}(y)$ depends only on the residue of $hf_0^{-1}(y)$
modulo $p^{e-1}$. Hence, the map is a permutation if and only if
$u\mapsto u+h\gamma_z(u)$ is a permutation of $\Z_{p^{e-1}}$ for every $z$.
Finally, $m_0+ph$ runs over the residue class of $m_0$ modulo $p$ as
$h$ runs over $\Z_{p^{e-1}}$.
\end{proof}

For $e\ge2$, we call $c\in\F_p$ \emph{regular} for $f$ if
$D_{p^e}(f)$ contains a residue class $m+p^2\Z_{p^e}$ with
$\bar m=c$. Thus, $c$ is regular exactly when the finite direction class
of $c$ is regular.

\begin{lemma}\label{lem:regular}
Let $e\ge2$, let $f$ be a permutation of $\Z_{p^e}$ that is not
congruence-preserving, and let $c\in\F_p$ with $d_c(f)>0$. Form the maps
$\gamma_z$ of Lemma~\ref{lem:cubefibres} at a good slope $m_0$ with
$\bar m_0=c$.
\begin{enumerate}
\item If $c$ is regular, then every $\gamma_z$ is
congruence-preserving. Moreover, the good slopes with residue $c$ form a
union of at most $p-1$ residue classes modulo $p^2$. Thus, $d_c(f)$ is a
multiple of $p^{e-2}$, and $d_c(f)\le(p-1)p^{e-2}$.
\item If $c$ is not regular, then some $\gamma_z$ is not
congruence-preserving, and $d_c(f)\le|H(\gamma_z)|$ for this map. In
particular, for $e=3$, we have $d_c(f)\le p(p-1)/2$ if $p\ge5$, and
$d_c(f)\le4$ if $p=3$.
\end{enumerate}
\end{lemma}

\begin{proof}
Put $H=\bigcap_zH(\gamma_z)$. By Lemma~\ref{lem:cubefibres}, the slope
$m_0+ph$ is good if and only if $h\in H$, and $d_c(f)=|H|$. The slopes
$m_0+ph$ with $h$ in a fixed residue class modulo $p$ form a residue
class modulo $p^2$. Hence, $c$ is regular if and only if $H$ contains a
full residue class modulo $p$.

Suppose that $c$ is regular. Then, every $H(\gamma_z)$ contains a full
residue class modulo $p$, and Lemma~\ref{lem:Hmaps}(ii) shows that every
$\gamma_z$ is congruence-preserving. By Lemma~\ref{lem:Hmaps}(iii), every
$H(\gamma_z)$ is a union of residue classes modulo $p$, and so is $H$.
If $H$ were all of $\Z_{p^{e-1}}$, then $D_{p^e}(f)$ would contain the
residue class of $m_0$ modulo $p$, and $f$ would be
congruence-preserving by Proposition~\ref{prop:fullclass}. Hence, $H$ is
a union of at most $p-1$ residue classes modulo $p$, each with $p^{e-2}$
elements. This proves (i).

Suppose that $c$ is not regular. If every $\gamma_z$ were
congruence-preserving, then $H$ would be a union of residue classes
modulo $p$, as above. Since $0\in H$, it would contain a full class, and
$c$ would be regular. Hence, some $\gamma_z$ is not
congruence-preserving, and $d_c(f)=|H|\le|H(\gamma_z)|$. For $e=3$, this
map acts on $\Z_{p^2}$, and Lemma~\ref{lem:Hbound} gives the last
assertion.
\end{proof}

The next lemma bounds $d_c(f)$ for a regular residue $c$ when the
residue pattern of $f$ is not constant.

\begin{lemma}\label{lem:regularpattern}
Let $e\ge3$, and let $f$ be a permutation of $\Z_{p^e}$ that loses
congruences at level $1$. If $w_f$ is not constant, then every
regular residue $c$ satisfies
$d_c(f)\le\frac{(p-1)p^{e-2}}2.$
\end{lemma}

\begin{proof}
Choose $m_0$ with $\bar m_0=c$ and
$m_0+p^2\Z_{p^e}\subseteq D_{p^e}(f)$, and let $f_0$ and $\gamma_z$ be as
in Lemma~\ref{lem:cubefibres}, formed at $m_0$. By
Lemma~\ref{lem:regular}(i), every $\gamma_z$ is congruence-preserving.
Let $\psi_z:\F_p\to\F_p$ be the map induced by $\gamma_z$ modulo $p$. We
claim that
$$
w_f(r,z-cr)=p^{e-2}\bigl|\{t\in\F_p:\psi_z(t)=r\}\bigr|
$$
for all $r,z\in\F_p$. Indeed, $\overline{f_0(x)}=\overline{f(x)}+c\bar x$.
Hence, $w_f(r,z-cr)$ is the number of elements $x$ with $\bar x=r$ and
$\overline{f_0(x)}=z$. The elements $x$ with
$\overline{f_0(x)}=z$ are the elements $f_0^{-1}(z+pu)$ with
$u\in\Z_{p^{e-1}}$, and they satisfy $\bar x=\psi_z(\bar u)$. Each
residue $\bar u$ has $p^{e-2}$ lifts $u$.

Suppose that every $\psi_z$ is affine. An affine map of $\F_p$ is
constant or bijective. Hence, for each $z$, either $w_f(r,z-cr)=p^{e-2}$
for all $r$, or $w_f(r,z-cr)$ equals $p^{e-1}$ for one value $r_z$ of
$r$ and vanishes for the others. Suppose that the second alternative
occurs for some $z$. As the row sums of $w_f$ equal $p^{e-1}$, the row
of $r_z$ vanishes outside the entry $(r_z,z-cr_z)$. In particular,
$w_f(r_z,z'-cr_z)=0$ for every $z'\neq z$, and the second alternative
occurs for every $z'$. The rows $r_z$ are distinct, again because the
row sums equal $p^{e-1}$. Thus, $w_f$ takes only the values $0$ and
$p^{e-1}$, and $f$ preserves congruences modulo $p$. Otherwise, the
first alternative occurs for every $z$, and $w_f$ is constant. Both
conclusions contradict the hypotheses. Hence, some $\psi_z$ is not
affine.

For $h\in\F_p^\times$, we have $t+h\psi_z(t)=h\bigl(\psi_z(t)+h^{-1}t\bigr)$.
Thus, $t\mapsto t+h\psi_z(t)$ is a permutation if and only if
$h^{-1}\in D_{\F_p}(\psi_z)$. By Theorem~\ref{thm:fields}(i), we obtain
$|H(\psi_z)|\le1+(p-3)/2=(p-1)/2$. By Lemma~\ref{lem:Hmaps}(iii), every
element of $H(\gamma_z)$ reduces modulo $p$ to an element of
$H(\psi_z)$, and each residue has $p^{e-2}$ lifts. Therefore,
$d_c(f)\le|H(\gamma_z)|\le(p-1)p^{e-2}/2$.
\end{proof}

Next, we treat permutations with a constant residue pattern. For a
permutation $f$ of $\Z_{p^e}$ and $r,a\in\F_p$, put
$$
E_{r,a}=\{x\in\Z_{p^e}:\bar x=r,\ \overline{f(x)}=a\}.
$$
These sets partition $\Z_{p^e}$, and $|E_{r,a}|=w_f(r,a)$.

\begin{lemma}\label{lem:regularsets}
Let $e\ge3$, let $f$ be a permutation of $\Z_{p^e}$ with
$w_f(r,a)=p^{e-2}$ for all $r,a$, and let $c$ be regular.
Choose $m$ with $\bar m=c$ and
$m+p^2\Z_{p^e}\subseteq D_{p^e}(f)$, and put $g(x)=f(x)+mx$.
Then, $g$ maps every $E_{r,a}$ onto a residue class modulo $p^2$.
\end{lemma}

\begin{proof}
Since $w_f$ is constant and $e\ge3$, the permutation $f$ does not
preserve congruences modulo $p$. By Lemma~\ref{lem:regular}(i), each map
$u\mapsto g^{-1}(z+pu)\bmod p^{e-1}$ is congruence-preserving. Hence,
the residue of $g^{-1}(y)$ modulo $p$ depends only on the residue of $y$
modulo $p^2$. Thus, the inverse image under $g$ of a residue class
modulo $p^2$ consists of $p^{e-2}$ elements with the same residue
modulo $p$.

Fix $z\in\F_p$, and let $X_z$ be the inverse image of $z+p\Z_{p^e}$
under $g$. An element $x$ lies in $X_z$ if and only if
$\overline{f(x)}=z-c\bar x$. Hence, $X_z$ meets the residue class of $r$
modulo $p$ in the set $E_{r,z-cr}$, which has $p^{e-2}$ elements. On the
other hand, $X_z$ is the disjoint union of the inverse images of the $p$
residue classes modulo $p^2$ contained in $z+p\Z_{p^e}$. Each of them
has $p^{e-2}$ elements and lies in a single residue class modulo $p$.
Comparing sizes, we see that each $E_{r,z-cr}$ is one of these inverse
images. As $r$ and $z$ vary, the pair $(r,z-cr)$ runs over $\F_p^2$,
and the lemma follows.
\end{proof}

Two regular residues force a preserved congruence modulo $p^2$.

\begin{lemma}\label{lem:tworegular}
Let $e\ge3$, and let $f$ be a permutation of $\Z_{p^e}$ with
$w_f(r,a)=p^{e-2}$ for all $r,a$. If two distinct residues are
regular, then $f$ preserves congruences modulo $p^2$.
\end{lemma}

\begin{proof}
Let $c_1$ and $c_2$ be the two regular residues. For $i=1,2$, choose
$m_i$ with $\bar m_i=c_i$ and $m_i+p^2\Z_{p^e}\subseteq D_{p^e}(f)$, and
put $g_i(x)=f(x)+m_ix$. By Lemma~\ref{lem:regularsets}, both $g_1$ and
$g_2$ are constant modulo $p^2$ on every set $E_{r,a}$. Hence, so is
$g_2(x)-g_1(x)=(m_2-m_1)x$. As $m_2-m_1$ is a unit, every $E_{r,a}$ lies
in a residue class modulo $p^2$. Since $E_{r,a}$ and every residue class
modulo $p^2$ have $p^{e-2}$ elements, the sets $E_{r,a}$ are exactly the
residue classes modulo $p^2$. By Lemma~\ref{lem:regularsets}, $g_1$ maps
each of them onto a residue class modulo $p^2$. Thus, $g_1$ preserves
congruences modulo $p^2$, and so does $f(x)=g_1(x)-m_1x$.
\end{proof}

A regular residue with a large count imposes further restrictions. To
transfer them from one digit to the next, we need three elementary facts
about digit maps. We say that a congruence-preserving map of $\Z_{p^r}$
has \emph{affine digits} if all its digit maps are affine. The
multipliers of these affine maps may be zero. For instance, the
identity, the constant maps and the maps $x\mapsto bx$ with
$b\in\Z_{p^r}$ have affine digits.

\begin{lemma}\label{lem:affine-digits}
Let $p$ be an odd prime, and let $r\ge1$.
\begin{enumerate}
\item The maps of $\Z_{p^r}$ with affine digits are closed under
addition, under multiplication by elements of $\Z_{p^r}$, and under
composition.
\item If a permutation of $\Z_{p^r}$ has affine digits, then so does its
inverse.
\item If $\eta$ is a map of $\Z_{p^r}$ with affine digits, then
$\sum_{u=0}^{p^t-1}\eta(u)\equiv0\pmod{p^t}$ for $1\le t\le r$.
\end{enumerate}
\end{lemma}

\begin{proof}
Let $\eta_1$ and $\eta_2$ have affine digits. The maps $\eta_1+\eta_2$ and
$\eta_1\circ\eta_2$ are congruence-preserving. Fix $j$ and the input digits
$x_0,\dots,x_{j-1}$. Then, the lower $j$ digits of $\eta_1(x)$ and of
$\eta_2(x)$ are fixed. Hence, the carry into position $j$ of
$\eta_1(x)+\eta_2(x)$ does not depend on $x_j$. Thus, the digit map of
$\eta_1+\eta_2$ is the sum of the digit maps of $\eta_1$ and $\eta_2$ and a constant,
and it is affine. Multiplication by an element of $\Z_{p^r}$ is repeated
addition. For the composition, the lower $j$ digits of $\eta_2(x)$ are
fixed, and they determine which digit map of $\eta_1$ acts at position $j$.
Hence, the digit map of $\eta_1\circ\eta_2$ is a composition of two affine
maps, and it is affine. This proves (i).

Let $\eta$ be a permutation with affine digits. By Lemma~\ref{lem:digits}(i),
its digit maps are permutations of $\F_p$. The inverse of $\eta$ is
congruence-preserving, because $\eta$ maps every residue class modulo
$p^j$ onto a residue class. At a given prefix of output digits, the
digit map of $\eta^{-1}$ is the inverse of the digit map of $\eta$ at the
corresponding prefix of input digits. This inverse is affine, which
proves (ii).

For (iii), the residue of $\eta(u)$ modulo $p^t$ depends only on $u$ modulo
$p^t$. Write it in base $p$ with digits $y_j(u)\in\{0,\dots,p-1\}$, where
$0\le j<t$. Fix $j$. The digit $y_j(u)$ is the value of the digit map
$\sigma_{j,v}$ at the $j$th digit of $u$, where $v$ is the prefix of the
lower $j$ digits of $u$. The digits of $u$ above position $j$ are
arbitrary. Hence,
$$
\sum_{u=0}^{p^t-1}y_j(u)=p^{t-j-1}\sum_v\sum_{a=0}^{p-1}\sigma_{j,v}(a),
$$
where the values of $\sigma_{j,v}$ are represented in
$\{0,\dots,p-1\}$. An affine map of $\F_p$ is constant or bijective. As
$p$ is odd, the sum $0+1+\dots+(p-1)$ is divisible by $p$. Thus, each
inner sum is divisible by $p$, and $p^j\sum_uy_j(u)$ is divisible by
$p^t$. Summing over $j$ proves (iii).
\end{proof}

\begin{lemma}\label{lem:cycle}
Let $p$ be an odd prime, let $r\ge1$, and let $b\in\Z_{p^r}^\times$. Let
$\eta$ be a map of $\Z_{p^r}$ with affine digits. Then, the map
$A(x)=x+b+p\eta(x)$ is a permutation of $\Z_{p^r}$ consisting of a single
cycle.
\end{lemma}

\begin{proof}
Let $x\neq y$, and put $v=\nu(x-y)$. Since $\eta$ is congruence-preserving,
$p\bigl(\eta(x)-\eta(y)\bigr)$ is divisible by $p^{v+1}$. Hence,
$A(x)-A(y)\equiv x-y\pmod{p^{v+1}}$. Thus, $A$ is a permutation, and it
preserves congruences modulo every power of $p$. Modulo $p$, the map $A$
is the translation by the nonzero residue of $b$, and hence it is a
single cycle.

Suppose that $A$ induces a single cycle on $\Z_{p^t}$, where $1\le t<r$,
and let $x\in\Z_{p^r}$. The elements $x,A(x),\dots,A^{p^t-1}(x)$ then
represent every residue class modulo $p^t$ exactly once. As $\eta$ is
congruence-preserving, Lemma~\ref{lem:affine-digits}(iii) shows that the
sum of the values of $\eta$ at these elements is divisible by $p^t$. Adding
the increments $A(y)-y=b+p\eta(y)$ along this orbit segment, we obtain
$$
A^{p^t}(x)-x\equiv p^tb\pmod{p^{t+1}}.
$$
Hence, $A^{p^t}$ maps the residue class of $x$ modulo $p^t$ onto itself,
and it acts on the $p$ residue classes modulo $p^{t+1}$ contained in it
as the translation by $p^tb$. This translation is a single cycle.
Therefore, $A$ induces a single cycle on $\Z_{p^{t+1}}$, and induction on
$t$ proves the lemma.
\end{proof}

\begin{lemma}\label{lem:connected}
Let $p$ be an odd prime, let $r\ge1$, and let $b\in\Z_{p^r}^\times$. For
each $z\in\{0,\dots,p-1\}$, let $\theta_z$ be a map of $\Z_{p^r}$ with
affine digits, and put
$$
T_z(U)=bU+z+p\theta_z(U)
$$
for $U\in\Z_{p^r}$. Form the bipartite graph whose vertices are two
copies of $\Z_{p^r}$, in which $U$ in the first copy is joined to
$T_z(U)$ in the second copy for every $z$. Then, this graph is
connected.
\end{lemma}

\begin{proof}
By Lemma~\ref{lem:affine-digits}(i), every $T_z$ has affine digits. Let
$U\neq V$, and put $v=\nu(U-V)$. As $\theta_z$ is congruence-preserving,
we have $T_z(U)-T_z(V)\equiv b(U-V)\pmod{p^{v+1}}$. Hence, $T_z$ is
injective, and it is a permutation. By Lemma~\ref{lem:affine-digits},
the permutation $A=T_0^{-1}\circ T_1$ has affine digits. The identity
$T_0(A(U))=T_1(U)$ reads
$$
bA(U)+p\theta_0(A(U))=bU+1+p\theta_1(U).
$$
Put $\beta=b^{-1}$. Then, $A(U)=U+\beta+p\eta(U)$, where
$\eta(U)=\beta\bigl(\theta_1(U)-\theta_0(A(U))\bigr)$. By
Lemma~\ref{lem:affine-digits}(i), the map $\eta$ has affine digits. By
Lemma~\ref{lem:cycle}, $A$ consists of a single cycle.

For every $U$, the vertex $T_1(U)=T_0(A(U))$ is joined to $U$ and to
$A(U)$. Since $A$ is a single cycle, all vertices of the first copy lie
in one connected component. Every vertex of the second copy has the form
$T_0(U)$, and it is joined to $U$. Hence, the graph is connected.
\end{proof}

The next lemma compares the maps of Lemma~\ref{lem:cubefibres} formed
at two good slopes whose difference is a unit. If the maps formed at
the first slope are congruence-preserving permutations, it determines
the good slopes in $p\Z_{p^{e-1}}$ of the maps formed at the second
slope.

\begin{lemma}\label{lem:zerotransfer}
Let $e\ge3$, let $g$ be a permutation of $\Z_{p^e}$, and put $F=g^{-1}$.
Suppose that the maps
$\alpha_z(u)=F(z+pu)\bmod p^{e-1},$
where $z\in\{0,1,\dots,p-1\}$ and $u\in\Z_{p^{e-1}}$, are
congruence-preserving permutations of $\Z_{p^{e-1}}$. Let
$b\in\Z_{p^e}^\times$ be such that $x\mapsto g(x)+bx$ is a permutation,
and put
$$
\gamma_w(v)=(g+b\,\mathrm{id})^{-1}(w+pv)\bmod p^{e-1}
$$
for $w\in\{0,1,\dots,p-1\}$ and $v\in\Z_{p^{e-1}}$. Then, every
$\gamma_w$ is a permutation of $\Z_{p^{e-1}}$. Its good slopes in
$p\Z_{p^{e-1}}$ are the elements $ph$ with $1+\bar b\bar h\neq0$. In
particular, $d_0(\gamma_w)=(p-1)p^{e-3}$.
\end{lemma}

\begin{proof}
Fix $w$ and $t\in\Z_{p^{e-1}}$. Let $r\in\{0,1,\dots,p-1\}$ be the
residue of $t$ modulo $p$, and let $z\in\{0,1,\dots,p-1\}$ be the
residue of $w-br$. The elements $x$ with $g(x)\equiv z\pmod p$ are the
elements $F(z+pu)$, and their residues modulo $p^{e-1}$ are the values
$\alpha_z(u)$. Since $\alpha_z$ is a permutation, exactly one of them,
which we call $x_t$, satisfies $x_t\equiv t\pmod{p^{e-1}}$. It satisfies
$g(x_t)+bx_t\equiv z+br\equiv w\pmod p$. Every other element
$x\equiv t\pmod{p^{e-1}}$ satisfies $g(x)\not\equiv z\pmod p$, and hence
$g(x)+bx\not\equiv w\pmod p$. Thus, the elements $x$ with
$g(x)+bx\equiv w\pmod p$ are the elements $x_t$, one in each residue
class modulo $p^{e-1}$. This proves that $\gamma_w$ is a permutation,
and $\gamma_w^{-1}(t)$ is determined by
$g(x_t)+bx_t=w+p\gamma_w^{-1}(t)$.

Put $q_z=\alpha_z^{-1}$. Since $\alpha_z$ maps every residue class onto
a residue class, $q_z$ is congruence-preserving. We have
$g(x_t)=z+pq_z(t)$. As $x_t\equiv t\pmod{p^{e-1}}$, we obtain
$$
\gamma_w^{-1}(t)\equiv q_z(t)+\frac{z+bt-w}p\pmod{p^{e-2}},
$$
where $t$ is represented by an integer. Write $t=r+pa$ with
$a\in\Z_{p^{e-2}}$. Since $q_z$ is congruence-preserving, there is a
congruence-preserving map $Q_r$ of $\Z_{p^{e-2}}$ with
$q_z(r+pa)=q_z(r)+pQ_r(a)$. Hence,
$\gamma_w^{-1}(r+pa)\equiv ba+pQ_r(a)+C_r\pmod{p^{e-2}},$
where $C_r$ does not depend on $a$.

Let $h\in\Z_{p^{e-2}}$. Since $D_{p^{e-1}}(\gamma_w)=H(\gamma_w^{-1})$,
the element $ph$ is a good slope of $\gamma_w$ if and only if
$t\mapsto t+ph\gamma_w^{-1}(t)$ is a permutation. This map preserves
residues modulo $p$, and it maps $r+pa$ to
$r+p\bigl(a+h\gamma_w^{-1}(r+pa)\bigr)$. Hence, it is a permutation if
and only if, for every $r$, the map
$$
a\longmapsto(1+bh)a+hC_r+phQ_r(a)
$$
is a permutation of $\Z_{p^{e-2}}$. Suppose that $1+bh$ is a unit, and
let $a\neq a'$ with $s=\nu(a-a')$. Then, $(1+bh)(a-a')$ has valuation
$s$, while $ph\bigl(Q_r(a)-Q_r(a')\bigr)$ is divisible by $p^{s+1}$.
Hence, the map is injective. If $p$ divides $1+bh$, then the map is
constant modulo $p$, and it is not a permutation, because $e-2\ge1$.
Thus, $ph$ is a good slope of $\gamma_w$ if and only if
$1+\bar b\bar h\neq0$. Exactly $(p-1)p^{e-3}$ elements $h$ of
$\Z_{p^{e-2}}$ satisfy this condition.
\end{proof}

We now show that a regular residue with a large count controls the maps
of Lemma~\ref{lem:cubefibres} formed at the other residues.

\begin{lemma}\label{lem:maxregular}
Let $e\ge3$, and let $f$ be a permutation of $\Z_{p^e}$ that loses
congruences at every level $1,\dots,e-1$. Suppose that $w_f$ is constant
and that some regular residue $c_1$ satisfies $d_{c_1}(f)>(p-1)p^{e-2}/2$.
Let $c\neq c_1$ satisfy $d_c(f)>0$, and choose $m_1$ with
$\bar m_1=c_1$ and $m_1+p^2\Z_{p^e}\subseteq D_{p^e}(f)$. Then, the maps
of Lemma~\ref{lem:cubefibres} formed at $m_1$ are congruence-preserving
permutations. Moreover, the maps of that lemma formed at any good slope
with residue $c$ are permutations of $\Z_{p^{e-1}}$, and at least one of
them is not congruence-preserving.
\end{lemma}

\begin{proof}
Put $g(x)=f(x)+m_1x$ and $F=g^{-1}$. The maps
$\gamma'_z(u)=F(z+pu)\bmod p^{e-1}$ are the maps of
Lemma~\ref{lem:cubefibres} formed at $m_1$. By Lemma~\ref{lem:regular}(i),
they are congruence-preserving, and the set
$H=\bigcap_zH(\gamma'_z)$ is a union of residue classes modulo $p$. Let
$C\subseteq\F_p$ be the set of residues of the elements of $H$. Then,
$0\in C$ and $|C|=d_{c_1}(f)/p^{e-2}>(p-1)/2$. Let $\sigma$ be a digit
map of some $\gamma'_z$. By Lemma~\ref{lem:Hmaps}(iii), the map
$t\mapsto t+h\sigma(t)$ is a permutation of $\F_p$ for every $h\in C$.
For $h\neq0$, this map is a permutation exactly when
$h^{-1}\in D_{\F_p}(\sigma)$. Hence, $|D_{\F_p}(\sigma)|\ge|C|-1>(p-3)/2$,
and Theorem~\ref{thm:fields}(i) shows that $\sigma$ is affine. Thus, every
$\gamma'_z$ has affine digits. At this point, some multipliers may be
zero. Let $\psi_z$ be the map induced by $\gamma'_z$ modulo $p$, which is
its digit map for $j=0$. The counting identity in the proof of
Lemma~\ref{lem:regularpattern} shows that every fibre of $\psi_z$ has
$w_f(r,z-c_1r)/p^{e-2}=1$ element. Thus, $\psi_z$ is a permutation, and
its multiplier is nonzero.

Now fix a good slope $m_0$ with $\bar m_0=c$, and put $b=m_0-m_1$. Then,
$b$ is a unit. The map $T(y)=y+bF(y)$ satisfies
$T=(g+b\,\mathrm{id})\circ F$. Both $F$ and
$g+b\,\mathrm{id}=f+m_0\,\mathrm{id}$ are permutations. Hence, $T$ is a
permutation. We show by
induction on $j$ that every digit map $\sigma_{j,u}$ of every
$\gamma'_z$ has a nonzero multiplier. This holds for $j=0$. Let
$1\le j\le e-2$, and suppose that it holds for the digits
$0,\dots,j-1$.

Let $B$ be a residue class of $\Z_{p^e}$ modulo $p^{j+1}$, and write its
elements as $y=z+pu$. Then, $z$ and the residue of $u$ modulo $p^j$ are
fixed on $B$. Hence, $F(y)$ modulo $p^j$ is constant on $B$. Digit $j$
of $F(y)$ is digit $j$ of $\gamma'_z(u)$. Since $\gamma'_z$ has affine
digits, it is an affine function of digit $j$ of $u$. If its multiplier
is $0$, then $F(y)$ modulo $p^{j+1}$ is constant on $B$. As $F$ is
injective and $B$ has $p^{e-j-1}$ elements, $F$ maps $B$ onto a whole
residue class modulo $p^{j+1}$. We then say that $B$ is of the
\emph{constant type}. If the multiplier is nonzero, then $F(B)$ meets
each of the $p$ residue classes modulo $p^{j+1}$ above the residue of
$F(y)$ modulo $p^j$. We then say that $B$ is of the \emph{distributed
type}. The element $y$ is constant modulo $p^{j+1}$ on $B$, and $b$ is a
unit. Therefore, if $B$ is of the constant type, then $T$ maps $B$ onto
a whole residue class modulo $p^{j+1}$. If $B$ is of the distributed
type, then $T(B)$ meets all $p$ classes above the residue of $T(y)$
modulo $p^j$.

Let $B$ and $B'$ be two such classes with the same residue of $F(y)$
modulo $p^j$. If $B$ were of the constant type and $B'$ of the
distributed type, then $F(B')$ would meet the residue class $F(B)$,
contrary to the injectivity of $F$. Thus, $B$ and $B'$ are of the same
type. In the same way, classes with the same residue of $T(y)$ modulo
$p^j$ are of the same type.

For fixed $z$, the first $j$ digits of $\gamma'_z$ define a
congruence-preserving map of $\Z_{p^j}$ with affine digits. By the
induction hypothesis, its digit maps have nonzero multipliers. Hence, it
is a permutation by Lemma~\ref{lem:digits}(i). By
Lemma~\ref{lem:affine-digits}(ii), its inverse $\theta_z$ also has affine
digits. The classes $B$ correspond bijectively
to the pairs $(z,U)$, where $U\in\Z_{p^j}$ is the residue of $F(y)$
modulo $p^j$. Then, $\theta_z(U)$ is the residue of $u$ modulo $p^j$,
and the residue of $T(y)$ modulo $p^j$ is
$$
W=bU+z+p\theta_z(U).
$$
Consider the bipartite graph of Lemma~\ref{lem:connected}, with $r=j$.
Each class $B$ is an edge joining $U$ to $W$. Two edges with a common
endpoint belong to classes of the same type. As the graph is connected,
all classes $B$ have the same type. If they were all of the constant
type, $F$ would preserve congruences modulo $p^{j+1}$. Then, so would
$g=F^{-1}$ and $f=g-m_1x$, which is excluded. Hence, every digit map at
digit $j$ has a nonzero multiplier. This completes the induction. In
particular, every $\gamma'_z$ is a permutation of $\Z_{p^{e-1}}$, by
Lemma~\ref{lem:digits}(i).

Since the maps $\gamma'_z$ are congruence-preserving permutations,
Lemma~\ref{lem:zerotransfer} applies to $g$ and $b$, with
$\alpha_z=\gamma'_z$. Since $g(x)+bx=f(x)+m_0x$, the maps $\gamma_w$ of
that lemma are the maps $\gamma_z$ of Lemma~\ref{lem:cubefibres} formed
at $m_0$. Hence, these maps are permutations.
By Lemma~\ref{lem:tworegular}, the residue $c$ is not regular, since
$f$ does not preserve congruences modulo $p^2$. By
Lemma~\ref{lem:regular}(ii), at least one $\gamma_z$ is not
congruence-preserving.
\end{proof}

\begin{corollary}\label{cor:regulartransfer}
Let $e\ge3$, put $P=p^{e-2}$, and suppose that every permutation $\pi$
of $\Z_{p^{e-1}}$ that is not congruence-preserving satisfies
$N(\pi)\le M$. Let $f$ be a permutation of $\Z_{p^e}$ that loses
congruences at every level, and suppose that $w_f$ is constant. Suppose
also that a regular residue $c_1$ satisfies $d_{c_1}(f)>(p-1)P/2$. Then,
every residue $c\neq c_1$ with $d_c(f)>0$ satisfies
$
d_c(f)\le M-\frac{(p-1)P}p.
$
\end{corollary}

\begin{proof}
Choose $m_1$ with $\bar m_1=c_1$ and
$m_1+p^2\Z_{p^e}\subseteq D_{p^e}(f)$, and put $g(x)=f(x)+m_1x$. Choose
a good slope $m_0$ with residue $c$, and put $b=m_0-m_1$. Then, $b$ is a
unit, and $g(x)+bx=f(x)+m_0x$ is a permutation. By
Lemma~\ref{lem:maxregular}, the maps of Lemma~\ref{lem:cubefibres}
formed at $m_1$ are congruence-preserving permutations. They are the maps
$\alpha_z$ of Lemma~\ref{lem:zerotransfer} for $g$. The maps $\gamma_w$
of Lemma~\ref{lem:zerotransfer} are the maps of
Lemma~\ref{lem:cubefibres} formed at $m_0$. By
Lemma~\ref{lem:maxregular}, one of them is not congruence-preserving.
Let $\pi$ be its inverse. Then, $\pi$ is not congruence-preserving, and
Lemma~\ref{lem:zerotransfer} gives
$d_0(\pi^{-1})=d_0(\gamma_w)=(p-1)P/p$. Since
$H(\gamma_w)=D_{p^{e-1}}(\pi)$, we obtain
$|H(\gamma_w)|=N(\pi)-d_0(\pi^{-1})\le M-\frac{(p-1)P}p.$
Finally, $d_c(f)\le|H(\gamma_w)|$ by Lemma~\ref{lem:cubefibres}.
\end{proof}

We can now prove the main reduction lemma for all odd primes.

\begin{lemma}\label{lem:strongpropagate}
Let $e\ge3$, and put $P=p^{e-2}$. Suppose that every permutation $\pi$
of $\Z_{p^{e-1}}$ that is not congruence-preserving satisfies
$N(\pi)\le M$, where
$$
\max\Bigl\{P-1,\ \frac{(p-1)P}2,\ \frac{2(p-1)P}p\Bigr\}\le M\le(p-1)P.
$$
If $f$ is a permutation of $\Z_{p^e}$ that loses congruences at every
level, then
$|D_{p^e}(f)|\le(p-1)M$ and $N(f)\le pM.$
\end{lemma}

\begin{proof}
Since $M\ge P-1$, Lemma~\ref{lem:Hnormal} shows that every map $\gamma$
of $\Z_{p^{e-1}}$ that is not congruence-preserving satisfies
$|H(\gamma)|\le M$. Hence, by Lemma~\ref{lem:regular}, every direction
class that is not regular has count at most $M$, and every regular class
has count at most $(p-1)P$. The permutation $f^{-1}$ also loses
congruences at every level, and the count of the vertical class of $f$
is $d_0(f^{-1})$. Hence, the same bounds hold for the vertical class. By
Theorem~\ref{thm:universal}, at most $p-1$ finite classes and at most
$p$ classes altogether have a nonzero count.

Suppose first that $w_f$ is not constant. By
Lemma~\ref{lem:regularpattern}, applied to $f$ and to $f^{-1}$, every
regular count is at most $(p-1)P/2\le M$. Hence, every count is at most
$M$, and we obtain $|D_{p^e}(f)|\le(p-1)M$ and $N(f)\le pM$.

From now on, let $w_f$ be constant, and put $K=(p-1)P/p$. By
Lemma~\ref{lem:tworegular}, at most one finite class is regular. If no
finite class is regular with a count greater than $(p-1)P/2$, then every
finite count is at most $M$, and $|D_{p^e}(f)|\le(p-1)M$. Otherwise,
Corollary~\ref{cor:regulartransfer} bounds every other finite count by
$M-K$. Hence,
$$
|D_{p^e}(f)|\le(p-1)P+(p-2)(M-K)
=(p-1)M-\Bigl(M-\frac{2(p-1)P}p\Bigr)\le(p-1)M.
$$

For the bound on $N(f)$, suppose first that at most two classes have a
nonzero count. Then, $N(f)\le2(p-1)P\le pM$. Suppose next that at least
three classes have a nonzero count. Then, at most one of them is regular.
Otherwise, take a good form from a third class with nonzero count as the
first coordinate and a good form from one of two regular classes as the
second coordinate. Lemma~\ref{lem:projective}(ii) gives a permutation
with constant residue pattern that loses congruences at every level and
has two regular finite classes. This contradicts
Lemma~\ref{lem:tworegular}. If no regular class has a count greater
than $(p-1)P/2$, then every count is at most $M$, and $N(f)\le pM$.
Otherwise, let $\delta$ be a class with nonzero count other than the
regular class. Take a good form from a third class with nonzero count as
the first coordinate and a good form from the regular class as the second
coordinate. In the resulting permutation, the regular class is the finite
class of $0$, and $\delta$ is another finite class. By
Corollary~\ref{cor:regulartransfer}, the count of $\delta$ is at most
$M-K$. At most $p-1$ classes other than the regular one have a nonzero
count. Hence,
$N(f)\le(p-1)P+(p-1)(M-K)=(p-1)M+K\le pM,$
because $K\le(p-1)P/2\le M$.
\end{proof}

\section{Moment identities and sharper bounds}\label{sec:sharpsquares}

We sharpen the bounds for $\Z_{p^2}$ by summing the values of a good
linear form over the points at which it has a fixed residue modulo $p$.
The estimates obtained here
complete the base case for Theorem~\ref{thm:allpowers}.
Throughout this section, $p$ is an odd prime.

Power sums over finite fields are a classical tool in the theory of
permutation polynomials, for instance through Hermite's criterion
\cite{LN}. For a map $h:\F_p\to\F_p$, an integer $i\ge1$ and
$0\le j\le p-2$, the mixed sum $\sum_{x\in\F_p}x^jh(x)^i$ is, up to
sign, the coefficient of $x^{p-1-j}$ in the polynomial of degree less
than $p$ that represents $h^i$. Ball and G\'acs use bounds on the
degrees of these polynomials to study the graph of $h$ \cite{BG}. The
sums $\Lambda_k$ in the proof of Lemma~\ref{lem:tworigid} are sums of
the same kind over $\F_p^2$.

We use the language of Section~\ref{sec:cube} for arbitrary sets of
points. Let $e\ge2$, and let $\Gamma\subseteq\Z_{p^e}^2$ be a set of
$p^e$ points. A class $[a:b]$ of primitive pairs is \emph{good} on
$\Gamma$ if the linear form $aX+bY$ maps $\Gamma$ bijectively onto
$\Z_{p^e}$. For a direction $\delta$ of $\F_p^2$, the classes $[a:b]$
for which $\bar aX+\bar bY$ is constant on the lines with direction
$\delta$ form a direction class. The \emph{count} of $\delta$ is the
number of good classes in it. A direction class is \emph{regular} if it
contains a whole fibre of the reduction modulo $p^2$ consisting of good
classes. If $\Gamma$ is the graph of a permutation $\pi$ of $\Z_{p^e}$,
the count of $\langle(1,-c)\rangle$ is $d_c(\pi)$, the count of the
vertical direction is $d_0(\pi^{-1})$, and the sum of all counts is
$N(\pi)$. Let $A\in\mathrm{GL}_2(\Z_{p^e})$. As in the proof of
Lemma~\ref{lem:projective}(ii), the form $aX+bY$ is good on $\Gamma$ if
and only if the form with coefficient vector $(a,b)A^{-1}$ is good on
$A\Gamma$. Hence, the reduction of $A$ maps the directions of $\Gamma$
onto the directions of $A\Gamma$, and it preserves counts and
regularity. A translation of $\Gamma$ adds a constant to every form,
and it also preserves these notions. If $e=2$ and $\Gamma$ is a lifted
plane, then so is $A\Gamma$, and the rigid directions of $\Gamma$ are
mapped onto those of $A\Gamma$, as in the proof of
Lemma~\ref{lem:tworigid}. For $i,j\in\F_p$, we write $C_{i,j}$ for the
set of points of $\Gamma$ that reduce to $(i,j)$.

The next lemma computes the counts of two rigid directions of a lifted
plane, and it restricts the other directions with a nonzero count.

\begin{lemma}\label{lem:rigidmoments}
Let $p\ge5$ be prime, and let $\Gamma$ be a lifted plane for which the
directions $\langle(0,1)\rangle$ and $\langle(1,0)\rangle$ are rigid.
Write the points of $\Gamma$ as in Lemma~\ref{lem:normalform}, and put
$T=\sum_it_i$ and $S=\sum_js_j$ in $\F_p$.
\begin{enumerate}
\item The count of $\langle(0,1)\rangle$ is $p-|\{t_i:i\in\F_p\}|$, and
the count of $\langle(1,0)\rangle$ is $p-|\{s_j:j\in\F_p\}|$.
\item Let $c\in\F_p^\times$. If the direction $\langle(1,-c)\rangle$ has
a nonzero count, then $Tc^2-c+S=0$.
\item At most two directions other than $\langle(0,1)\rangle$ and
$\langle(1,0)\rangle$ have a nonzero count. If one of the two counts in
(i) equals $p-1$, at most one such direction has a nonzero count. If
both counts in (i) equal $p-1$, no such direction has a nonzero count.
\end{enumerate}
\end{lemma}

\begin{proof}
The direction class of $\langle(0,1)\rangle$ consists of the classes of
the forms $X+p\lambda Y$ with $\lambda\in\F_p$. At $P_{ij}$, such a form
takes the value $\tau(i)+p\bigl(t_ij+\nu_i+\lambda j\bigr)$. Its residue
modulo $p$ is $i$. For fixed $i$, the coefficient of $p$ is an affine
function of $j$ with multiplier $t_i+\lambda$. Hence, the form is good
if and only if $\lambda\neq-t_i$ for every $i$. This gives the first
count. The second count is obtained in the same way from the forms
$Y+p\lambda X$, whose values are $\tau(j)+p(s_ji+\mu_j+\lambda i)$.

For (ii), the direction class of $\langle(1,-c)\rangle$ consists of the
classes of the forms $Y+(\tau(c)+ph)X$ with $h\in\F_p$. Suppose that
such a form is good, and fix $z\in\F_p$. The points at which the form
has residue $z$ are the points $P_{ij}$ with $j=z-ci$. Since $\tau$ is
multiplicative, Lemma~\ref{lem:carry} shows that the value of the form
at such a point is
$\tau(z)+p\bigl(E(j,ci)+s_ji+\mu_j+ct_ij+c\nu_i+hi\bigr).$
These $p$ values are the $p$ elements of $\Z_{p^2}$ with residue $z$.
Hence, the coefficients of $p$ run over $\F_p$, and their sum is zero.
As $i$ runs over $\F_p$, so do $z-ci$ and $ci$. By
Lemma~\ref{lem:carry}, we therefore have
$$
p\sum_iE(z-ci,ci)=\sum_i\bigl(\tau(z-ci)+\tau(ci)-\tau(z)\bigr)
=-p\tau(z)
$$
in $\Z_{p^2}$, and hence $\sum_iE(z-ci,ci)=-z$ in $\F_p$. We also
have $\sum_ii=0$. Writing $i=(z-j)/c$ in the terms involving $s_j$ and
$\mu_j$, we obtain
$$
0=-z+\frac{zS}c-\frac1c\sum_jjs_j+\sum_j\mu_j+czT
-c^2\sum_iit_i+c\sum_i\nu_i.
$$
Both sides are affine functions of $z$, and the identity holds for
every $z\in\F_p$. Comparing the coefficients of $z$ gives
$-1+S/c+cT=0$. Multiplication by $c$ proves (ii).

For (iii), note that the directions other than $\langle(0,1)\rangle$
and $\langle(1,0)\rangle$ are the directions $\langle(1,-c)\rangle$ with
$c\neq0$. The coefficient of $c$ in $Tc^2-c+S$ is $-1$. Hence, the
equation $Tc^2-c+S=0$ has at most two solutions $c\in\F_p$. If a count
in (i) equals $p-1$, then the $t_i$ or the $s_j$ are all equal. In the
first case, $T=0$, and the only solution is $c=S$. In the second case,
$S=0$, and the equation $c(Tc-1)=0$ has at most one nonzero solution.
If both counts equal $p-1$, then $T=S=0$, and $c=0$ is the only
solution. Now (ii) proves (iii).
\end{proof}

For spread permutations, Lemma~\ref{lem:rigidmoments} improves the bound
of Theorem~\ref{thm:p2}.

\begin{proposition}\label{prop:spreadimproved}
Let $p\ge5$ be prime, and let $f$ be a spread permutation of
$\Z_{p^2}$. Then,
$$
N(f)\le\frac{(p-1)^2}2,\quad |D_{p^2}(f)|\le\frac{(p-1)^2}2-1.
$$
Moreover, every permutation of $\Z_{p^2}$ that is not
congruence-preserving has at most $(p-1)^2/2$ good slopes.
\end{proposition}

\begin{proof}
Let $\Gamma$ be the graph of $f$, and put $b=(p-3)/2$. As in the proof
of Proposition~\ref{prop:directions}, $N(f)$ is the sum of the counts of
the $p+1$ directions of $\F_p^2$. Every count is at most $p-1$, and the
count of a direction that is not rigid is at most $b$. By
Lemma~\ref{lem:projective}(i), some count vanishes. By
Lemma~\ref{lem:tworigid}, at most two directions are rigid.

If at most one direction is rigid, then
$$
N(f)\le(p-1)+(p-1)b=\frac{(p-1)^2}2.
$$
Suppose that two directions are rigid. Choose
$A\in\mathrm{GL}_2(\Z_{p^2})$ whose reduction maps them to
$\langle(0,1)\rangle$ and $\langle(1,0)\rangle$. Then, $A\Gamma$ is a
lifted plane with these two rigid directions, and its counts are those
of $\Gamma$, in a different order. We apply
Lemma~\ref{lem:rigidmoments} to $A\Gamma$. If both rigid directions have
count $p-1$, all other counts vanish, and $N(f)=2p-2$. If exactly one
of them has count $p-1$, their counts sum to at most $2p-3$, and at most
one other direction has a nonzero count. This count is at most $b$.
Hence, $N(f)\le2p-3+b=(5p-9)/2$. If neither of them has count $p-1$,
their counts sum to at most $2p-4$, and at most two other directions
have a nonzero count. Hence, $N(f)\le2p-4+2b=3p-7$. The differences
between $(p-1)^2/2$ and these three bounds are $(p-1)(p-5)/2$,
$(p-2)(p-5)/2$ and $(p-3)(p-5)/2$, respectively. Since $p\ge5$, all
three bounds are at most $(p-1)^2/2$. This proves the bound on
$N(f)$.

Since $0$ is a good slope of $f^{-1}$, we have $d_0(f^{-1})\ge1$.
Hence, $|D_{p^2}(f)|=N(f)-d_0(f^{-1})\le(p-1)^2/2-1$. Finally, a
permutation that is neither congruence-preserving nor spread has at
most $(p-1)^2/2$ good slopes by Corollary~\ref{cor:notspread}.
\end{proof}

\section{Bounds for all odd prime powers}\label{sec:higher}

We now prove Theorem~\ref{thm:allpowers}. The estimates for $p^2$
come from the preceding sections. The reduction through a preserved
congruence and Lemma~\ref{lem:strongpropagate} handle larger exponents.

\begin{proof}[Proof of Theorem~\ref{thm:allpowers}]
For $p=3$, put $U_e=3^{e-1}$ and $V_e=4\cdot3^{e-2}$. For $p\ge5$,
put $U_e=(p-1)^2p^{e-2}/2$ and $V_e=p(p-1)p^{e-2}/2$. We show by
induction on $e$ that $|D_{p^e}(f)|\le U_e$ and $N(f)\le V_e$. Let
$e=2$. The first bound is Theorem~\ref{thm:p2} if $p=3$, and
Proposition~\ref{prop:spreadimproved} if $p\ge5$. The second bound is
Proposition~\ref{prop:directions}.

Let $e\ge3$. Suppose first that $f$ preserves congruences modulo
$p^i$ for some $i$ with $1\le i<e$. By Lemma~\ref{lem:split}(iii), the
quotient permutation or one of the residue-class permutations of $f$ is
not congruence-preserving. We denote it by $h$, and we let $\Z_{p^j}$
be the ring on which it acts. Then, $2\le j<e$, because every
permutation of $\Z_p$ is congruence-preserving. By
Lemmas~\ref{lem:split}(iii) and~\ref{lem:projective}(iii) and the
induction hypothesis,
$$
|D_{p^e}(f)|\le p^{e-j}U_j=U_e,\quad
N(f)\le p^{e-j}V_j=V_e.
$$

Hence, we may assume that $f$ loses congruences at every level. Put
$P=p^{e-2}$ and $M=V_{e-1}$. By the induction hypothesis, every
permutation $\pi$ of $\Z_{p^{e-1}}$ that is not congruence-preserving
satisfies $N(\pi)\le M$. For $p\ge5$, we have $M=(p-1)P/2$. Then, the
hypotheses of Lemma~\ref{lem:strongpropagate} hold, because
$P-1<(p-1)P/2$ and $2(p-1)P/p\le(p-1)P/2$. The lemma gives
$|D_{p^e}(f)|\le(p-1)M=U_e$ and $N(f)\le pM=V_e$. For $p=3$, we have
$M=4P/3$, and the hypotheses of Lemma~\ref{lem:strongpropagate} also
hold. The lemma gives $|D_{3^e}(f)|\le8P/3<3P=U_e$ and
$N(f)\le4P=V_e$. This completes the induction. Finally,
$U_e=\varphi(3^e)/2$ for $p=3$, and $U_e<\varphi(p^e)/2$ for $p\ge5$.

For sharpness when $p=3$, let $f$ be the permutation of
Example~\ref{ex:Z9spread} if $e=2$, and the permutation of
Example~\ref{ex:lift} if $e\ge3$. Then, $|D_{3^e}(f)|=3^{e-1}$. At
$e=2$, the permutation is an involution with $D_9(f)=\{0,2,5\}$, and
hence $N(f)=3+1=4$. For $e\ge3$, the permutation $f^{-1}$ preserves
congruences modulo $9$. Its quotient permutation is the inverse of the
permutation of Example~\ref{ex:Z9spread}, and its residue-class
permutations are $t\mapsto2^{-1}t$. By Lemma~\ref{lem:split-general},
an element $m\in3\Z_{3^e}$ is a good slope of $f^{-1}$ if and only if
$m\equiv0\pmod9$. Hence, $d_0(f^{-1})=3^{e-2}$ and
$N(f)=3^{e-1}+3^{e-2}=4\cdot3^{e-2}$.
\end{proof}

\begin{example}\label{ex:81}
For $n=81$, every permutation with more than $27$ good slopes has
exactly $54$. Write $x=x_0+3x_1+9x_2+27x_3$ and
$f(x)=y_0+3y_1+9y_2+27y_3$ in base $3$.
By Proposition~\ref{thm:classification}, these are exactly the
permutations given by
$$
\begin{aligned}
y_0&=\alpha x_0+\beta,& y_1&=\alpha x_1+h_1(x_0),\\
y_2&=\alpha x_2+h_2(x_0,x_1),&
y_3&=\alpha x_3+h_3(x_0,x_1,x_2),
\end{aligned}
$$
where the equations are in $\F_3$, $\alpha\in\{1,2\}$,
$\beta\in\F_3$, and the maps $h_j$ are arbitrary.
There are $2\cdot3^{40}$ such permutations.
At the boundary, the lift in Example~\ref{ex:lift} has
$$
D_{81}(f)=\{m\in\Z_{81}:m\bmod9\in\{0,2,5\}\}.
$$
It does not preserve congruences modulo $3$ and has $27$ good slopes.
The set $\{0,2\}$ is admissible, and $\mu_{81}(f)=2$. This shows
that the threshold $3^{e-1}$ in Corollary~\ref{cor:threehalf} cannot be
lowered.
\end{example}

\section{Sharp bounds for powers of \texorpdfstring{$5$}{5}}\label{sec:five}

In this section, we determine the largest number of good slopes of a
permutation of $\Z_{5^e}$ that is not congruence-preserving for
$e=2,3,4$, and we bound it for larger $e$. For $p=5$,
Proposition~\ref{prop:spreadimproved} bounds the number of good slopes
of a spread permutation by seven, and Example~\ref{ex:Z25seven} attains
this bound. We first treat the permutations of $\Z_{25}$ that are
neither congruence-preserving nor spread.

We also use the following counting principle from \cite{Extremal}. For
a prime $r$ and a map $w:\F_r^2\to\Z$, the \emph{weight} of a line of
$\F_r^2$ is the sum of the values of $w$ on its points.

\begin{lemma}[{\cite[Lemma~6.2]{Extremal}}]\label{lem:uniform}
Let $r$ be a prime.
\begin{enumerate}
\item Let $w:\F_r^2\to\Z$ be a map and $L\in\Z$. Suppose that, for $r$
of the $r+1$ directions of $\F_r^2$, every line with that direction has
weight $L$. Then, $w$ is constant on every line with the
remaining direction.
\item Let $T$ be a finite set of size $L$, and let $A:\F_r\times T\to\F_r$
be a map. Suppose that, for every $a\in\F_r^\times$ and $c\in\F_r$,
exactly $L$ pairs $(s,t)$ satisfy $A(s,t)+as=c$. Then, for each
$z\in\F_r$, the number $|\{t\in T:A(s,t)=z\}|$ is independent of $s$.
\end{enumerate}
\end{lemma}

A set $S\subseteq\F_p^2$ is \emph{centrally symmetric} about a point
$P$ if $Q\in S$ implies $2P-Q\in S$.

\begin{lemma}\label{lem:fivepoints}
Every permutation of $\F_5$ is affine or has the form
$t\mapsto A(t-b)^3+d$, where $A\neq0$. Consequently, every set of five
points of $\F_5^2$ that is injective in two distinct directions is
centrally symmetric about one of its points. If such a set is not a
line, then, after an affine change of the first coordinate, the
multiplicities $h(i)$ of the first coordinates of its points are
$h(i)=1$ for every $i$, or
$$
h(0)=3,\ h(1)=h(-1)=1\quad\text{or}\quad h(0)=1,\ h(1)=h(-1)=2.
$$
In the last two cases, $h(2)=h(-2)=0$.
\end{lemma}

\begin{proof}
Let $\pi$ be a permutation of $\F_5$, and represent it by a polynomial
of degree at most $4$. Over $\F_5$, the power sum $\sum_tt^k$ equals
$-1$ if $k\ge1$ and $4\mid k$, and it equals $0$ otherwise. Since
$\sum_t\pi(t)=\sum_tt=0$, the coefficient of $t^4$ vanishes. A
polynomial of degree $2$ is not injective in odd characteristic. Thus,
if $\pi$ is not affine, it has degree $3$. After a translation of the
variable, $\pi(t)=At^3+Ct+D$ with $A\neq0$. Since $\pi$ is a
permutation, $\sum_t\pi(t)^2=\sum_tt^2=0$. Among the exponents
$0,\dots,6$, only $4$ has a nonzero power sum. Hence, expanding
$\pi(t)^2$ gives $\sum_t\pi(t)^2=-2AC$, and therefore $C=0$.

Let $S$ be a set of five points that is injective in two distinct
directions $\delta_1$ and $\delta_2$. Choose linear forms $L_1$ and
$L_2$ on $\F_5^2$ that are constant on the lines with directions
$\delta_1$ and $\delta_2$, respectively. They are linearly independent,
and each of them maps $S$ bijectively onto $\F_5$. Hence, in the
coordinates $(L_1,L_2)$, the set $S$ is the graph of a permutation of
$\F_5$. An affine graph is centrally symmetric about each of its points,
and the graph of $A(t-b)^3+d$ is centrally symmetric about $(b,d)$.
Since an invertible linear change of coordinates preserves central
symmetry, $S$ is centrally symmetric about one of its points
$(i_0,v_0)$.

The reflection $Q\mapsto2(i_0,v_0)-Q$ fixes $(i_0,v_0)$ and pairs the
other four points of $S$. Hence, the multiplicity of $i_0$ among the
first coordinates is odd, and the multiplicities of $i_0+i$ and $i_0-i$
are equal. If it is $5$, then $S$ is a vertical line. If it is $3$, the
remaining two points have first coordinates $i_0\pm i$ for some
$i\neq0$. If it is $1$, the remaining four points have first
coordinates $i_0\pm i$ and $i_0\pm i'$, where either $i'=\pm i$ or the
four coordinates are distinct. The affine change $X\mapsto(X-i_0)/i$
gives the stated forms.
\end{proof}

We next prove a nonexistence result. We state it for every prime
$p\ge5$ and every exponent, because Lemma~\ref{lem:fivepattern} uses it
for higher powers of $5$.

\begin{lemma}\label{lem:threelines}
Let $p\ge5$ be prime, let $e\ge2$, put $P=p^{e-2}$, and let
$1\le k\le(p-1)/2$. Define $h:\F_p\to\Z$ by $h(0)=p-2k$,
$h(1)=h(-1)=k$, and $h(i)=0$ for $i\notin\{0,1,-1\}$. Let $\tau(i)$ be
integer representatives of the elements of Lemma~\ref{lem:carry}, with
$\tau(0)=0$, $\tau(1)=1$ and $\tau(-1)=-1$. Then, there is no set
$\Gamma\subseteq\Z_{p^e}^2$ with the following properties.
\begin{enumerate}
\item For all $i,j\in\F_p$, exactly $Ph(i)$ points of $\Gamma$ reduce to
$(i,j)$.
\item For every $c\in\F_p$, some form $Y+(\tau(c)+pb)X$ with
$b\in\Z_{p^{e-1}}$ is good on $\Gamma$.
\item For all $i,j\in\F_p$, both sums
$\sum_{(X,Y)\in C_{i,j}}(X-\tau(i))/p$ and
$\sum_{(X,Y)\in C_{i,j}}(Y-\tau(j))/p$ lie in $P\Z_{p^{e-1}}$.
\end{enumerate}
\end{lemma}

For $e=2$, we have $P=1$, and condition (iii) holds automatically.

\begin{proof}
Suppose that such a set $\Gamma$ exists. By (i), it has $p^e$ points.
Multiplication by $P$ is an isomorphism of additive groups
from $\F_p$ onto $P\Z_{p^{e-1}}$. Hence, by (iii), we may define $A_i(j),B_i(j)\in\F_p$
for $i\in\{0,1,-1\}$ by
$$
P\,h(i)A_i(j)=\sum_{(X,Y)\in C_{i,j}}\frac{X-\tau(i)}p,\quad
P\,h(i)B_i(j)=\sum_{(X,Y)\in C_{i,j}}\frac{Y-\tau(j)}p.
$$
This is possible, because $0<h(i)<p$.

Let $c,z\in\F_p$, and choose a good form $L=Y+(\tau(c)+pb)X$ as in (ii).
Write a point of $C_{i,j}$ as $(\tau(i)+pU,\tau(j)+pV)$, where
$U,V\in\Z_{p^{e-1}}$. If $j=z-ci$, then
$$
L=\tau(z)+p\bigl(V+\tau(c)U+b\tau(i)+pbU+T_{i,j}\bigr),\quad
T_{i,j}=\frac{\tau(j)+\tau(c)\tau(i)-\tau(z)}p.
$$
The number $T_{i,j}$ is an integer, and $T_{i,j}\equiv E(j,ci)\pmod p$
by Lemma~\ref{lem:carry}. Since $L$ is good, its values with residue $z$
are the $p^{e-1}$ elements $\tau(z)+pt$ with $t\in\Z_{p^{e-1}}$. The sum
of these elements $t$ is zero in $\Z_{p^{e-1}}$, because $p$ is odd.

We sum the coefficient of $p$ in the formula for $L$ over the points of
the sets $C_{i,z-ci}$. By (iii), every sum of $U$ over a set $C_{i,j}$
lies in $P\Z_{p^{e-1}}$. Hence, the terms $pbU$ contribute $0$, and
$\tau(c)$ may be replaced by $c$ in the terms $\tau(c)U$. The terms
$b\tau(i)$ contribute $bP\sum_ih(i)\tau(i)=0$, by the choice of
$\tau(0)$ and $\tau(\pm1)$. The terms $T_{i,j}$ contribute
$P\sum_ih(i)T_{i,z-ci}$. Dividing the resulting identity by $P$, we
obtain the identity
$$
\sum_{i\in\{0,1,-1\}}h(i)\bigl(B_i(z-ci)+cA_i(z-ci)+E(z-ci,ci)\bigr)=0
$$
in $\F_p$. We have $h(0)\equiv-2k\pmod p$ and $E(z,0)=0$. Dividing by
$k$, we obtain
$$
\begin{aligned}
&B_1(z-c)+B_{-1}(z+c)+c\bigl(A_1(z-c)+A_{-1}(z+c)\bigr)
+E(z-c,c)+E(z+c,-c)\\
&\quad=2B_0(z)+2cA_0(z).
\end{aligned}
$$
This identity holds for all $c,z\in\F_p$.

We substitute $u=z-c$ and $v=z+c$. This is an invertible change of
variables, because $p$ is odd. Define $\kappa\in\F_p$ by
$\tau(2)\equiv2+p\kappa\pmod{p^2}$. By Lemma~\ref{lem:carry} and
$\tau(c)+\tau(-c)\equiv0\pmod{p^2}$, we have
$p\bigl(E(z-c,c)+E(z+c,-c)\bigr)\equiv\tau(u)+\tau(v)-2\tau(z)
\equiv\tau(2z)-2\tau(z)+pE(u,v)\pmod{p^2}.$
Moreover, $\tau(2z)\equiv\tau(2)\tau(z)\equiv2\tau(z)+p\kappa z
\pmod{p^2}$. Hence, $E(z-c,c)+E(z+c,-c)=E(u,v)+\kappa(u+v)/2$, and we
obtain the identity
$$
\begin{aligned}
0={}&B_1(u)+B_{-1}(v)+\frac{v-u}2\bigl(A_1(u)+A_{-1}(v)\bigr)
+E(u,v)\\
&+\kappa\,\frac{u+v}2-2B_0\Bigl(\frac{u+v}2\Bigr)
-(v-u)\,A_0\Bigl(\frac{u+v}2\Bigr)
\end{aligned}
$$
for all $u,v\in\F_p$.

We represent every function of one variable by a polynomial of degree
at most $p-1$, and we expand the right-hand side as a polynomial in $u$
and $v$. Every exponent is at most $p$, and $u^p$ and $v^p$ occur only
as pure powers. We replace them by $u$ and $v$, respectively. The
result $F(u,v)$ has degree at most $p-1$ in each variable and vanishes
on $\F_p^2$. For fixed $v$, it is a polynomial in $u$ of degree at most
$p-1$ with $p$ roots. Hence, its coefficients, as polynomials in $v$,
vanish on $\F_p$, and $F$ is the zero polynomial.

We compare the coefficients of $u^2v^{p-2}$ and $u^{p-2}v^2$ in $F$.
These monomials are distinct, because $p\ge5$. By the definition of $E$
in Lemma~\ref{lem:carry}, the polynomial $E(u,v)$ contributes $1/2$ to
each of them, since $(-1)^{p-2}/(p-2)=1/2$ in $\F_p$. The term
$-(v-u)A_0((u+v)/2)$ changes sign when $u$ and $v$ are exchanged.
Hence, its contributions to the two coefficients are opposite. The term
$-2B_0((u+v)/2)$ has total degree at most $p-1$, and the terms $B_1(u)$,
$B_{-1}(v)$ and $\kappa(u+v)/2$ have no mixed monomials. Every mixed
monomial of $(v-u)\bigl(A_1(u)+A_{-1}(v)\bigr)$ has an exponent equal to
$1$. Thus, these terms contribute nothing to the two coefficients.
Consequently, the sum of the two coefficients of $F$ is $1$, which
contradicts $F=0$.
\end{proof}

We can now determine the maximum for $\Z_{25}$.

\begin{theorem}\label{thm:exact25}
The largest number of good slopes of a permutation of $\Z_{25}$ that is
not congruence-preserving is $7$. More precisely, every such permutation
$f$ satisfies
$$
N(f)\le8,\quad |D_{25}(f)|\le7,
$$
and $|D_{25}(f)|\le6$ if $f$ is not spread.
\end{theorem}

\begin{proof}
If $f$ is spread, Proposition~\ref{prop:spreadimproved} gives both
bounds. Suppose that $f$ is not spread, and let $\Gamma$ be its graph.
By Corollary~\ref{cor:notspread}, $d_c(f)\le2$ for every $c\in\F_5$. As
in the proof of Proposition~\ref{prop:directions}, $f^{-1}$ is neither
congruence-preserving nor spread, and the same corollary gives
$d_0(f^{-1})\le2$. Thus, every direction of $\F_5^2$ has count at most
$2$. By Lemma~\ref{lem:projective}(i), at most five directions have a
nonzero count.

We claim that, if five directions have a nonzero count, then all these
counts equal $1$. Suppose otherwise. Let $\delta_0$ be the direction
with count $0$, and let $\delta_1$ be a direction with count $2$. For
$Q\in\F_5^2$, let $w(Q)$ be the number of points of $\Gamma$ that reduce
to $Q$. A good form takes every value of $\Z_{25}$ exactly once on
$\Gamma$. Hence, its reduction takes every value of $\F_5$ exactly five
times. Thus, every line with one of the five directions different from
$\delta_0$ has weight $5$. By Lemma~\ref{lem:uniform}(i), the function
$w$ is constant on every line with direction $\delta_0$.

Choose a primitive linear form $X'$ whose reduction is constant on the
lines with direction $\delta_0$, and a good form $Y'$ in the direction
class of $\delta_1$. Their reductions are linearly independent. Hence,
$(X',Y')$ is an invertible linear change of coordinates over $\Z_{25}$.
Let $\Gamma'$ be the image of $\Gamma$. For $i,j\in\F_5$, the number of
points of $\Gamma'$ that reduce to $(i,j)$ is a number $h(i)$ that does
not depend on $j$. In the new coordinates, $\delta_0$ is the vertical
direction, and $\delta_1$ is the direction $\langle(1,0)\rangle$. The
other five directions are $\langle(1,-c)\rangle$ with $c\in\F_5$, and
they all have a nonzero count. Thus, for every $c\in\F_5$, some form
$Y'+(\tau(c)+5b)X'$ is good on $\Gamma'$.

Fix $j\in\F_5$. Since $Y'$ is good, the five points of $\Gamma'$ whose
second coordinate has residue $j$ have second coordinates $\tau(j)+5v$
with five distinct $v\in\F_5$. Let $S_j$ be the set of the pairs
$(i,v)$, where $i$ is the residue of the first coordinate of the point.
The count of $\langle(1,0)\rangle$ is $2$. Hence, the forms $Y'+5\lambda
X'$ are good for two distinct values of $\lambda\in\F_5$. Such a form
takes the value $\tau(j)+5(v+\lambda i)$ at the point corresponding to
$(i,v)$. Thus, $S_j$ is injective in the directions
$\langle(1,-\lambda)\rangle$ for these two values of $\lambda$. The
multiplicity of $i$ among the first coordinates of $S_j$ is $h(i)$.

If $h(i)=1$ for every $i$, then $\Gamma'$ is a lifted plane. Then,
$\Gamma$ is also a lifted plane, and $f$ is spread. This happens in
particular when $S_j$ is a line that is not vertical. If $S_j$ is a
vertical line, then $h$ is supported at a single point $i_0$. In this
case, the reduction of $\Gamma$ lies on one line with direction
$\delta_0$. This line is not vertical, because $\bar x$ takes every
value as $x$ runs over $\Z_{25}$. Hence, $\overline{f(x)}$ is an affine
function of $\bar x$, and $f$ preserves congruences modulo $5$. Both
cases are excluded. Therefore, $S_j$ is not a line, and $h$ is not
constant. By Lemma~\ref{lem:fivepoints}, after an affine change of the
first coordinate, we have
$$
h(0)=3,\ h(1)=h(-1)=1\quad\text{or}\quad h(0)=1,\ h(1)=h(-1)=2.
$$
Let this change of the first coordinate be $X''=s(X'-x_0)$, where
$s\in\Z_{25}^\times$ and $x_0\in\Z_{25}$, and let $\Gamma''$ be the
image of $\Gamma'$ under $(X',Y')\mapsto(X'',Y')$. Then, $\Gamma''$
satisfies condition (i) of Lemma~\ref{lem:threelines} with $p=5$ and
$k=1$ or $k=2$. Let $c\in\F_5$, and choose a good form
$Y'+(\tau(c')+5b)X'$ on $\Gamma'$, where $c'\equiv cs\pmod5$. In the
coordinates $(X'',Y')$, it becomes the form $Y'+(\tau(c')+5b)s^{-1}X''$
plus a constant. This form is good on $\Gamma''$, and its coefficient of
$X''$ is congruent to $c$ modulo $5$. Hence, $\Gamma''$ also satisfies
condition (ii). Condition (iii) holds
automatically, since $e=2$. This contradicts Lemma~\ref{lem:threelines}
and proves the claim.

If some direction has count $2$, at most four directions have a nonzero
count by the claim, and $N(f)\le8$. Otherwise, $N(f)\le5$. Now suppose
that $|D_{25}(f)|\ge7$. Since every finite direction has count
at most $2$, at least four finite directions have a nonzero count, and
one of them has count $2$. The vertical direction also has a nonzero
count, because $0\in D_{25}(f^{-1})$. This contradicts the claim.
Hence, $|D_{25}(f)|\le6$. Finally, the permutation of
Example~\ref{ex:Z25seven} is not congruence-preserving and has seven good
slopes.
\end{proof}

\begin{example}\label{ex:Z25complete}
By Theorem~\ref{thm:exact25}, every permutation of $\Z_{25}$ with more
than seven good slopes is congruence-preserving. By
Lemma~\ref{lem:digits}, its number of good slopes is then a multiple of
$5$. Hence, it has $10$, $15$ or $20$ good slopes, and
Corollary~\ref{cor:fullthreshold} describes it.
Corollary~\ref{cor:counthalf} gives
$62500$ permutations with $20$ good slopes and $5812500$ with $15$. For
$10$ good slopes, three of the four nonzero multipliers occur among the
six digit maps. Corollary~\ref{cor:counthalf} gives
$$
\binom43 5^6\bigl(3^6-3\cdot2^6+3\bigr)=33750000
$$
such permutations $f$, and each of them satisfies $\mu_{25}(f)=2$. In
particular, no permutation of $\Z_{25}$ has $8$, $9$, $11$, $12$, $13$, $14$, $16$,
$17$, $18$ or $19$ good slopes.
\end{example}

The exact bound over $\Z_{25}$ also settles permutations whose
lost congruence levels are isolated.

\begin{corollary}\label{cor:isolatedsharp}
Let $e\ge2$, and let $f$ be a permutation of $\Z_{5^e}$ that is not
congruence-preserving. Suppose that $f$ does not lose congruences at two
consecutive levels. Then, $|D_{5^e}(f)|\le7\cdot5^{e-2}$, and this bound
is attained for every $e\ge2$.
\end{corollary}

\begin{proof}
The hypothesis gives $L(f)=1$.
Proposition~\ref{prop:block} gives a permutation $h$ of $\Z_{25}$ that
is not congruence-preserving and satisfies
$|D_{5^e}(f)|\le5^{e-2}|D_{25}(h)|$. By Theorem~\ref{thm:exact25},
$|D_{25}(h)|\le7$. The bound is attained by the permutations of
Corollary~\ref{cor:twolinesLift}. For $e\ge3$, they have the form
$f(x)=g(\hat x)+(x-\hat x)$ of Proposition~\ref{prop:slopeLift}. If
$x\equiv y\pmod{5^i}$ with $i\ge2$, then $\hat x=\hat y$ and $f(x)-f(y)=x-y$.
Hence, they lose congruences only at level $1$.
\end{proof}

For higher powers of $p$, the moment identities in the proof of
Lemma~\ref{lem:threelines} require control of the sums over the sets
$C_{i,j}$. A regular direction
class with a large count provides this control.

\begin{lemma}\label{lem:cellsums}
Let $e\ge3$, put $P=p^{e-2}$, and let $\Gamma\subseteq\Z_{p^e}^2$ be a
set of $p^e$ points. Suppose that the forms $Y+p^2hX$ with
$h\in\Z_{p^{e-2}}$ are good on $\Gamma$, and that the direction class of
$Y$ has count at least $2P$. For $i\in\F_p$, let $\tau(i)$ be an integer
representative of the element of Lemma~\ref{lem:carry}.
\begin{enumerate}
\item For all $i,j\in\F_p$, the number $|C_{i,j}|$ is a multiple of
$P$, and
$$
\sum_{(X,Y)\in C_{i,j}}\frac{X-\tau(i)}p\in P\Z_{p^{e-1}},\quad
\sum_{(X,Y)\in C_{i,j}}\frac{Y-\tau(j)}p\in P\Z_{p^{e-1}}.
$$
\item For every $j\in\F_p$, there is a set of $p$ points of $\F_p^2$
that is injective in two distinct directions, and in which the first
coordinate $i$ occurs exactly $|C_{i,j}|/P$ times, for every
$i\in\F_p$.
\end{enumerate}
\end{lemma}

\begin{proof}
The form $Y$ is good. Hence, we may write the points of $\Gamma$
uniquely as $(X(y),y)$ with $y\in\Z_{p^e}$. For $j\in\{0,\dots,p-1\}$,
put $\alpha_j(u)=X(j+pu)\bmod p^{e-1}$, where $u\in\Z_{p^{e-1}}$. For
$h\in\Z_{p^{e-1}}$ and $y=j+pu$, we have
$y+phX(y)=j+p\bigl(u+h\alpha_j(u)\bigr)$. Thus, the form $Y+phX$ is good
if and only if $h\in H(\alpha_j)$ for every $j$, as in
Lemma~\ref{lem:cubefibres}. By hypothesis, $p\Z_{p^{e-1}}\subseteq
H(\alpha_j)$, and Lemma~\ref{lem:Hmaps}(i) shows that every $\alpha_j$
is congruence-preserving. By Lemma~\ref{lem:Hmaps}(iii), the set
$\bigcap_jH(\alpha_j)$ is a union of residue classes modulo $p$. Its
size is the count of the direction class of $Y$, which is at least
$2P$. Hence, it contains a unit $b$.

Fix $j$ and $u_0\in\{0,\dots,p-1\}$, and let $i$ be the residue of
$\alpha_j(u_0)$ modulo $p$. Since $\alpha_j$ is congruence-preserving,
there is a congruence-preserving map $\sigma$ of $\Z_P$ with
$\alpha_j(u_0+pv)=i+p\sigma(v)$ for $v\in\Z_P$. The permutation
$u\mapsto u+b\alpha_j(u)$ is congruence-preserving. Hence, it maps the
residue class of $u_0$ modulo $p$ onto a residue class, and the induced
map $v\mapsto v+b\sigma(v)+a$ is a permutation of $\Z_P$ for some
constant $a$. Summing over $v\in\Z_P$, we obtain $Pa+b\sum_v\sigma(v)=0$
in $\Z_P$. Since $b$ is a unit, $\sum_v\sigma(v)=0$ in $\Z_P$.

The $P$ points with second coordinates $j+p(u_0+pv)$, where
$v\in\Z_P$, lie in $C_{i,j}$. Their first coordinates satisfy
$$
\sum_vX\bigl(j+p(u_0+pv)\bigr)\equiv Pi+p\sum_v\sigma(v)\equiv Pi
\pmod{p^{e-1}}.
$$
Since $P\bigl(i-\tau(i)\bigr)$ is divisible by $pP=p^{e-1}$, the sum of
$(X-\tau(i))/p$ over these points lies in $P\Z_{p^{e-1}}$. The sum of
$(Y-\tau(j))/p$ over the same points is
$$
P\,\frac{j-\tau(j)}p+Pu_0+p\sum_{v=0}^{P-1}v,
$$
and it also lies in $P\Z_{p^{e-1}}$, because $P$ is odd. As $u_0$ varies,
these sets of $P$ points partition the sets $C_{i,j}$ with this $j$.
This proves (i).

Let $\psi_j$ be the map of $\F_p$ induced by $\alpha_j$, and put
$S_j=\{(\psi_j(u_0),u_0):u_0\in\F_p\}$. By the preceding paragraph,
the first coordinate $i$ occurs exactly $|C_{i,j}|/P$ times in $S_j$.
The forms $(s,u)\mapsto u$ and $(s,u)\mapsto u+\bar bs$ take every value
of $\F_p$ exactly once on $S_j$, because $u\mapsto u+b\alpha_j(u)$
induces a permutation of $\F_p$. These forms are constant on the lines
of two distinct directions, since $\bar b\neq0$. This proves (ii).
\end{proof}

For higher powers of $5$, Lemmas~\ref{lem:fivepoints},
\ref{lem:threelines} and~\ref{lem:cellsums} restrict the regular
direction classes when the residue pattern is not constant.

\begin{lemma}\label{lem:fivepattern}
Let $e\ge3$, put $P=5^{e-2}$, and let $f$ be a permutation of
$\Z_{5^e}$ that loses congruences at level $1$. Suppose that $w_f$ is
not constant and that five direction classes have a nonzero count.
Then, every regular direction class has count at most $P$.
\end{lemma}

\begin{proof}
By Lemma~\ref{lem:regular}(i), applied to $f$ or to $f^{-1}$, the count
of a regular class is a positive multiple of $P$. Suppose that some
regular class has count at least $2P$. By Lemma~\ref{lem:projective}(i),
exactly one of the six direction classes has count $0$. Choose a
primitive form $X'$ in this class. Choose a form $Y'$ in the regular
class such that all classes congruent to the class of $Y'$ modulo $25$
are good. Their coefficient vectors are linearly independent modulo
$5$. Let $\Gamma'$ be the image of the graph of $f$ in the coordinates
$(X',Y')$. Then, every form $Y'+25hX'$ is good on $\Gamma'$, and the
direction class of $Y'$ has count at least $2P$.

For $i,j\in\F_5$, let $w(i,j)=|C_{i,j}|$ for $\Gamma'$. By
Lemma~\ref{lem:cellsums}(i), every $w(i,j)$ is divisible by $P$. The
reduction of a good form takes every value of $\F_5$ exactly
$5^{e-1}=5P$ times on $\Gamma'$. Thus, every line with one of the five
directions of nonzero count has weight $5P$. By
Lemma~\ref{lem:uniform}(i), the function $w$ is constant on every line
with the remaining direction, which is vertical in the coordinates
$(X',Y')$. Hence, $w(i,j)=Ph(i)$ for nonnegative integers $h(i)$ with
$\sum_ih(i)=5$.

By Lemma~\ref{lem:cellsums}(ii), the numbers $h(i)$ are the
multiplicities of the first coordinates of a set of five points of
$\F_5^2$ that is injective in two distinct directions. If $h(i)=1$ for
every $i$, then $w$ is constant, and so is $w_f$, because the change of
coordinates is invertible modulo $5$. If $h$ is supported at a single
point, then the reduction of the graph of $f$ lies on one line. This
line is not vertical, because $\bar x$ takes every value as $x$ runs
over $\Z_{5^e}$. Hence, $\overline{f(x)}$ is an affine function of
$\bar x$, and $f$ preserves congruences modulo $5$. Both cases are
excluded. By Lemma~\ref{lem:fivepoints}, after an affine change of the
first coordinate, we have
$$
h(0)=3,\ h(1)=h(-1)=1\quad\text{or}\quad h(0)=1,\ h(1)=h(-1)=2,
$$
and $h(2)=h(-2)=0$. As in the proof of Theorem~\ref{thm:exact25}, this
change is induced by a map $X'\mapsto s(X'-x_0)$ with
$s\in\Z_{5^e}^\times$. It keeps the forms $Y'+25hX'$ good, up to
additive constants, and the five directions other than the vertical
one keep a nonzero count. Hence, the new set satisfies conditions (i)
and (ii) of Lemma~\ref{lem:threelines}, with $p=5$ and $k=1$ or $k=2$.
It also satisfies condition (iii), by Lemma~\ref{lem:cellsums}(i). This
contradicts Lemma~\ref{lem:threelines}. Hence, every regular class has
count at most $P$.
\end{proof}

For the powers of $5$, Lemma~\ref{lem:fivepattern} gives the following
recursion.

\begin{lemma}\label{lem:fivepropagate}
Let $e\ge3$, and put $P=5^{e-2}$. Suppose that every permutation $\pi$
of $\Z_{5^{e-1}}$ that is not congruence-preserving satisfies
$N(\pi)\le M$, where $8P/5\le M\le2P$. If $f$ is a permutation of
$\Z_{5^e}$ that loses congruences at every level, then
$$
|D_{5^e}(f)|\le2P+3M,\quad N(f)\le2P+4M.
$$
\end{lemma}

\begin{proof}
As in the proof of Lemma~\ref{lem:strongpropagate}, every direction
class that is not regular has count at most $M$, and every regular class
has count at most $4P$. At most five classes have a nonzero count, and
the vertical class is one of them, because $0\in D_{5^e}(f^{-1})$.

Suppose first that $w_f$ is not constant. By
Lemma~\ref{lem:regularpattern}, applied to $f$ and to $f^{-1}$, every
regular count is at most $2P$. If five classes have a nonzero count,
Lemma~\ref{lem:fivepattern} improves this bound to $P\le M$. Then,
every count is at most $M$, and $|D_{5^e}(f)|\le4M$ and $N(f)\le5M$.
Otherwise, at most four classes have a nonzero count, and at most three
of them are finite. Every count is at most $2P$, because $M\le2P$.
Hence, $|D_{5^e}(f)|\le6P$ and $N(f)\le8P$. Since $8P/5\le M\le2P$, we
have
$$
\max\{4M,6P\}\le2P+3M,\quad \max\{5M,8P\}\le2P+4M.
$$

Suppose now that $w_f$ is constant. By Lemma~\ref{lem:tworegular}, at
most one finite class is regular. If there is no regular finite class,
or if its count is at most $2P$, then
$$
|D_{5^e}(f)|\le\max\{4M,2P+3M\}=2P+3M.
$$
Otherwise, Corollary~\ref{cor:regulartransfer} bounds every other
finite count by $M-4P/5$, and
$$
|D_{5^e}(f)|\le4P+3\Bigl(M-\frac{4P}5\Bigr)=\frac{8P}5+3M\le2P+3M.
$$
For the bound on $N(f)$, suppose first that at most two classes have a
nonzero count. Then, $N(f)\le8P\le2P+4M$. Otherwise, at most one class
is regular, as in the proof of Lemma~\ref{lem:strongpropagate}. If no
class is regular, or if the regular count is at most $2P$, then
$N(f)\le\max\{5M,2P+4M\}=2P+4M$. Otherwise, as in the proof of
Lemma~\ref{lem:strongpropagate}, every other class has count at most
$M-4P/5$. There are at most four such classes, and hence
$$
N(f)\le4P+4\Bigl(M-\frac{4P}5\Bigr)=\frac{4P}5+4M\le2P+4M.
$$
\end{proof}

For higher powers of $5$, we use the structure of permutations with
large values of $N$. The next lemma extends
Lemma~\ref{lem:regularsets} to sets of points. We use the sets
$C_{i,j}$ from Section~\ref{sec:sharpsquares}.

\begin{lemma}\label{lem:regularcells}
Let $k\ge2$, and let $\Gamma\subseteq\Z_{p^k}^2$ be a set of points such
that every set $C_{i,j}$ has exactly $p^{k-2}$ points. Let $L$ be a
primitive linear form such that every class congruent to $[L]$ modulo
$p^2$ is good on $\Gamma$. Then, $L$ is constant modulo $p^2$ on every
set $C_{i,j}$.
\end{lemma}

\begin{proof}
Choose a linear form $X$ whose reduction modulo $p$ is not a multiple of
that of $L$, and use $(X,Y)=(X,L)$ as new coordinates. This change of
coordinates is invertible modulo $p$, and hence it permutes the sets
$C_{i,j}$. The classes congruent to $[Y]$ modulo $p^2$ are the classes
of the forms $Y+p^2hX$, and all of them are good. Since $Y$ is good, we
may write the points of $\Gamma$ uniquely as $(X(y),y)$ with
$y\in\Z_{p^k}$. For $z\in\{0,\dots,p-1\}$ and $u\in\Z_{p^{k-1}}$, put
$\alpha_z(u)=X(z+pu)\bmod p^{k-1}$. As in the proof of
Lemma~\ref{lem:cubefibres}, the form $Y+phX$ is good if and only if
$h\in H(\alpha_z)$ for every $z$. Hence,
$p\Z_{p^{k-1}}\subseteq H(\alpha_z)$, and every $\alpha_z$ is
congruence-preserving by Lemma~\ref{lem:Hmaps}(i). Consequently, the
residue of $X(y)$ modulo $p$ depends only on the residue of $y$ modulo
$p^2$.

The $p^{k-2}$ points with a fixed residue of $y$ modulo $p^2$ therefore
lie in one set $C_{i,j}$. This set also has $p^{k-2}$ points, and hence
the two sets are equal. Thus, $Y=L$ is constant modulo $p^2$ on every
set $C_{i,j}$.
\end{proof}

The next lemma describes the permutations of $\Z_{5^k}$ with a large
value of $N$.

\begin{lemma}\label{lem:fiveprofile}
Let $k\ge3$, put $P=5^{k-2}$, and let $8P/5\le M\le2P$.
Suppose that every permutation $\pi$ of $\Z_{5^{k-1}}$ that is not
congruence-preserving satisfies $N(\pi)\le M$. Let $f$ be a
permutation of $\Z_{5^k}$ that loses congruences at every level.
If $N(f)>5M$, then $w_f$ is constant, exactly five direction classes
have nonzero counts, and exactly one of them is regular. The count of
this regular class is $2P$.
\end{lemma}

\begin{proof}
As in the proof of Lemma~\ref{lem:fivepropagate}, every direction class
that is not regular has count at most $M$. If $w_f$ is not constant,
that proof gives $N(f)\le\max\{5M,8P\}=5M$. Hence, $w_f$ is constant. If
at most two classes have nonzero counts, then $N(f)\le8P\le5M$. If at
least three classes have nonzero counts, then at most one of them is
regular, by the argument in the proof of Lemma~\ref{lem:strongpropagate}.
If none of them is regular, then $N(f)\le5M$.

Suppose that a regular class has count greater than $2P$. Then,
Corollary~\ref{cor:regulartransfer}, applied after the changes of
coordinates in the proof of Lemma~\ref{lem:strongpropagate}, bounds
every other count by $M-4P/5$. Hence,
$$
N(f)\le4P+4\Bigl(M-\frac{4P}5\Bigr)=\frac{4P}5+4M\le5M,
$$
which is excluded. By Lemma~\ref{lem:regular}(i), applied to $f$ or to
$f^{-1}$, the count of a regular class is a positive multiple of $P$. If
it equals $P$, then $N(f)\le P+4M\le5M$. Thus, it equals $2P$. If at
most four classes had nonzero counts, we would obtain
$N(f)\le2P+3M\le5M$, because $M\ge P$. Therefore, exactly five classes
have nonzero counts.
\end{proof}

The structure in Lemma~\ref{lem:fiveprofile} bounds the counts of the
residues other than a regular one.

\begin{lemma}\label{lem:secondtransfer}
Let $e\ge4$, put $k=e-1$, and let $Q\ge8\cdot5^{k-2}$.
Suppose that every permutation $\pi$ of $\Z_{5^k}$ that is not
congruence-preserving and satisfies $N(\pi)>Q$ has constant residue
pattern, exactly five direction classes with nonzero counts, and a
regular direction class. Let $f$ be a permutation of $\Z_{5^e}$ that
loses congruences at every level, and suppose that $w_f$ is constant.
If $c_1$ is regular, then $d_c(f)\le Q$ for every $c\ne c_1$.
\end{lemma}

\begin{proof}
We first consider a map $\gamma$ of $\Z_{5^k}$ that is not
congruence-preserving and satisfies $|H(\gamma)|>Q$, and we put
$S=\{(u,\gamma(u)):u\in\Z_{5^k}\}$. Since
$Q\ge8\cdot5^{k-2}>5^{k-1}-1$, Lemma~\ref{lem:Hnormal} gives a
permutation $\pi$ that is not congruence-preserving, and its proof maps
$S$ onto the graph of $\pi$ by an invertible linear change of
coordinates. Hence, the number of good classes on $S$ is
$N(\pi)\ge|H(\gamma)|>Q$. Invertible linear changes of coordinates
permute the sets $C_{i,j}$ and preserve counts and regularity, as
explained at the beginning of Section~\ref{sec:sharpsquares}. Hence, by
hypothesis, every set $C_{i,j}$ of $S$ has $5^{k-2}$ points, exactly five
direction classes have nonzero counts on $S$, and one of them is
regular.

In the coordinates $(u,t)$ of $S$, a class $[1:h]$ is good exactly when
$h\in H(\gamma)$. These classes form all direction classes except the
one that reduces to $[0:1]$. Suppose that this exceptional direction
class contained a good class. Take a good form from it as the first
coordinate, and a good form from another direction class as the
second. This gives the graph of a permutation $\rho$, which is not
congruence-preserving by Lemma~\ref{lem:projective}(ii), applied
through the graph of $\pi$. The finite good classes of $\rho$ are
exactly the good classes outside the exceptional direction class.
Therefore, $|D_{5^k}(\rho)|=|H(\gamma)|>Q$. This contradicts
Theorem~\ref{thm:allpowers}, which gives
$|D_{5^k}(\rho)|\le8\cdot5^{k-2}\le Q$.

Hence, the exceptional direction class has count zero, and $H(\gamma)$
meets all five residue classes modulo $5$. More generally, every subset
$J\subseteq H(\gamma)$ with $|J|>Q$ meets all five residue classes.
Otherwise, take a good form from the direction class of a residue
missed by $J$ as the first coordinate, and a good form from another
direction class as the second. Every element of $J$ then gives a finite
good class of a permutation $\rho$ that is not congruence-preserving,
and $|J|\le|D_{5^k}(\rho)|\le8\cdot5^{k-2}\le Q$, a contradiction.

Now suppose that $c\ne c_1$ and $d_c(f)>Q$. Choose a good slope $m_1$
with residue $c_1$ such that $m_1+25\Z_{5^e}\subseteq D_{5^e}(f)$, put
$g(x)=f(x)+m_1x$, and let $F=g^{-1}$. By Lemma~\ref{lem:regular}(i), the
maps
$$
\alpha_z(u)=F(z+5u)\bmod5^k
$$
are congruence-preserving. Choose a good slope $m_0$ with residue $c$,
and put $b=m_0-m_1$. Then, $b$ is a unit. The map $L(y)=y+bF(y)$
satisfies $L=(g+b\,\mathrm{id})\circ F$. Both factors are
permutations, because $g+b\,\mathrm{id}=f+m_0\,\mathrm{id}$ and
$m_0\in D_{5^e}(f)$. For $w\in\{0,\dots,4\}$, put
$$
\gamma_w(v)=(g+b\,\mathrm{id})^{-1}(w+5v)\bmod5^k,
\quad J=\bigcap_wH(\gamma_w).
$$
By Lemma~\ref{lem:cubefibres}, $|J|=d_c(f)>Q$. By
Lemma~\ref{lem:tworegular}, the residue $c$ is not regular. Hence, some
$\gamma_w$ is not congruence-preserving, by Lemma~\ref{lem:regular}(ii).
The preceding paragraph, applied to this map and to $J$, shows that $J$
meets every residue class modulo $5$.

We claim that no $\gamma_w$ is congruence-preserving. For every $w$ and
$r\in\F_5$, exactly $5^{e-2}$ values of $v$ satisfy
$\gamma_w(v)\equiv r\pmod5$. Indeed, these values correspond to the
elements $x$ with $\bar x=r$ and $\overline{f(x)}=w-cr$, and $w_f$ is
constant. If $\gamma_w$ were congruence-preserving, its induced map
$\psi_w$ on $\F_5$ would therefore be a permutation, and by
Lemma~\ref{lem:Hmaps}(iii) the residues of the elements of $J$ would lie
in $H(\psi_w)$. However, $H(\psi_w)\ne\F_5$. Indeed, for distinct $a,a'$,
the element $h=-(a-a')/(\psi_w(a)-\psi_w(a'))$ satisfies
$a+h\psi_w(a)=a'+h\psi_w(a')$. This proves the claim.

Thus, the first two paragraphs apply to every $\gamma_w$. Every set
$C_{i,j}$ of $S_w=\{(v,\gamma_w(v)):v\in\Z_{5^k}\}$ has $5^{e-3}$ points,
the exceptional direction class of $S_w$ has count zero, and $S_w$ has a
regular direction class. As every good class on $S_w$ has a
representative $[1:h]$, we may choose $h_w$ such that every class
congruent modulo $25$ to the class of the form $v+h_wt$ is good on
$S_w$. By Lemma~\ref{lem:regularcells}, this form is constant modulo
$25$ on every set $C_{i,j}$ of $S_w$.

Fix a residue class $B$ modulo $125$. For $y\in B$, the residue
$z=y\bmod5$ and the residue of $(y-z)/5$ modulo $25$ are fixed. Since
$\alpha_z$ is congruence-preserving, the residue of $F(y)$ modulo $25$
is fixed as well. Consequently, $L(y)$ is constant modulo $25$ on $B$,
and in particular $w=L(y)\bmod5$ is fixed. The points
$$
\left(\frac{L(y)-w}{5},\ F(y)\bmod5^k\right),\quad y\in B,
$$
belong to $S_w$, because $(g+b\,\mathrm{id})(F(y))=L(y)$. Since $L(y)$
is fixed modulo $25$ and $F(y)$ is fixed modulo $5$, these points lie in
one set $C_{i,j}$ of $S_w$. They are distinct, because $L$ is a
permutation. As both this set of points and
$C_{i,j}$ have $5^{e-3}$ elements, they are equal. It follows that
$$
y+(b+5h_w)F(y)-w=5\left(\frac{L(y)-w}{5}+h_wF(y)\right)
$$
is constant modulo $125$ on $B$. Since $y$ is constant modulo $125$ on
$B$ and $b+5h_w$ is a unit, $F(y)$ is constant modulo $125$ on $B$. This
holds for every $B$. Thus, $F$ preserves congruences modulo $125$, and so
do $g=F^{-1}$ and $f(x)=g(x)-m_1x$. Since $e\ge4$, this contradicts the
hypothesis that $f$ loses congruences at every level.
\end{proof}

Combining these lemmas, we obtain a recursion over two consecutive
exponents.

\begin{lemma}\label{lem:twolevelpropagate}
Let $e\ge4$, put $P=5^{e-2}$, and let $8P/5\le M\le2P$.
Suppose that every permutation $\pi$ of $\Z_{5^{e-1}}$ that is not
congruence-preserving satisfies $N(\pi)\le M$. Let $Q\ge8P/5$
satisfy the hypothesis about permutations with $N(\pi)>Q$ in
Lemma~\ref{lem:secondtransfer}. If $f$ is a permutation of
$\Z_{5^e}$ that loses congruences at every level, then
$$
|D_{5^e}(f)|\le\max\{4M,\ 2P+3Q\},\quad
N(f)\le\max\{5M,\ 2P+4Q\}.
$$
\end{lemma}

\begin{proof}
If $w_f$ is not constant, the proof of
Lemma~\ref{lem:fivepropagate} gives
$|D_{5^e}(f)|\le\max\{4M,6P\}=4M$ and
$N(f)\le\max\{5M,8P\}=5M$.
Suppose that $w_f$ is constant. By Lemma~\ref{lem:tworegular}, at most
one finite direction class is regular. If there is none, then
$|D_{5^e}(f)|\le4M$. If the regular finite class has count at most
$2P$, Lemma~\ref{lem:secondtransfer} bounds the other finite counts by
$Q$, and $|D_{5^e}(f)|\le2P+3Q$. If its count is greater than $2P$, then
Corollary~\ref{cor:regulartransfer} gives
$$
|D_{5^e}(f)|\le4P+3\Bigl(M-\frac{4P}5\Bigr)
=\frac{8P}5+3M\le4M.
$$

For $N(f)$, if at most two classes have nonzero counts, then
$N(f)\le8P\le5M$. If at least three classes have nonzero counts, then
at most one of them is regular. If none of them is regular, then
$N(f)\le5M$. If the regular class has count at most $2P$, the changes
of coordinates used in
Lemma~\ref{lem:strongpropagate} allow
Lemma~\ref{lem:secondtransfer} to be applied to every other class.
Thus, $N(f)\le2P+4Q$. If the regular count is greater than $2P$, the
same changes of coordinates and Corollary~\ref{cor:regulartransfer}
give
$N(f)\le4P+4\Bigl(M-\frac{4P}5\Bigr)
=\frac{4P}5+4M\le5M.$
These inequalities prove both bounds.
\end{proof}

We can now prove the bounds for the powers of $5$.

\begin{corollary}\label{cor:fivebound}
For $e\ge2$, put
$B_e=5^{e-2}\left(10-2\left(\frac45\right)^{\lfloor(e-1)/2\rfloor}\right).$
For $e\ge3$, put
$$
A'_e=5^{e-2}\left(8-\frac65
\left(\frac45\right)^{\lfloor(e-3)/2\rfloor}\right),
\quad A_e=\max\{7\cdot5^{e-2},A'_e\},
$$
and put $A_2=7$. Let $e\ge2$. Every permutation $f$ of $\Z_{5^e}$ that is not
congruence-preserving satisfies
$|D_{5^e}(f)|\le A_e$ and $N(f)\le B_e.$
If $e\ge3$ and $f$ loses congruences at every level, then
$|D_{5^e}(f)|\le A'_e$. In particular, $A'_3=34$ and $A'_4=170$, and the
first values of $A_e$ and $B_e$ are as follows.
$$
\begin{array}{|c|rrrr|}
\hline
e&2&3&4&5\\ \hline
A_e&7&35&175&880\\
B_e&8&42&210&1090\\ \hline
\end{array}
$$
Moreover, a permutation of $\Z_{5^e}$
has more than $A_e$ good slopes if and only if it is
congruence-preserving and all its digit maps are affine with at most
three distinct multipliers.
\end{corollary}

\begin{proof}
We prove the bounds simultaneously by induction on $e$. We also prove
that, for $e\ge3$, every permutation that is not
congruence-preserving and satisfies $N(f)>5B_{e-1}$ has constant
residue pattern, exactly five direction classes with nonzero counts,
and a regular direction class.

The ratios $B_e/5^{e-2}$ do not decrease with $e$, because
$\lfloor(e-1)/2\rfloor$ does not decrease. The same argument, with
$\lfloor(e-3)/2\rfloor$, shows that $A'_e/5^{e-2}$ does not decrease
for $e\ge3$. Thus, the ratios $A_e/5^{e-2}$ do not decrease either.
Moreover, $A_2=A_3/5=A_4/25=7$, and $A_e=A'_e$ for $e\ge5$.

For $e=2$, Theorem~\ref{thm:exact25} gives $A_2=7$ and $B_2=8$.
For $e=3$, suppose first that $f$ preserves a nontrivial congruence.
The reductions of Lemmas~\ref{lem:split} and~\ref{lem:projective}
give $|D_{125}(f)|\le5A_2=35$ and $N(f)\le5B_2=40$.
Otherwise, Lemma~\ref{lem:fivepropagate}, with $P=5$ and $M=8$,
gives $|D_{125}(f)|\le34=A'_3$ and $N(f)\le42=B_3$.
If $N(f)>40$, the permutation loses congruences at every level, and
Lemma~\ref{lem:fiveprofile} proves the additional assertion.

Let $e\ge4$, and assume the assertions for smaller exponents.
Suppose first that $f$ preserves a nontrivial congruence. As in the
proof of Theorem~\ref{thm:allpowers}, some quotient or residue-class
permutation is not congruence-preserving and acts on $\Z_{5^j}$,
where $2\le j<e$. By induction and by the two reduction lemmas,
$$
|D_{5^e}(f)|\le5^{e-j}A_j\le A_e,\quad
N(f)\le5^{e-j}B_j\le5B_{e-1}\le B_e.
$$
The monotonicity of the ratios gives all the last inequalities.
In particular, such a permutation does not satisfy
$N(f)>5B_{e-1}$.

Now suppose that $f$ loses congruences at every level. Put
$P=5^{e-2}$, $M=B_{e-1}$ and $Q=5B_{e-2}$. Since
$8\le B_j/5^{j-2}<10$ for every $j\ge2$, we have
$8P/5\le M<2P$ and $Q\ge8P/5$. The induction hypothesis at exponent
$e-1$ gives the assertion about permutations with $N>Q$ required in
Lemma~\ref{lem:secondtransfer}. Hence,
Lemma~\ref{lem:twolevelpropagate} gives
$$
|D_{5^e}(f)|\le\max\{4M,2P+3Q\},\quad
N(f)\le\max\{5M,2P+4Q\}.
$$
For odd $e$, the formulas give $M=Q<2P$. Therefore, the two maxima
are $2P+3Q=A'_e$ and $2P+4Q=B_e$. For even $e$, the formulas give
$5M=2P+4Q$ and $4M\le2P+3Q$. Again, the maxima are $A'_e$ and $B_e$.
This proves the bounds. If $N(f)>5B_{e-1}=5M$,
Lemma~\ref{lem:fiveprofile} proves the additional assertion.
The induction is complete.

Finally, suppose that $|D_{5^e}(f)|>A_e$. Then, $f$ is
congruence-preserving. By Lemma~\ref{lem:digits},
$|D_{5^e}(f)|=5^{e-1}|T(f)|$. The formulas give
$5^{e-1}<A_e<2\cdot5^{e-1}$. Hence, $|T(f)|\ge2$, and
Theorem~\ref{thm:fields}(i) shows that every digit map is affine.
If $s$ is the number of distinct multipliers of the digit maps, then
$|T(f)|=5-s$, and hence $s\le3$. Conversely, these conditions give
$|D_{5^e}(f)|\ge2\cdot5^{e-1}>A_e$.
\end{proof}

For $e=3$ and $e=4$, the bounds of Corollary~\ref{cor:fivebound} are
attained.

\begin{theorem}\label{thm:exact125}
The largest number of good slopes of a permutation of $\Z_{125}$
that is not congruence-preserving is $35$. Such a permutation has at
most $34$ good slopes if it loses congruences modulo both $5$ and $25$.
\end{theorem}

\begin{proof}
Both upper bounds follow from Corollary~\ref{cor:fivebound} with $e=3$.
Let $g$ be the permutation of Example~\ref{ex:Z25seven}. For
$x\in\Z_{125}$, let $\hat x\in\{0,\dots,24\}$ be the residue of $x$ modulo
$25$, and put $f(x)=g(\hat x)+(x-\hat x)$. The good slopes
$0,1,5,6,10,11,15$ of $g$ are not congruent to $-1$ modulo $5$. Hence,
Proposition~\ref{prop:slopeLift} applies with $c=1$, and it gives
$$
D_{125}(f)=\{a+25t:a\in\{0,1,5,6,10,11,15\},\ 0\le t<5\}.
$$
Thus, $f$ is a permutation with $35$ good slopes, and
$\mu_{125}(f)=\mu_{25}(g)=2$. Since $f(0)=0$ and $f(5)=1$, it does not
preserve congruences modulo $5$.
\end{proof}

\begin{theorem}\label{thm:exact625}
The largest number of good slopes of a permutation of $\Z_{625}$
that is not congruence-preserving is $175$. A permutation that loses
congruences modulo $5$, $25$ and $125$ has at most $170$ good slopes.
\end{theorem}

\begin{proof}
The upper bounds are $A_4=175$ and $A'_4=170$ in
Corollary~\ref{cor:fivebound}. To attain $175$, let $g$ be the
permutation of Example~\ref{ex:Z25seven}. For $x\in\Z_{625}$, let
$\hat x\in\{0,\dots,24\}$ be the residue of $x$ modulo $25$, and put
$f(x)=g(\hat x)+(x-\hat x)$. By Proposition~\ref{prop:slopeLift},
$$
D_{625}(f)=\{a+25t:a\in\{0,1,5,6,10,11,15\},\ 0\le t<25\}.
$$
Thus, $f$ is a permutation with $175$ good slopes. These slopes meet
exactly two residue classes modulo $5$, and
Lemma~\ref{lem:admissible}(iii) gives $\mu_{625}(f)=2$. The values
$f(0)=0$ and $f(5)=1$ show that $f$ does not preserve congruences
modulo $5$.
\end{proof}

Above these exact maxima, Corollary~\ref{cor:fivebound} determines
the possible numbers of good slopes. We record these numbers and
some corresponding permutations for $125$ and $625$.

\begin{example}\label{ex:125complete}
For $n=125$, Theorem~\ref{thm:exact125} shows that every permutation
with more than $35$ good slopes is congruence-preserving. By
Corollary~\ref{cor:fivebound}, it has $50$, $75$ or $100$ good slopes.
The corresponding values of $\mu_{125}$ are $2$, $3$ and $4$. Each such
permutation has $J=1+5+25=31$ digit maps. By
Corollary~\ref{cor:counthalf}, the numbers of these permutations, in the
same order, are
$$
4\cdot5^{31}(3^{31}-3\cdot2^{31}+3),\quad
6\cdot5^{31}(2^{31}-2),\quad 4\cdot5^{31}.
$$
In particular, no permutation has $36,\dots,49$, $51,\dots,74$
or $76,\dots,99$ good slopes.

For a concrete permutation with $50$ good slopes, write
$x=x_0+5x_1+25x_2$, with $0\le x_j<5$, and put
$f(x)=x_0+5\bigl(\lambda(x_0)x_1\bmod5\bigr)+25x_2,$
where the middle digit is represented in $\{0,\dots,4\}$,
$\lambda(0)=2$, $\lambda(1)=3$, and
$\lambda(2)=\lambda(3)=\lambda(4)=1$.
Each digit map is a permutation. Their multipliers form
$\{1,2,3\}$, and hence
$D_{125}(f)=\{m\in\Z_{125}:m\bmod5\in\{0,1\}\}.$
Thus, $|D_{125}(f)|=50$ and $\mu_{125}(f)=2$.
This permutation is not affine. Indeed, $f(0)=0$ and $f(1)=1$
would force an affine map to be the identity, whereas $f(5)=10$.
\end{example}

Table~\ref{tab:bounds} collects the upper bounds of
Theorem~\ref{thm:allpowers} and Corollary~\ref{cor:fivebound}.
The values for $p=3$ and the values for $p=5$, $e=2,3,4$, are attained.

\begin{table}[ht]
\centering
\renewcommand{\arraystretch}{1.2}
\begin{tabular}{|c|r|r|r|r|}
\hline
 & $e=2$ & $e=3$ & $e=4$ & $e=5$\\ \hline
$p=3$ & $3$ & $9$ & $27$ & $81$\\
$p=5$ & $7$ & $35$ & $175$ & $880$\\
$p=7$ & $18$ & $126$ & $882$ & $6174$\\
$p=11$ & $50$ & $550$ & $6050$ & $66550$\\ \hline
\end{tabular}
\caption{Upper bounds for $|D_{p^e}(f)|$ when $f$ is not
congruence-preserving.}\label{tab:bounds}
\end{table}

\begin{example}\label{ex:625complete}
For $n=625$, Theorem~\ref{thm:exact625} and
Corollary~\ref{cor:fivebound} show that every permutation with more
than $175$ good slopes is congruence-preserving. Its digit maps are
affine, and they have one, two or three distinct multipliers.
The corresponding numbers of good slopes are $500$, $375$ and $250$,
and the corresponding values of $\mu_{625}$ are $4$, $3$ and $2$.
Each possibility occurs. For instance, choose every digit map to
have multiplier $1$ for the first possibility. Introduce multiplier
$2$ at one digit map for the second possibility, and multiplier $3$
at a further digit map for the third. Choose all constant terms to
be zero. Lemma~\ref{lem:digits} verifies these assertions directly.
In particular, no permutation of $\Z_{625}$ has $176,\dots,249$,
$251,\dots,374$ or $376,\dots,499$ good slopes.
\end{example}

\section{Gaps below the maximum over composite moduli}\label{sec:gap}

For $n\ge2$, let $q$ be its largest prime divisor and put
$G(n)=\varphi(n)(q-2)/(q-1)$. We will prove that every non-extremal
permutation has at most $G(n)$ good slopes for several families of
odd moduli. We first recall the results of \cite{Extremal} that we
use.

\begin{proposition}[{\cite[Proposition~8.1]{Extremal}}]\label{prop:gapCP}
Let $f$ be a congruence-preserving permutation of $\Z_n$. If
$|D_n(f)|<\varphi(n)$, then $|D_n(f)|\le G(n)$.
\end{proposition}

\begin{theorem}[{\cite[Theorem~8.3]{Extremal}}]\label{thm:gappq}
Let $p<q$ be odd primes. If a permutation $f$ of $\Z_{pq}$ is not
extremal, then
$|D_{pq}(f)|\le(p-1)(q-2)=G(pq).$
\end{theorem}

\begin{theorem}[{\cite[Theorem~8.4]{Extremal}}]\label{thm:gapmq}
Let $n=mq$ be odd, where $m>1$ and $q$ is a prime with $q>m$ and
$\varphi(m)\le(q-1)/2$. Then, every non-extremal permutation $f$ of
$\Z_n$ satisfies $|D_n(f)|\le G(n)$.
\end{theorem}

\begin{theorem}[{\cite[Theorem~8.6]{Extremal}}]\label{thm:gapsplit}
Let $n$ be odd, let $n=dk$ with $d,k\ge2$, and let $f$ be a
permutation of $\Z_n$ that preserves congruences modulo $d$.
Let $g$ and $f_r$ be the permutations in
Lemma~\ref{lem:split-general}. Suppose that each of these permutations,
say $h$ on $\Z_\ell$ with $\ell\in\{d,k\}$, is extremal or satisfies
$|D_\ell(h)|\le G(\ell)$. Then, $f$ is extremal or satisfies
$|D_n(f)|\le G(n)$.
\end{theorem}

\begin{corollary}[{\cite[Corollary~8.7]{Extremal}}]\label{cor:minimalgap}
Let $n$ be odd and composite, and suppose that some non-extremal
permutation of $\Z_n$ has more than $G(n)$ good slopes. Then, there are
a composite divisor $\ell$ of $n$ and a non-extremal permutation $h$ of
$\Z_\ell$ with $|D_\ell(h)|>G(\ell)$ that preserves congruences modulo no
$d$ with $1<d<\ell$.
\end{corollary}

\begin{proposition}[{\cite[Proposition~8.8]{Extremal}}]\label{prop:gapsharp}
Let $n$ be odd and composite, and let $q$ be its largest prime divisor.
Then, there is a non-affine permutation $f$ of $\Z_n$ with
$|D_n(f)|=G(n)$. If $q^2\mid n$, such a permutation can be chosen to be
polynomial.
\end{proposition}

\begin{corollary}\label{cor:gapprime}
Let $p$ be an odd prime and $e\ge2$. Every non-extremal permutation
$f$ of $\Z_{p^e}$ satisfies
$|D_{p^e}(f)|\le G(p^e)=(p-2)p^{e-1}.$
\end{corollary}

\begin{proof}
For congruence-preserving permutations, this is
Proposition~\ref{prop:gapCP}. Otherwise, Theorem~\ref{thm:allpowers}
gives $|D_{p^e}(f)|\le(p-1)p^{e-1}/2$.
This is at most $(p-2)p^{e-1}$, since $p\ge3$.
\end{proof}

\begin{corollary}\label{cor:gapcomp}
Let $n$ be odd and composite, and let $f$ be a permutation of $\Z_n$
that preserves congruences modulo $p^e$ for every $p^e\,\|\,n$.
If $f$ is not extremal, then $|D_n(f)|\le G(n)$.
\end{corollary}

\begin{proof}
For each $p^e\,\|\,n$, let $f_p$ be the induced permutation of
$\Z_{p^e}$. Each output component depends only on its corresponding input
component. Hence, the Chinese remainder theorem gives
$|D_n(f)|=\prod_{p\mid n}|D_{p^e}(f_p)|$.
Every factor is at most $\varphi(p^e)$.
Since $f$ is not extremal, there is a prime $\ell\mid n$, with
$\ell^k\,\|\,n$ say, such that $|D_{\ell^k}(f_\ell)|<\varphi(\ell^k)$.
For $k=1$, this factor is at most $\ell-2$.
For $k\ge2$, Corollary~\ref{cor:gapprime} gives the bound
$(\ell-2)\ell^{k-1}$.
Writing $q$ for the largest prime divisor of $n$, we obtain
$$
|D_n(f)|\le\varphi(n)\frac{\ell-2}{\ell-1}
\le\varphi(n)\frac{q-2}{q-1}=G(n).
$$
\end{proof}

\begin{corollary}\label{cor:threefactors}
Let $n$ be odd and have exactly three prime factors, counted with
multiplicity, and let $f$ be a permutation of $\Z_n$. Suppose that $n$
is a prime power, or that $f$ preserves congruences modulo some divisor
$d$ with $1<d<n$. Then, either $f$ is extremal or $|D_n(f)|\le G(n)$.
\end{corollary}

\begin{proof}
If $n=p^3$, this is Corollary~\ref{cor:gapprime}.
Otherwise, each of $d$ and $n/d$ is a prime, a prime square, or a
product of two distinct odd primes. The gap bound holds for every
permutation of these moduli by Theorem~\ref{thm:universal},
Corollary~\ref{cor:gapprime} and Theorem~\ref{thm:gappq}, respectively.
Hence, Theorem~\ref{thm:gapsplit} applies.
\end{proof}

We now prove that every non-extremal permutation of $\Z_{9q}$ satisfies
the gap bound, where $q\ge5$ is prime. This includes $n=45$, $63$ and
$99$, to which Theorem~\ref{thm:gapmq} does not apply. For these moduli,
Corollary~\ref{cor:threefactors} covers only the permutations that
preserve a nontrivial congruence.
From now on, we write
$\Delta_1=\{1,4,7\}$ and $\Delta_2=\{2,5,8\}$ for the two residue classes
modulo $3$ of the units of $\Z_9$. For a set $P\subseteq\Z_9^2$ of nine
points $(s,y)$, let $\mathcal{A}(P)$ be the set of units $\alpha$ of
$\Z_9$ for which the map $(s,y)\mapsto y+\alpha s$ is injective on $P$.
We use the language of Section~\ref{sec:sharpsquares} for such sets.

\begin{lemma}\label{lem:ninepoints}
Let $P\subseteq\Z_9^2$ be a set of nine points such that
$\mathcal{A}(P)$ meets both $\Delta_1$ and $\Delta_2$.
\begin{enumerate}
\item If $\mathcal{A}(P)\ne\Z_9^\times$, then
$|\mathcal{A}(P)\cap\Delta_1|\le2$ and $|\mathcal{A}(P)\cap\Delta_2|\le2$.
In particular, $\mathcal{A}(P)=\Z_9^\times$ whenever $\mathcal{A}(P)$
contains $\Delta_1$ or $\Delta_2$.
\item If $\mathcal{A}(P)=\Z_9^\times$, then the residues modulo $3$ of
the points of $P$ are the three points of a line of $\F_3^2$ with
direction $\langle(1,0)\rangle$ or $\langle(0,1)\rangle$, and each of
them is the residue of exactly three points of $P$.
\item If $|\mathcal{A}(P)\cap\Delta_i|=2$ for some $i$, then the points of
$P$ have pairwise distinct residues modulo $3$. Moreover, if $\alpha$ is
the element of $\Delta_i\setminus\mathcal{A}(P)$, then the map
$(s,y)\mapsto y+\alpha s$ takes exactly three values on $P$, each of them
at three points.
\end{enumerate}
\end{lemma}

\begin{proof}
Choose $\alpha_1\in\mathcal{A}(P)\cap\Delta_1$ and
$\alpha_2\in\mathcal{A}(P)\cap\Delta_2$. The linear map
$A(s,y)=(y+\alpha_1s,\,y+\alpha_2s)$ has determinant $\alpha_1-\alpha_2$,
which is a unit. Both coordinates of $A$ are injective on $P$. Hence,
$AP$ is the graph of a permutation $\pi$ of $\Z_9$. By the remarks at the
beginning of Section~\ref{sec:sharpsquares}, the class of a form is good
on $P$ if and only if the corresponding class is good on $AP$, and this
correspondence preserves direction classes and their counts. For
$i=1,2$, the classes of the forms $y+\alpha s$ with $\alpha\in\Delta_i$
form one direction class, which has three elements.

Suppose first that $\pi$ is congruence-preserving. By
Lemma~\ref{lem:digits}(ii), $D_9(\pi)$ is a union of residue classes
modulo $3$. The permutation $\pi^{-1}$ is also congruence-preserving,
and Proposition~\ref{prop:fullclass} gives $3\Z_9\subseteq D_9(\pi^{-1})$.
Thus, every class of the vertical direction class is good for $\pi$.
Hence, the good classes of $\pi$ form a union of direction classes, and
so do the good classes on $P$. As $\mathcal{A}(P)$ meets $\Delta_1$ and
$\Delta_2$, we obtain $\mathcal{A}(P)=\Z_9^\times$. Moreover, $\pi$
induces a permutation of $\F_3$, which is affine, and every point of the
graph of this permutation is the residue of exactly three points of
$AP$. Since $A$ is invertible modulo $3$, the residues of the points of
$P$ are the three points of a line $\ell$ of $\F_3^2$, each of them
attained three times. The forms
$y+s$ and $y+2s$ take every value of $\Z_9$ once on $P$. Hence, their
reductions take every value of $\F_3$ three times, and they are not
constant on $\ell$. This excludes the directions $\langle(1,2)\rangle$ and
$\langle(1,1)\rangle$, and proves (ii).

Suppose now that $\pi$ is not congruence-preserving. If $\pi$ is not
spread, then Corollary~\ref{cor:notspread}, applied to $\pi$ and to
$\pi^{-1}$, shows that every direction class has count at most $1$. If
$\pi$ is spread, then Lemma~\ref{lem:rigidcount}, applied to $\pi$ and to
$\pi^{-1}$, shows that every count is at most $2$. In both cases,
$|\mathcal{A}(P)\cap\Delta_i|\le2$ for $i=1,2$, and
$\mathcal{A}(P)\ne\Z_9^\times$. Together with the preceding paragraph,
this proves (i).

Finally, let $|\mathcal{A}(P)\cap\Delta_i|=2$. Then,
$\mathcal{A}(P)\ne\Z_9^\times$, and hence $\pi$ is not
congruence-preserving, as shown above. The direction class formed by the
classes of the forms $y+\alpha s$ with $\alpha\in\Delta_i$ has count $2$
on $P$. Hence, $\pi$ is spread, since otherwise every count would be at
most $1$. Thus, the points of $AP$, and therefore those of $P$, have
pairwise distinct residues modulo $3$. Replacing $\pi$ by $\pi^{-1}$ if
necessary, we may assume
that this direction class corresponds to the direction class of
$\langle(1,-c)\rangle$ for $\pi$, for some $c\in\F_3$. Then, $d_c(\pi)=2$.
By Lemma~\ref{lem:rigidcount}, the residue $c$ is rigid for $\pi$, and
all local maps $\psi_{c,z}$, defined with a fixed lift $\hat c\in\Z_9$ of
$c$, have the same multiplier $\lambda$. Thus,
$\hat c-3\lambda$ is the only slope with residue $c$ that is not good.
As in the proof of Lemma~\ref{lem:liftedplane}, the map
$x\mapsto\pi(x)+(\hat c-3\lambda)x$ is constant on each of the three sets
$B_z$ of Lemma~\ref{lem:fibres}, with distinct values for distinct $z$.
These sets have three elements each. The corresponding form on $P$ is
$y+\alpha s$ with $\alpha\in\Delta_i\setminus\mathcal{A}(P)$, and this
proves (iii).
\end{proof}

The next lemma concerns the maps $L\mapsto L+\kappa\theta(L)$ of $\Z_9$
that occur in the proof of Theorem~\ref{thm:gap9q}.

\begin{lemma}\label{lem:digitshift}
Let $\theta$ be a congruence-preserving map of $\Z_9$, and let
$\kappa\in\Z_9$. Suppose that the map $L\mapsto L+\kappa\theta(L)$ of
$\Z_9$ has a fibre with at least five elements and that all its other
nonempty fibres have one element. Then, the map $L\mapsto L+(\kappa+3)\theta(L)$
does not have this property.
\end{lemma}

\begin{proof}
Write $L=r+3l$ with $r,l\in\{0,1,2\}$, and call $\{r,r+3,r+6\}$ the
block of $r$. Since $\theta$ is congruence-preserving, we can write
$\theta(r+3l)=\theta_0(r)+3\tau_r(l)$ with $\theta_0(r)\in\{0,1,2\}$ and
integers $\tau_r(l)$. Let $\lambda\in\{\kappa,\kappa+3\}$, and write
$r+\lambda\theta_0(r)=\rho(r)+3c_\lambda(r)$ with $\rho(r)\in\{0,1,2\}$
and an integer $c_\lambda(r)$. Then, $\rho(r)$ is the same for both
values of $\lambda$, and $c_{\kappa+3}(r)=c_\kappa(r)+\theta_0(r)$. Since
$3\lambda\tau_r(l)\equiv3\kappa\tau_r(l)\pmod9$, we obtain
$$
L+\lambda\theta(L)\equiv\rho(r)+3\bigl(c_\lambda(r)+l+\kappa\tau_r(l)\bigr)
\pmod9.
$$
Thus, the map $\delta_\lambda(L)=L+\lambda\theta(L)$ sends the block of
$r$ into the residue class of $\rho(r)$ modulo $3$. We call the residue
of $c_\lambda(r)+l+\kappa\tau_r(l)$ modulo $3$ the second digit of
$\delta_\lambda(r+3l)$. Passing from $\kappa$ to $\kappa+3$ adds
$\theta_0(r)$ to all second digits in the block of $r$.

Suppose that $\delta_\lambda$ has a fibre $B$ with at least five elements
while all its other nonempty fibres are singletons. Let $R$ be the set of $r$
whose blocks are mapped into the residue class of the value of $B$. A
residue class modulo $3$ contains three elements of $\Z_9$, and the
elements of these blocks outside $B$ have distinct values different from
the value of $B$. Hence, $|R|\ge2$ and $|B|\ge3|R|-2$, and at least two
blocks of $R$ contain at least two elements of $B$ each. We say that a
block has a repeated digit if two of its elements have the same second
digit. A block has at most one repeated digit, and the elements with this
digit lie in $B$. Therefore, every repeated digit is the second digit of
the value of $B$, and at least two blocks of $R$ have a repeated digit.
Moreover, if $r\ne r'$ lie in $R$, then
$r+\kappa\theta_0(r)\equiv r'+\kappa\theta_0(r')\pmod3$, and hence
$\theta_0(r)\ne\theta_0(r')$. Finally, $R$ is the only set of at least two
blocks that are mapped into one residue class. Hence, $R$ does not
depend on $\lambda$.

Now suppose that both $\delta_\kappa$ and $\delta_{\kappa+3}$ have the
property. The blocks with a repeated digit are the same for both maps,
since passing from $\kappa$ to $\kappa+3$ shifts all second digits of a
block by the same amount. For $\delta_\kappa$, all repeated digits are
equal. For $\delta_{\kappa+3}$, two blocks $r\ne r'$ of $R$ have the same
repeated digit. Their repeated digits for $\delta_\kappa$ are equal, and
they change by $\theta_0(r)$ and $\theta_0(r')$, respectively. Hence,
$\theta_0(r)=\theta_0(r')$, a contradiction.
\end{proof}

We now collect the properties of a hypothetical counterexample over
$\Z_{9q}$.

\begin{lemma}\label{lem:ninecols}
Let $q\ge5$ be a prime, put $n=9q$, and identify $\Z_n$ with
$\Z_9\times\F_q$. Let $f$ be a permutation of $\Z_n$ with
$G(n)<|D_n(f)|<\varphi(n)$. Choose a unit $c$ with
$D_n(f)\subseteq\Z_n^\times-c$, and write
$f(x)-cx=\bigl(g_1(s,t),g_2(s,t)\bigr)$ for $x=(s,t)$. For
$\alpha\in\Z_9^\times$ and $\beta\in\F_q^\times$, put
$$
F_\alpha(s,t)=g_1(s,t)+\alpha s,\quad G_\beta(s,t)=g_2(s,t)+\beta t.
$$
Let $\mathcal{M}$ be the set of pairs $(\alpha,\beta)$ for which
$(F_\alpha,G_\beta)$ is not a permutation of $\Z_9\times\F_q$, and let
$\mathcal{C}_\beta$ be the set of $\alpha$ with
$(\alpha,\beta)\notin\mathcal{M}$. Then, $1\le|\mathcal{M}|\le5$, and the
following hold.
\begin{enumerate}
\item Every $G_\beta$ takes each value exactly nine times. For
$z\in\F_q$, the number $r_z$ of elements $s\in\Z_9$ with $g_2(s,t)=z$
does not depend on $t$.
\item Let $\beta\in\F_q^\times$. For every fibre $\Psi$ of $G_\beta$, the
set $P_\Psi=\{(s,g_1(s,t)):(s,t)\in\Psi\}$ has nine points, and
$\mathcal{C}_\beta$ is the intersection of the sets $\mathcal{A}(P_\Psi)$.
If $\mathcal{C}_\beta$ contains one of $\Delta_1,\Delta_2$ and meets the
other, then $\mathcal{C}_\beta=\Z_9^\times$. In particular, if
$\mathcal{C}_\beta\ne\Z_9^\times$, then $\mathcal{M}$ contains at least two
pairs with second entry $\beta$.
\item Let $(\alpha_0,\beta_0)\in\mathcal{M}$, and let $x\ne x'$ satisfy
$F_{\alpha_0}(x)=F_{\alpha_0}(x')$ and $G_{\beta_0}(x)=G_{\beta_0}(x')$.
Write $x-x'=(\sigma,\tau)$. Then, $\sigma\ne0$ and $\tau\ne0$. If
$\sigma\in3\Z_9$, then $(\alpha,\beta_0)\in\mathcal{M}$ for every
$\alpha\equiv\alpha_0\pmod3$.
\item Either $\mathcal{M}=\Delta_i\times\{\beta_0\}$ for some $i$ and
$\beta_0$, or $q=5$ and
$\mathcal{M}=\{\alpha_1,\alpha_2\}\times\{\beta_1,\beta_2\}$ with
$\alpha_1\in\Delta_1$, $\alpha_2\in\Delta_2$ and $\beta_1\ne\beta_2$.
\end{enumerate}
\end{lemma}

\begin{proof}
For a unit $u=(\alpha,\beta)$, we have
$f(x)-cx+ux=\bigl(F_\alpha(x),G_\beta(x)\bigr)$. Hence,
$(\alpha,\beta)\notin\mathcal{M}$ if and only if $u-c\in D_n(f)$. As
$D_n(f)\subseteq\Z_n^\times-c$, this gives
$|\mathcal{M}|=\varphi(n)-|D_n(f)|$. Since $\varphi(n)=6(q-1)$ and
$G(n)=6(q-2)$, we obtain $1\le|\mathcal{M}|\le5$.

For every $\beta$, the six pairs $(\alpha,\beta)$ are not all in
$\mathcal{M}$. Hence, some $(F_\alpha,G_\beta)$ is a permutation, and
$G_\beta$ takes each value nine times. Lemma~\ref{lem:uniform}(ii),
applied over $\F_q$ with first variable $t$ and index set $\Z_9$, gives
the second assertion of (i).

For (ii), choose $\alpha\in\mathcal{C}_\beta$. Since
$(F_\alpha,G_\beta)$ is a permutation, $F_\alpha$ is injective on every
fibre $\Psi$ of $G_\beta$, and $P_\Psi$ has nine points. For every
$\alpha'$, the map $(F_{\alpha'},G_\beta)$ is a permutation if and only
if $F_{\alpha'}$ is injective on every fibre of $G_\beta$, that is, if
and only if $\alpha'\in\mathcal{A}(P_\Psi)$ for every $\Psi$. If
$\mathcal{C}_\beta$ contains one class and meets the other, then so does
every $\mathcal{A}(P_\Psi)$. Lemma~\ref{lem:ninepoints}(i) gives
$\mathcal{A}(P_\Psi)=\Z_9^\times$ for every $\Psi$, and hence
$\mathcal{C}_\beta=\Z_9^\times$. In particular, $\mathcal{C}_\beta$ cannot
have exactly five elements, which proves the last assertion.

For (iii), let $(\alpha,\beta)\notin\mathcal{M}$. Since
$(F_\alpha,G_\beta)$ is a permutation, the element
$$
(F_\alpha,G_\beta)(x)-(F_\alpha,G_\beta)(x')
=\bigl((\alpha-\alpha_0)\sigma,\,(\beta-\beta_0)\tau\bigr)
$$
is nonzero. If $\sigma=0$, it follows that $\beta\ne\beta_0$ for every pair
outside $\mathcal{M}$. Then, $\mathcal{M}$ contains the six pairs with
second entry $\beta_0$, which is impossible. Suppose that $\tau=0$. If
$\sigma$ is a unit, then $\alpha\ne\alpha_0$ for every pair outside
$\mathcal{M}$. Then, $\mathcal{M}$ contains a pair with second entry
$\beta$ for every $\beta$, and by (ii) it contains at least
$2(q-1)\ge8$ pairs, a contradiction. If $\sigma\in3\Z_9$, then
$\alpha\not\equiv\alpha_0\pmod3$ for every pair outside $\mathcal{M}$.
This gives $|D_n(f)|\le3(q-1)\le G(n)$, again a contradiction. Finally,
let $\sigma\in3\Z_9$ and $\tau\ne0$. For $\alpha\equiv\alpha_0\pmod3$ and
$\beta=\beta_0$, the element above vanishes. Hence, these pairs lie in
$\mathcal{M}$.

For (iv), let $\mathcal{M}_\alpha$ be the set of $\beta$ with
$(\alpha,\beta)\in\mathcal{M}$, for $\alpha\in\Z_9^\times$. Suppose first
that $1\le|\mathcal{M}_\alpha|<(q-1)/2$. Choose $\beta'\notin\mathcal{M}_\alpha$.
As $(F_\alpha,G_{\beta'})$ is a permutation, every fibre $\Phi$ of
$F_\alpha$ has $q$ elements, and $(s,t)\mapsto(t,g_2(s,t))$ is injective
on $\Phi$. Let $S_\Phi\subseteq\F_q^2$ be its image. For
$\beta\notin\mathcal{M}_\alpha$, the form $z+\beta t$ is injective on
$S_\Phi$. Hence, $S_\Phi$ is injective in more than $(q-1)/2$
directions, and it is a line by Lemma~\ref{lem:pointsets}(ii). Now let
$\beta_0\in\mathcal{M}_\alpha$. Then, $G_{\beta_0}$ is not injective on
some fibre $\Phi$. A linear form is either injective or constant on a
line. Hence, $z+\beta_0t$ is constant on the line $S_\Phi$, and $t$ takes
$q$ distinct values on $\Phi$. Since $q>3$, two points
$(s,t)$ and $(s',t')$ of $\Phi$ satisfy $s\equiv s'\pmod3$. By (iii),
$s-s'$ is a nonzero element of $3\Z_9$, and $(\alpha',\beta_0)\in\mathcal{M}$
for every $\alpha'\equiv\alpha\pmod3$. As $\beta_0$ was arbitrary,
$|\mathcal{M}|\ge3|\mathcal{M}_\alpha|$, and therefore
$|\mathcal{M}_\alpha|=1$. Suppose next that
$|\mathcal{M}_\alpha|\ge(q-1)/2$. By (ii), $\mathcal{M}$ contains at least
two pairs with second entry $\beta$ for every $\beta\in\mathcal{M}_\alpha$.
Hence, $|\mathcal{M}|\ge q-1$. Since $|\mathcal{M}|\le5$, we get $q=5$
and $|\mathcal{M}_\alpha|=2$.

Suppose that $\mathcal{M}_\alpha=\{\beta_0\}$ for some $\alpha$, and let
$\Delta_i$ be the class of $\alpha$. By the preceding paragraph,
$\Delta_i\times\{\beta_0\}\subseteq\mathcal{M}$. Suppose that
$\mathcal{M}$ contains a further pair $(\alpha',\beta')$. If
$|\mathcal{M}_{\alpha'}|=1$, then $\Delta_j\times\{\beta'\}\subseteq\mathcal{M}$,
where $\Delta_j$ is the class of $\alpha'$. This set is disjoint from
$\Delta_i\times\{\beta_0\}$, and $|\mathcal{M}|\ge6$, which is impossible.
Hence, $q=5$ and $|\mathcal{M}_{\alpha'}|=2$. If $\alpha'\notin\Delta_i$,
then $\mathcal{M}$ consists of the three pairs in
$\Delta_i\times\{\beta_0\}$ and the two pairs with first entry $\alpha'$.
One of the latter has a second entry different from $\beta_0$, and it is
the only pair of $\mathcal{M}$ with this second entry, contrary to (ii).
If $\alpha'\in\Delta_i$, then $\mathcal{M}_{\alpha'}=\{\beta_0,\beta''\}$
with $\beta''\ne\beta_0$. By (ii), $\mathcal{M}$ contains a further pair
$(\alpha'',\beta'')$, and so $|\mathcal{M}|=5$. If $\alpha''\notin\Delta_i$,
then $\mathcal{M}_{\alpha''}=\{\beta''\}$, which was excluded above.
Hence, $\alpha',\alpha''\in\Delta_i$, and $\mathcal{C}_{\beta''}$
contains the other class and the third element of $\Delta_i$. Now (ii)
gives $\mathcal{C}_{\beta''}=\Z_9^\times$, a contradiction. Thus,
$\mathcal{M}=\Delta_i\times\{\beta_0\}$.

Suppose finally that no $\alpha$ satisfies $|\mathcal{M}_\alpha|=1$. Then,
$q=5$, and $|\mathcal{M}_\alpha|\in\{0,2\}$ for every $\alpha$. As
$1\le|\mathcal{M}|\le5$, there are one or two values of $\alpha$ with
$|\mathcal{M}_\alpha|=2$. By (ii), every $\beta$ that occurs in
$\mathcal{M}$ occurs in at least two pairs of $\mathcal{M}$. Hence, there
are exactly two such values $\alpha_1,\alpha_2$, and
$\mathcal{M}_{\alpha_1}=\mathcal{M}_{\alpha_2}=\{\beta_1,\beta_2\}$. If
$\alpha_1$ and $\alpha_2$ were in the same class, then
$\mathcal{C}_{\beta_1}$ would contain the other class and meet this one,
and (ii) would give $\mathcal{C}_{\beta_1}=\Z_9^\times$. Hence, after
renumbering, $\alpha_1\in\Delta_1$ and $\alpha_2\in\Delta_2$.
\end{proof}

We can now prove the gap bound for the moduli $9q$.

\begin{theorem}\label{thm:gap9q}
Let $q\ge5$ be a prime. Every non-extremal permutation $f$ of $\Z_{9q}$
satisfies $|D_{9q}(f)|\le G(9q)=6(q-2)$.
\end{theorem}

\begin{proof}
Suppose that $G(9q)<|D_{9q}(f)|<\varphi(9q)$, and use the notation of
Lemma~\ref{lem:ninecols} and of its proof. We show that both cases of
Lemma~\ref{lem:ninecols}(iv) are impossible.

Suppose first that $\mathcal{M}=\Delta_i\times\{\beta_0\}$, and let
$\Delta_j$ be the other class. Fix $\alpha\in\Delta_i$. Then,
$\mathcal{M}_\alpha=\{\beta_0\}$, and, as shown in the proof of
Lemma~\ref{lem:ninecols}(iv), every set $S_\Phi$ attached to a fibre
$\Phi$ of $F_\alpha$ is a line. For $\beta\ne\beta_0$, the form
$z+\beta t$ is injective on every $S_\Phi$. Hence, every $S_\Phi$ is
vertical, horizontal, or a line on which $z+\beta_0t$ is constant, and
the last case occurs. For $t\in\F_q$, let $v_t$ be the number of
vertical lines $S_\Phi$ contained in $\{t\}\times\F_q$. For $z\in\F_q$,
let $h_z$ be the number of horizontal lines $S_\Phi$ contained in
$\F_q\times\{z\}$. For $\mu\in\F_q$, let $\ell_\mu$ be the number of
lines $S_\Phi$ on which $z+\beta_0t=\mu$. Counting the elements $s$ with
$g_2(s,t)=z$, we obtain
$$
r_z=v_t+h_z+\ell_{z+\beta_0t}.
$$
Summing over $z$ gives $9=qv_t+\sum_zh_z+\sum_\mu\ell_\mu$. Hence, $v_t$
does not depend on $t$. By Lemma~\ref{lem:ninecols}(i), $\ell_{z+\beta_0t}$
does not depend on $t$ either. As $t$ varies, $z+\beta_0t$ runs over
$\F_q$, and so $\ell_\mu=\ell$ for all $\mu$, where $\ell\ge1$. Now
$9\ge q(v_t+\ell)\ge q$. Therefore, $q\in\{5,7\}$, $\ell=1$ and $v_t=0$.
Thus, for each $\mu\in\F_q$, exactly one fibre of $F_\alpha$ is contained
in the fibre $\Psi_\mu$ of $G_{\beta_0}$ with value $\mu$, and every other
fibre of $F_\alpha$ meets $\Psi_\mu$ in at most one point. Consequently,
$F_\alpha$ takes one value at $q\ge5$ points of $\Psi_\mu$, and each of its
other attained values at one point. This holds for every
$\alpha\in\Delta_i$, since $\alpha$ was arbitrary.

Fix $\mu$ and $b\in\Delta_j$, and put $P=P_{\Psi_\mu}$. Since
$\Delta_j\subseteq\mathcal{C}_{\beta_0}$, the form $y+bs$ is injective on
$P$. Hence, there is a map $\theta$ of $\Z_9$ with $s=\theta(y+bs)$ for
every $(s,y)\in P$. Put $L=y+bs$. For $k\in\{1,2\}$, the form
$y+(b+3k)s=L+3k\theta(L)$ is injective on $P$ as well. Thus,
$3k\in H(\theta)$ for $k=1,2$, and $\theta$ is congruence-preserving by
Lemma~\ref{lem:Hmaps}(i). For $\alpha\in\Delta_i$, the element
$\kappa=\alpha-b$ satisfies $y+\alpha s=L+\kappa\theta(L)$. As $\alpha$
runs over $\Delta_i$, the element $\kappa$ runs over a residue class
modulo $3$. By the preceding paragraph, each map
$L\mapsto L+\kappa\theta(L)$ obtained in this way has a fibre with at least
five elements, and all its other nonempty fibres have one element. This
contradicts Lemma~\ref{lem:digitshift}.

Suppose now that $q=5$ and
$\mathcal{M}=\{\alpha_1,\alpha_2\}\times\{\beta_1,\beta_2\}$ with
$\alpha_1\in\Delta_1$ and $\alpha_2\in\Delta_2$, and let $\beta_3,\beta_4$
be the other two elements of $\F_5^\times$. Put
$\mathcal{C}=\Z_9^\times\setminus\{\alpha_1,\alpha_2\}$. For
$\beta\in\F_5^\times$, we call a fibre $\Psi$ of $G_\beta$ full if
$\mathcal{A}(P_\Psi)=\Z_9^\times$, and lifted otherwise. By
Lemma~\ref{lem:ninepoints}(ii), the residues of the points of
$P_\Psi$ for a full fibre $\Psi$ form a line of $\F_3^2$ with direction
$\langle(1,0)\rangle$ or $\langle(0,1)\rangle$, each point of the line
being attained three times. In the first case, $g_1$ is constant modulo
$3$ on $\Psi$, and we call $\Psi$ horizontal. In the second case, $s$ is
constant modulo $3$ on $\Psi$, and we call $\Psi$ vertical. All fibres of
$G_{\beta_3}$ and $G_{\beta_4}$ are full. Let $\Psi$ be a lifted fibre of
$G_{\beta_1}$ or $G_{\beta_2}$. Then, $\mathcal{A}(P_\Psi)$ contains
$\mathcal{C}$, which has two elements in each class. By
Lemma~\ref{lem:ninepoints}(i), $\mathcal{A}(P_\Psi)=\mathcal{C}$. By
Lemma~\ref{lem:ninepoints}(iii), the points of $P_\Psi$ have pairwise
distinct residues modulo $3$, and each of $F_{\alpha_1}$ and
$F_{\alpha_2}$ takes three values on $\Psi$, each at three points. Since
$\alpha_1\notin\mathcal{C}_{\beta_1}$, the map $G_{\beta_1}$ has a lifted
fibre, and so does $G_{\beta_2}$.

For $\bar s,\bar y\in\F_3$, let $w(\bar s,\bar y)$ be the number of
elements $(s,t)$ with $s\equiv\bar s$ and $g_1(s,t)\equiv\bar y\pmod3$.
Every line of $\F_3^2$ on which $\bar s$ is constant has weight $15$.
Choose $\alpha\in\mathcal{C}\cap\Delta_1$. Since $(F_\alpha,G_{\beta_1})$
is a permutation, $F_\alpha$ takes each value of $\Z_9$ five times. Hence,
$\bar y+\bar s$ takes each value of $\F_3$ fifteen times, and every line
on which $\bar y+\bar s$ is constant has weight $15$. In the same way,
using $\mathcal{C}\cap\Delta_2$, every line on which $\bar y+2\bar s$ is
constant has weight $15$. By Lemma~\ref{lem:uniform}(i), with $r=3$, the
number $w(\bar s,\bar y)=w_{\bar y}$ does not depend on $\bar s$. Now fix
$\beta$, and suppose that $G_\beta$ has $\lambda$ lifted fibres,
$h_{\bar y}$ horizontal fibres with $g_1\equiv\bar y$, and $v_{\bar s}$
vertical fibres with $s\equiv\bar s$. A lifted fibre contains one point
with each pair of residues. A horizontal fibre with $g_1\equiv\bar y$
contains three points with $s\equiv\bar s$ for every $\bar s$, and a
vertical fibre with $s\equiv\bar s$ contains three points with
$g_1\equiv\bar y$ for every $\bar y$. Hence,
$$
w_{\bar y}=\lambda+3h_{\bar y}+3v_{\bar s}.
$$
Thus, $v_{\bar s}=v$ does not depend on $\bar s$, and
$\lambda+\sum_{\bar y}h_{\bar y}+3v=5$. For $\beta=\beta_3$, we have
$\lambda=0$, and so $w_{\bar y}\equiv0\pmod3$. For $\beta=\beta_1$, this
gives $\lambda\equiv0\pmod3$, and $\lambda\ge1$. Therefore, $\lambda=3$,
$v=0$ and $\sum_{\bar y}h_{\bar y}=2$. Hence, $G_{\beta_1}$ has three
lifted fibres and two horizontal fibres, and the same holds for
$G_{\beta_2}$.

Let $\Phi$ be a fibre of $F_{\alpha_1}$. It has five points, and
$(s,t)\mapsto(t,g_2(s,t))$ is injective on $\Phi$, because
$(F_{\alpha_1},G_{\beta_3})$ is a permutation. If $\Psi$ is a lifted fibre
of $G_{\beta_1}$, then $F_{\alpha_1}$ takes three values on $\Psi$, each
three times, and so $\Phi\cap\Psi$ has zero or three points. If $\Psi$ is
a full fibre of $G_{\beta_1}$, then $F_{\alpha_1}$ takes nine distinct
values on $\Psi$, and $\Phi\cap\Psi$ has at most one point. As
$G_{\beta_1}$ has two full fibres and $\Phi$ has five points, $\Phi$ meets
exactly one lifted fibre of $G_{\beta_1}$, in a set $T_\Phi$ of three
points. In the same way,
$\Phi$ meets exactly one lifted fibre of $G_{\beta_2}$, in a set
$T'_\Phi$ of three points. The images of $T_\Phi$ and $T'_\Phi$ under
$(s,t)\mapsto(t,g_2(s,t))$ lie on lines of $\F_5^2$ on which $z+\beta_1t$
and $z+\beta_2t$, respectively, are constant. These lines have distinct
directions, so they meet in at most one point. Since $\Phi$ has five
points, $T_\Phi\cap T'_\Phi$ consists of one point, and
$\Phi=T_\Phi\cup T'_\Phi$.

Let $\Omega_1$ and $\Omega_2$ be the unions of the lifted fibres of
$G_{\beta_1}$ and of $G_{\beta_2}$, respectively, and put
$\Omega=\Omega_1\cap\Omega_2$. Since $\Phi\cap\Omega_1=T_\Phi$ and
$\Phi\cap\Omega_2=T'_\Phi$ for every fibre $\Phi$ of $F_{\alpha_1}$, every
point lies in $\Omega_1\cup\Omega_2$, and $\Omega$ meets every fibre of
$F_{\alpha_1}$ in exactly one point. A lifted fibre $\Psi$ of
$G_{\beta_1}$ is the union of the three sets $T_\Phi$ that it contains.
Hence, it contains exactly three points of $\Omega$. Its other six points
are not in $\Omega_2$, so they lie in the two horizontal fibres of
$G_{\beta_2}$. Let $\bar y_1$ and $\bar y_2$ be the residues of $g_1$ on
these two fibres. The six points have six distinct pairs of residues
$(\bar s,\bar y)$, all with $\bar y\in\{\bar y_1,\bar y_2\}$. Hence,
$\bar y_1\ne\bar y_2$, and the three points of $\Psi\cap\Omega$ have the
residue pairs $(0,\bar y_3)$, $(1,\bar y_3)$ and $(2,\bar y_3)$, where
$\bar y_3$ is the third residue. In the same way, every lifted fibre of
$G_{\beta_2}$ contains exactly three points of $\Omega$. Since
$\Omega\subseteq\Omega_1$, we conclude that $|\Omega|=9$, that
$g_1\equiv\bar y_3$ on $\Omega$, and that $\Omega$ contains exactly three
points with $s\equiv\bar s$ for every $\bar s$. For $\beta=\beta_2$, the
formula above gives $w_{\bar y_3}=3$. Hence, exactly nine points satisfy
$g_1\equiv\bar y_3$, and $\Omega$ is the set of these points.

Let $\beta\in\{\beta_3,\beta_4\}$. All fibres of $G_\beta$ are full.
Thus, $w_{\bar y}=3h_{\bar y}+3v$ and $\sum_{\bar y}h_{\bar y}+3v=5$, where
$v\in\{0,1\}$, and $w_{\bar y_3}=3$ gives $h_{\bar y_3}+v=1$. If $v=0$,
then some horizontal fibre of $G_\beta$ has $g_1\equiv\bar y_3$, and it
is equal to $\Omega$. Hence, $G_\beta$ is constant on $\Omega$. If $v=1$,
then no horizontal fibre has $g_1\equiv\bar y_3$, and $\Omega$ is
contained in the three vertical fibres of $G_\beta$, one for each residue
of $s$. The vertical fibre with $s\equiv\bar s$ contains exactly three
points with $g_1\equiv\bar y_3$, namely the three points of $\Omega$ with
$s\equiv\bar s$. Hence, $G_\beta$ is constant on
$\Omega_{\bar s}=\{(s,t)\in\Omega:s\equiv\bar s\}$ for every $\bar s$.

If $v=1$ for both $\beta_3$ and $\beta_4$, then $g_2+\beta_3t$ and
$g_2+\beta_4t$ are constant on every $\Omega_{\bar s}$. Hence, $t$ and
$g_2$ are constant on $\Omega_{\bar s}$, and so is $G_{\beta_1}$. The three
points of $\Omega_{\bar s}$ then lie in one lifted fibre of $G_{\beta_1}$.
This is impossible, because the points of $\Omega$ in a lifted fibre of
$G_{\beta_1}$ have distinct residues of $s$.

Hence, we may assume that $G_{\beta_3}$ is constant on $\Omega$, with
value $\kappa$. On $\Omega$, we have $g_2=\kappa-\beta_3t$ and
$G_{\beta_1}=\kappa+(\beta_1-\beta_3)t$. The three lifted fibres of
$G_{\beta_1}$ cover $\Omega$, and each of them contains three points of
$\Omega$. Hence, $t$ takes three values on $\Omega$, each at three points.
Let $T_0$ be the set of these values, and write
$\F_5\setminus T_0=\{a_1,a_2\}$. For $t_0\in T_0$, the three points
$(s,t_0)$ of $\Omega$ satisfy $g_2(s,t_0)=\kappa-\beta_3t_0$. Thus,
$r_{\kappa-\beta_3t_0}\ge3$ for the three values $t_0\in T_0$. As
$\sum_zr_z=9$, we get $r_{\kappa-\beta_3a_1}=r_{\kappa-\beta_3a_2}=0$.
Moreover, the values of $G_{\beta_1}$ on its lifted fibres are the
elements $\kappa+(\beta_1-\beta_3)t_0$ with $t_0\in T_0$. Hence, its full
fibres have the values $\kappa+(\beta_1-\beta_3)a_1$ and
$\kappa+(\beta_1-\beta_3)a_2$. In the same way, the full fibres of
$G_{\beta_2}$ have the values $\kappa+(\beta_2-\beta_3)a_1$ and
$\kappa+(\beta_2-\beta_3)a_2$. Let $(t^*,z^*)\in\F_5^2$ satisfy
$$
z^*+\beta_1t^*=\kappa+(\beta_1-\beta_3)a_1,\quad
z^*+\beta_2t^*=\kappa+(\beta_2-\beta_3)a_2.
$$
If $r_{z^*}>0$, then $g_2(s,t^*)=z^*$ for some $s$, by
Lemma~\ref{lem:ninecols}(i). The point $(s,t^*)$ then lies in a full fibre
of $G_{\beta_1}$ and in a full fibre of $G_{\beta_2}$. This contradicts
$\Omega_1\cup\Omega_2=\Z_9\times\F_5$. Hence, $z^*=\kappa-\beta_3a_k$ for
some $k\in\{1,2\}$. If $k=1$, then the first equation gives $t^*=a_1$,
and the second gives $(\beta_2-\beta_3)a_1=(\beta_2-\beta_3)a_2$. If
$k=2$, then the second equation gives $t^*=a_2$, and the first gives
$(\beta_1-\beta_3)a_2=(\beta_1-\beta_3)a_1$. Both conclusions contradict
$a_1\ne a_2$. This completes the proof.
\end{proof}

We isolate two facts needed when the largest prime divisor occurs
with a higher exponent. In both of them, we write the points of
$\Z_{q^e}^2$ and of $\F_q^2$ as $(t,z)$.

\begin{lemma}\label{lem:largeformset}
Let $q\ge5$ be prime, let $e\ge2$, and let
$S\subseteq\Z_{q^e}^2$ have $q^e$ points. Suppose that
$B\subseteq\Z_{q^e}$ has more than $\varphi(q^e)/2$ elements and that
$(t,z)\mapsto z+\beta t$ is injective on $S$ for every $\beta\in B$.
Then, the reduction of $S$ modulo $q$ is a line of $\F_q^2$.
Each point of this line is the residue of exactly $q^{e-1}$ points
of $S$.
\end{lemma}

\begin{proof}
Each residue class modulo $q$ has $q^{e-1}$ elements, and
$|B|>\varphi(q^e)/2\ge2q^{e-1}$ because $q\ge5$. Hence, $B$ meets at
least two residue classes modulo $q$, and there are
$\beta_0,\beta_1\in B$ whose difference is a unit. The coordinates
$u=z+\beta_0t$ and $v=z+\beta_1t$ define an invertible linear change
of coordinates. Both are bijective on $S$, and the image of $S$ is
therefore the graph of a permutation $\pi$ of $\Z_{q^e}$.
The forms $z+\beta t$ with $\beta\in B$ belong to distinct classes,
and these classes are good on $S$. As explained at the beginning of
Section~\ref{sec:sharpsquares}, the change of coordinates maps the
good classes on $S$ bijectively onto the good classes on the graph of
$\pi$. Hence, $N(\pi)\ge|B|>\varphi(q^e)/2$, and
Theorem~\ref{thm:allpowers}(ii) shows that $\pi$ is
congruence-preserving.

Consequently, the graph of $\pi$ reduces modulo $q$ to $q$ distinct
points, each of them attained $q^{e-1}$ times. The same holds for $S$,
because the change of coordinates is invertible modulo $q$. Let
$\overline S$ be this set of $q$ points. Let $\beta\in B$. The form
$z+\beta t$ takes every value of $\Z_{q^e}$ once on $S$. Hence, its
reduction takes every value of $\F_q$ exactly $q^{e-1}$ times on $S$.
As every point of $\overline S$ is the residue of $q^{e-1}$ points of
$S$, the reduction is injective on $\overline S$. Thus, $\overline S$
is injective in the direction $\langle(1,-\bar\beta)\rangle$. The set
$B$ meets more than $(q-1)/2$ residue classes modulo $q$. Therefore,
$\overline S$ is injective in more than $(q-1)/2$ directions, and it
is a line by Lemma~\ref{lem:pointsets}(ii).
\end{proof}

\begin{lemma}\label{lem:horizontalcover}
Let $q$ be prime, and let $\ell_1,\dots,\ell_r$ be lines of $\F_q^2$,
not necessarily distinct, where $r<q$. Suppose that the number of
indices $i$ with $(t,z)\in\ell_i$ depends only on $z$. Then, every
$\ell_i$ is horizontal.
\end{lemma}

\begin{proof}
Suppose that some $\ell_i$ is not horizontal, and choose a point
$(t_0,z_0)$ on it. For $t\in\F_q$, the number of indices $i$ with
$(t,z_0)\in\ell_i$ does not depend on $t$. The lines $\ell_i$ equal to
$\F_q\times\{z_0\}$ contribute the same amount at every $t$, and no
other horizontal line contains $(t,z_0)$. Hence, the number of
non-horizontal lines $\ell_i$ through $(t,z_0)$ does not depend on $t$,
and it is positive at $t=t_0$. Thus, every point of
$\F_q\times\{z_0\}$ lies on a non-horizontal line $\ell_i$. Each such
line meets $\F_q\times\{z_0\}$ in exactly one point. Therefore,
$r\ge q$, a contradiction.
\end{proof}

We can now prove the gap bound for a product $pq^e$ when $q$ is
sufficiently large relative to $p$.

\begin{theorem}\label{thm:gappqe}
Let $p<q$ be odd primes with $q\ge2p-1$, and let $e\ge1$. Every
non-extremal permutation $f$ of $\Z_{pq^e}$ satisfies
$|D_{pq^e}(f)|\le G(pq^e)$.
\end{theorem}

\begin{proof}
We use induction on $e$. For $e=1$, this is Theorem~\ref{thm:gappq}. Let
$e\ge2$, put $n=pq^e$, and suppose that $G(n)<|D_n(f)|<\varphi(n)$.
Identify $\Z_n$ with $\F_p\times\Z_{q^e}$, choose a unit $c=(c_1,c_2)$ with
$D_n(f)\subseteq\Z_n^\times-c$, and write
$f(x)-cx=\bigl(g_1(s,t),g_2(s,t)\bigr)$ for $x=(s,t)$. For
$\alpha\in\F_p^\times$ and $\beta\in\Z_{q^e}^\times$, put
$F_\alpha(s,t)=g_1(s,t)+\alpha s$ and $G_\beta(s,t)=g_2(s,t)+\beta t$, and
let $\mathcal{M}$ be the set of pairs $(\alpha,\beta)$ for which
$(F_\alpha,G_\beta)$ is not a permutation. As in the proof of
Lemma~\ref{lem:ninecols},
$$
|\mathcal{M}|=\varphi(n)-|D_n(f)|<\frac{\varphi(n)}{q-1}=(p-1)q^{e-1}.
$$
For $\alpha\in\F_p^\times$, let $\mathcal{R}_\alpha$ be the set of $\beta$
with $(\alpha,\beta)\notin\mathcal{M}$. Then,
$|\mathcal{R}_\alpha|>\varphi(q^e)-(p-1)q^{e-1}=q^{e-1}(q-p)$, and
$q^{e-1}(q-p)\ge\varphi(q^e)/2$ because $q\ge2p-1$.

Fix $\alpha\in\F_p^\times$. As $\mathcal{R}_\alpha\ne\emptyset$, every fibre
$\Phi$ of $F_\alpha$ has $q^e$ points, and $(s,t)\mapsto(t,g_2(s,t))$ is
injective on $\Phi$. Let $S_\Phi\subseteq\Z_{q^e}^2$ be its image.
For $\beta\in\mathcal{R}_\alpha$, the form $z+\beta t$ is injective
on $S_\Phi$. Since $|\mathcal{R}_\alpha|>\varphi(q^e)/2$,
Lemma~\ref{lem:largeformset} shows that the reduction of $S_\Phi$
modulo $q$ is a line $\ell_\Phi$ of $\F_q^2$, and that each point of
this line is the residue of exactly $q^{e-1}$ points of $S_\Phi$.

For $\bar t,\bar z\in\F_q$, let $W(\bar t,\bar z)$ be the number of
elements $(s,t)$ with $t\equiv\bar t$ and $g_2(s,t)\equiv\bar z\pmod q$.
Every line $\bar t=\mathrm{const}$ has weight $pq^{e-1}$. Let
$\bar\beta\in\F_q^\times$. If no pair $(\alpha',\beta)$ with
$\beta\equiv\bar\beta\pmod q$ were outside $\mathcal{M}$, then $\mathcal{M}$
would have at least $(p-1)q^{e-1}$ elements. Hence, some $G_\beta$ with
$\beta\equiv\bar\beta$ takes every value of $\Z_{q^e}$ exactly $p$ times,
and every line $\bar z+\bar\beta\bar t=\mathrm{const}$ has weight
$pq^{e-1}$. By Lemma~\ref{lem:uniform}(i), $W$ is constant on every line
$\bar z=\mathrm{const}$. On the other hand, every element $(s,t)$ lies
in exactly one fibre $\Phi$ of $F_\alpha$, and $(t,g_2(s,t))\in S_\Phi$.
Hence, $W(\bar t,\bar z)/q^{e-1}$ is the number of fibres $\Phi$ with
$(\bar t,\bar z)\in\ell_\Phi$. This number depends only on $\bar z$.
There are exactly $p<q$ fibres of $F_\alpha$. By
Lemma~\ref{lem:horizontalcover}, every $\ell_\Phi$ is horizontal.
Consequently, $g_2$ is constant modulo $q$ on every fibre of
$F_\alpha$. Let $h_\alpha:\F_p\to\F_q$ be the map with
$g_2\equiv h_\alpha\circ F_\alpha\pmod q$.

For $s,y\in\F_p$, let $N_0(s,y)$ be the number of elements $(s,t)$ with
$g_1(s,t)=y$. Every line $s=\mathrm{const}$ has weight $q^e$, and for
$\alpha\in\F_p^\times$, every line $y+\alpha s=\mathrm{const}$ has weight
$q^e$, since $F_\alpha$ takes every value $q^e$ times. By
Lemma~\ref{lem:uniform}(i), with $r=p$, the number $N_0(s,y)$ depends only
on $y$. Let $Y$ be the set of $y$ with $N_0(s,y)>0$. Suppose first that
$Y=\{y_0\}$. Then, $g_1=y_0$ everywhere, and the first coordinate of $f(x)$
is $y_0+c_1s$. Hence, $f$ preserves congruences modulo $p$. Suppose now that
$y\ne y'$ lie in $Y$. Choose distinct $\alpha,\alpha'\in\F_p^\times$, and
put $\varrho=\alpha'/\alpha$. For every $s\in\F_p$, some element with first
coordinate $s$ satisfies $g_1=y$, and so
$h_\alpha(y+\alpha s)=h_{\alpha'}(y+\alpha's)$. In other words,
$h_\alpha(a)=h_{\alpha'}(y+\varrho(a-y))$ for every $a\in\F_p$. The same
identity holds with $y'$ in place of $y$. If
$y'+\varrho(a'-y')=y+\varrho(a-y)$, then $h_\alpha(a')=h_\alpha(a)$. This
relation means that $a'=a+\delta$ with
$\delta=(y-y')(1-\varrho)/\varrho\ne0$. Hence, $h_\alpha$ is invariant under
the translation by $\delta$, and it is constant. Then, $h_{\alpha'}$ is
constant as well, with the same value. Thus, $g_2$ is constant modulo $q$,
and the second coordinate $g_2(x)+c_2t$ of $f(x)$ depends modulo $q$ only on
the residue of $t$ modulo $q$. Hence, $f$ preserves congruences modulo $q$.

In both cases, $f$ preserves congruences modulo some $d\in\{p,q\}$. The
permutations of Theorem~\ref{thm:gapsplit} act on $\Z_p$ and $\Z_{q^e}$,
or on $\Z_q$ and $\Z_{pq^{e-1}}$. Each of them is extremal or satisfies the
gap bound, by Theorem~\ref{thm:universal}, Corollary~\ref{cor:gapprime}
and the induction hypothesis. By Theorem~\ref{thm:gapsplit}, $f$ is
extremal or satisfies $|D_n(f)|\le G(n)$, a contradiction.
\end{proof}

For instance, Theorem~\ref{thm:gappqe} applies to $n=3q^e$ for every prime
$q\ge5$ and every $e\ge1$. Together with Theorem~\ref{thm:gappq},
Corollary~\ref{cor:gapprime} and Theorem~\ref{thm:gap9q}, it shows that
every non-extremal permutation of $\Z_n$ satisfies the gap bound for every
odd composite $n<105$.

Table~\ref{tab:gap} lists the odd $n\le175$ that have at least two
distinct prime divisors and are not products of two distinct primes. For
each of them, it gives $\varphi(n)$, $G(n)$ and the result that proves the
gap bound. Four entries are not settled by these results.
For each of them, the smaller moduli in
Theorem~\ref{thm:gapsplit} are settled. Hence, a counterexample at one of
these moduli would preserve congruences modulo no divisor $d$ with
$1<d<n$.

\begin{table}[ht]
\centering
\begin{tabular}{|c|c|c|c|}
\hline
$n$ & $\varphi(n)$ & $G(n)$ & Proof of the gap bound\\ \hline
$45=3^2\cdot5$ & $24$ & $18$ & Theorem~\ref{thm:gap9q}\\
$63=3^2\cdot7$ & $36$ & $30$ & Theorem~\ref{thm:gap9q}\\
$75=3\cdot5^2$ & $40$ & $30$ & Theorem~\ref{thm:gappqe}\\
$99=3^2\cdot11$ & $60$ & $54$ & Theorem~\ref{thm:gap9q}\\
$105=3\cdot5\cdot7$ & $48$ & $40$ & not settled here\\
$117=3^2\cdot13$ & $72$ & $66$ & Theorems~\ref{thm:gapmq} and~\ref{thm:gap9q}\\
$135=3^3\cdot5$ & $72$ & $54$ & not settled here\\
$147=3\cdot7^2$ & $84$ & $70$ & Theorem~\ref{thm:gappqe}\\
$153=3^2\cdot17$ & $96$ & $90$ & Theorems~\ref{thm:gapmq} and~\ref{thm:gap9q}\\
$165=3\cdot5\cdot11$ & $80$ & $72$ & not settled here\\
$171=3^2\cdot19$ & $108$ & $102$ & Theorems~\ref{thm:gapmq} and~\ref{thm:gap9q}\\
$175=5^2\cdot7$ & $120$ & $100$ & not settled here\\ \hline
\end{tabular}
\caption{The gap bound for small moduli and the scope of the results proved here.}\label{tab:gap}
\end{table}

\begin{example}\label{ex:gap45sharp}
Let $n=45$, and identify $\Z_{45}$ with $\Z_9\times\F_5$. The proof of
Proposition~\ref{prop:gapsharp}, given in \cite[Proposition~8.8]{Extremal},
uses the permutations $(s,t)\mapsto(s,\lambda(s)t)$, where $\lambda$ is a map
from $\Z_9$ to $\F_5^\times$ with exactly two values. We take
$\lambda(0)=2$ and $\lambda(s)=1$ for $s\neq0$. The resulting
permutation is
$$
f(x)=\begin{cases}2x&\text{if }9\mid x,\\x&\text{otherwise}.\end{cases}
$$
Its good slopes are the pairs $(\alpha,\beta)$ with
$1+\alpha\in\Z_9^\times$ and $\beta\notin\{3,4\}$. Hence,
$D_{45}(f)=\{m\in\Z_{45}:\ m\bmod15\in\{0,1,6,7,10,12\}\}$, $|D_{45}(f)|=18=G(45)$,
and $\mu_{45}(f)=2$. By Theorem~\ref{thm:gap9q}, no permutation of
$\Z_{45}$ has between $19$ and $23$ good slopes. Thus, $18$ is the largest
number of good slopes of a non-extremal permutation of $\Z_{45}$. The
permutation $f$ preserves congruences modulo $9$, but not modulo $5$.
\end{example}

\section{Remaining questions}\label{sec:questions}

The prime-power classification above $p^e/2$ is complete. We end
with the remaining maxima among permutations that are not
congruence-preserving, and with the gap question for odd moduli with at
least two distinct prime divisors.

\begin{question}\label{q:max}
Let $p\ge5$ be prime, and let $e\ge2$. What is the largest number of
good slopes of a permutation of $\Z_{p^e}$ that is not
congruence-preserving?
\end{question}

For $p=3$, the answer is $3^{e-1}$ by Theorem~\ref{thm:allpowers}. For
$p=5$, the answers for $e=2,3,4$ are $7$, $35$ and $175$, by
Theorems~\ref{thm:exact25}, \ref{thm:exact125} and~\ref{thm:exact625}.
For every prime $p\ge5$, Corollary~\ref{cor:twolinesLift} and
Theorem~\ref{thm:allpowers} show that the answer lies between
$(2p-3)p^{e-2}$ and $(p-1)^2p^{e-2}/2$. For $p=5$ and $e\ge5$,
Corollary~\ref{cor:fivebound} improves the upper bound, and we obtain
$$
7\cdot5^{e-2}\le\max_f|D_{5^e}(f)|\le
5^{e-2}\left(8-\frac65\left(\frac45\right)^{\lfloor(e-3)/2\rfloor}\right),
$$
where $f$ runs over the permutations that are not congruence-preserving.
For $p=5$, the first exponent not settled here is $e=5$.
The maximum over $\Z_{3125}$ lies between $875$ and $880$.
Any permutation of $\Z_{3125}$ that is not congruence-preserving
and has more than $875$ good slopes loses congruences at every level.
Indeed, by Lemma~\ref{lem:split} and
Corollary~\ref{cor:fivebound}, a preserved congruence gives the upper
bound $5^{5-j}A_j=875$ for some $j$ with $2\le j\le4$.

For $p=5$, Corollary~\ref{cor:isolatedsharp} also gives the exact
maximum $7\cdot5^{e-2}$ among permutations that are not
congruence-preserving and whose lost congruence levels are isolated.
These statements do not determine the exact
maximum in the remaining cases.

\begin{question}\label{q:gap}
Let $n$ be an odd integer with at least two distinct prime divisors,
and let $q$ be its largest prime divisor. Does every
non-extremal permutation $f$ of $\Z_n$ satisfy
$|D_n(f)|\le\varphi(n)\,\frac{q-2}{q-1}?$
\end{question}

This question was raised in \cite[Question~11.2]{Extremal}. By the results of
Section~\ref{sec:gap}, the answer is affirmative for
congruence-preserving permutations, for products of two distinct odd
primes, and for the moduli of Theorems~\ref{thm:gapmq}, \ref{thm:gap9q}
and~\ref{thm:gappqe}. Proposition~\ref{prop:gapsharp} shows that the
bound cannot be improved. The smallest modulus not covered by
these results is $105$. Table~\ref{tab:gap} records the other
unsettled entries up to $175$. By Corollary~\ref{cor:minimalgap}, a
counterexample at the smallest possible modulus would preserve no
nontrivial congruence.

The restriction to odd $n$ is necessary. For $n=12$, we have
$G(12)=2$. However, the permutation $f$ of $\Z_{12}$ with
$\bigl(f(0),f(1),\dots,f(11)\bigr)=(0,1,2,3,5,4,7,9,6,11,8,10)$
has exactly the good slopes $0$, $2$ and $8$. It is not extremal,
since $\varphi(12)=4$.

\section*{Acknowledgements}
The author is grateful to Simeon Ball for suggesting
the problems that motivated this work.

{\small
\bibliographystyle{abbrv}
\bibliography{ref}
}

\end{document}